\documentclass[11pt]{article}    
\usepackage{latexsym}
\usepackage{graphicx}
\usepackage{url}
\usepackage{amsmath}    
\usepackage{fixmath}
\usepackage{amssymb}
\usepackage{psfrag}
\usepackage{fullpage}
\usepackage{subfigure}
\usepackage{enumerate}
\usepackage{bbm}
\usepackage{color}
\usepackage[table]{xcolor}
\usepackage{booktabs}

\usepackage{multirow}
\usepackage{boxedminipage}
\usepackage{subfig}
\usepackage{algorithmic,algorithm}
\usepackage[american]{babel}
\usepackage{setspace}
\usepackage[authoryear,round,comma,compress]{natbib}
\bibpunct{(}{)}{,}{a}{}{,}

\long\def\symbolfootnote[#1]#2{\begingroup%
\def\thefootnote{\fnsymbol{footnote}}\footnote[#1]{#2}\endgroup}

\DeclareMathOperator*{\argmin}{arg\,min}

\newcommand {\beq}{\begin{equation}}
\newcommand {\eeq}{\end{equation}}
\newcommand {\beqn}{\begin{equation*}}
\newcommand {\eeqn}{\end{equation*}}
\newcommand {\bear}{\begin{eqnarray}}
\newcommand {\eear}{\end{eqnarray}}
\newcommand {\bearn}{\begin{eqnarray*}}
\newcommand {\eearn}{\end{eqnarray*}}
\newcommand{\Expect}{\mathbb{E}}

\newtheorem{theorem}{Theorem}

\newtheorem{lemma}{Lemma}
\newtheorem{prop}{Proposition}

\newtheorem{assumption}{Assumption}
\newtheorem{example}{Example}

\def\qed{ \ \vrule width.2cm height.2cm depth0cm\iffalse \smallskip\fi }

\def\spacen{1.3}  
\def\spaces{1.05}  

\normalsize

\begin{document}

\title{
\vspace{-0.7cm}
Non-stationary Stochastic Optimization}

\author{
{\sf Omar Besbes}
\\Columbia University\and
{\sf Yonatan Gur}
\\Stanford University\and
{\sf Assaf Zeevi}\thanks{This work is supported by NSF grant 0964170 and BSF grant 2010466. Correspondence: {\tt ob2105@columbia.edu}, {\tt ygur@stanford.edu}, {\tt assaf@gsb.columbia.edu}.}
\\Columbia University
}
\vspace{0.1cm}
\date{first version: July 2013, current version: \today}

\maketitle

\vspace{-0.1cm}
\begin{abstract}

We consider a non-stationary variant of a sequential stochastic optimization problem, in which the underlying cost functions may change along the horizon. We propose a measure, termed {\it variation budget}, that controls the extent of said change, and study how restrictions on this budget impact achievable performance. We identify sharp conditions under which it is possible to achieve long-run-average optimality and more refined performance measures such as rate optimality that fully characterize the complexity of such problems. In doing so, we also establish a strong connection between two rather disparate strands of literature: adversarial online convex optimization; and the more traditional stochastic approximation paradigm (couched in a non-stationary setting). This connection is the key to deriving well performing policies in the latter, by leveraging  structure of optimal policies in the former. Finally, tight bounds on the minimax regret allow us to quantify the ``price of non-stationarity," which mathematically captures the added complexity embedded in a temporally changing environment versus a stationary one. \\\vspace{-0.2cm}

\end{abstract}

\setstretch{1.43}
\vspace{-0.3cm}
\section{Introduction and Overview}\label{sec:intro}
\vspace{-0.2cm}
\textbf{\quad Background and motivation.} In the prototypical setting of sequential stochastic optimization, a decision maker selects at each epoch $t\in \left\{1,\ldots T\right\}$ a point $X_{t}$ that belongs (typically) to  some convex compact action set $\mathcal{X} \subset \mathbb{R}^{d}$, and incurs a {\it cost}  $f(X_t)$, where  $f(\cdot)$ is an a-priori unknown convex {\it cost function}. Subsequent to that, a {\it feedback}  $\phi_{t}\left(X_{t},f\right)$ is given to the decision maker;   representative feedback structures include a noisy realization of the cost and/or the gradient of the cost.  When the cost function  is assumed to be strongly convex, a typical objective is to minimize the mean-squared-error, $\mathbb{E}\left\|X_{T}-x^{\ast}\right\|^{2}$, where $x^{\ast}$ denotes the minimizer of $f(\cdot)$ in $\mathcal{X}$. When $f(\cdot)$ is only assumed  to be {\it weakly convex}, a more reasonable objective is to minimize $\mathbb{E}\left[f\left(X_{T}\right)-f\left(x^{\ast}\right)\right]$, the expected difference between the cost incurred at the terminal epoch $T$ and the minimal achievable cost. (This objective reduces to the  MSE criterion, up to a multiplicative constant, in the strongly convex case.)  The study of such problems originates with the pioneering work of \cite{Rob-Mon1951} which focuses on  stochastic estimation of a level crossing, and its counterpart studied by \cite{Kie-Vol1952} which focuses on stochastic estimation of the point of maximum; these methods are collectively known as stochastic approximation (SA), and with some abuse of terminology we will use this term to refer to both the methods as well as the problem area. Since the publication of these seminal papers, SA has been widely studied and applied to diverse problems in a variety of fields including Economics, Statistics, Operation Research, Engineering and Computer Science;  cf. books by \cite{Ben-Pri-Met1990} and \cite{Kus-Yin2003}, and a survey by \cite{Lai2003}.

A fundamental assumption in SA which has been adopted by almost all of the relevant literature (exceptions to be noted in what follows),  is that the cost function does not change throughout the horizon over which we seek to (sequentially) optimize it. Departure from this stationarity assumption brings forward many fundamental questions. Primarily, how to model temporal changes in a manner
that is ``rich" enough to capture a broad set of scenarios while still being mathematically tractable, and what is the performance that can be achieved in such settings in comparison to the stationary SA environment. Our paper is concerned with these questions.

\textbf{The non-stationary SA problem.}
Consider the stationary SA formulation outlined above with the following modifications: rather than a single unknown cost function, there is now a
\emph{sequence} of convex functions $\{f_t:t=1,\ldots,T\}$;  like the stationary setting, in every epoch $t=1,\ldots,T$ the  decision maker selects a point $X_{t}\in \mathcal{X}$ (this will be referred to as ``action" or ``decision" in what follows), and then observes a feedback, only now this signal, $\phi_{t}\left(X_{t},f_t\right)$, will depend on the particular function within the sequence. In this paper we consider two canonical feedback structures alluded to earlier, namely, noisy access to the function value $f(X_t)$, and noisy access to the gradient $\nabla f(X_t)$. Let $\{x^{\ast}_{t}:t=1,\ldots,T\}$ denote the sequence of minimizers corresponding to the sequence of cost functions.

In this ``moving target" formulation, a natural objective is to minimize the {\it cumulative} counterpart of the performance measure used in the stationary setting, for example, $\sum_{t=1}^{T}\mathbb{E}\left[f_{t}\left(X_{t}\right)-f_{t}\left(x^{\ast}_{t}\right)\right]$ in the general convex case. This is often referred to in the literature as the {\it regret}. It measures the quality of a policy, and the sequence of actions $\{X_1,\ldots,X_T\}$ it generates, by comparing its performance to a clairvoyant that knows the sequence of functions in advance, and hence selects the minimizer $x^{\ast}_t$ at each step $t$; we refer to this benchmark as a {\it dynamic oracle} for reasons that will become clear soon.\footnote{A more precise definition of an admissible policy will be advanced in the next section, but roughly speaking, we restrict attention to policies that are non-anticipating and adapted to past actions and observed feedback signals, allowing for auxiliary randomization; hence the expectation above is taken with respect to any randomness in the feedback, as well as in the policy's actions.}

To constrain temporal changes in the sequence of functions, this paper introduces the concept of a {\it temporal uncertainty set} $\mathcal{V}$, which is driven by a {\it variation budget} $V_T$:\vspace{-0.3cm}
\[
\mathcal{V}:= \left\{\left\{f_{1},\ldots,f_{T}\right\} \;:\;  \mbox{Var}(f_1,\ldots,f_T)\leq V_{T} \right\}.\vspace{-0.3cm}
\]
The precise definition of the variation functional Var$(\cdot)$ will be given in \S2; roughly speaking, it measures the extent to which functions can change from one time step to the next, and adds this up over the horizon $T$. As will be seen in \S2, the notion of variation we propose allows for a broad range of temporal changes in the sequence of functions and minimizers. Note that the variation budget is allowed to depend on the length of the horizon, and therefore measures the scales of variation relative to the latter.

For the purpose of outlining the flavor of our main analytical findings and key insights, let us further formalize the notion of {\it regret} of a policy $\pi$ relative to the above mentioned dynamic oracle:\vspace{-0.2cm}
\[
\mathcal{R}^{\pi}_{\phi}(\mathcal{V},T) = \sup_{f\in \mathcal{V}}\left\{\mathbb{E}^{\pi}\left[\sum_{t=1}^{T}f_{t}(X_{t})\right] -\sum_{t=1}^{T}f_{t}(x_{t}^{\ast})\right\}.\vspace{-0.2cm}
\]
In this set up, a policy $\pi$ is chosen and then nature (playing the role of the adversary) selects the sequence of functions $f:=\{f_t\}_{t=1,\ldots,T} \in \mathcal{V}$ that maximizes the regret; here we have made explicit the dependence of the regret and the expectation operator on the policy $\pi$, as well as its dependence on the feedback mechanism $\phi$ which governs the observations. The first order characteristic of a ``good'' policy is that it achieves \emph{sublinear} regret, namely,\vspace{-0.2cm}
\[
\frac{\mathcal{R}^{\pi}_{\phi}\left(\mathcal{V},T\right)}{T} \to 0 \quad \mbox{as }\; T \to \infty.\vspace{-0.2cm}
\]
A policy $\pi$ with the above characteristic is called {\it long-run-average optimal}, as the average cost it incurs (per period) asymptotically approaches the one incurred by the clairvoyant benchmark. Differentiating among such policies requires a more refined yardstick. Let $\mathcal{R}^{\ast}_{\phi}\left(\mathcal{V},T\right)$ denote the {\it minimax regret}: the minimal regret that can be achieved over the space of admissible policies subject to feedback signal~$\phi$, \emph{uniformly} over nature's choice of cost function sequences within the temporal uncertainty set $\mathcal{V}$. A policy is said to be \emph{rate optimal} if it achieves the minimax regret up to a constant multiplicative factor; this implies that, in terms of growth rate of regret, the policy's performance is essentially best possible.

\textbf{Overview of the main contributions.} Our main results and key qualitative insights can be summarized as follows:

{\it 1. Necessary and sufficient conditions for sublinear regret.} We first show that if the variation budget $V_T$ is {\it linear} in $T$, then, as one may expect, sublinear regret {\it cannot} be achieved  by any admissible policy. Conversely, we show that if $V_T$ 
is {\it sublinear} in $T$, long-run-average optimal policies exist. So, our notion of temporal uncertainty supports a sharp dichotomy in characterizing first-order optimality in the non-stationary SA problem.

{\it 2. Complexity characterization.}  We prove a sequence of results that characterizes the order of the minimax regret for both the convex as well as the strongly convex settings. This is done by deriving lower bounds on the regret that hold for {\it any} admissible policy, and then proving that the order of these lower bounds can be achieved by suitable (rate optimal) policies.   The essence of these results can be summarized by the following characterization of the minimax regret:\vspace{-0.4cm}
\[
 \mathcal{R}^{*}_{\phi}(\mathcal{V},T) \asymp V_T^\alpha T^{1-\alpha},\vspace{-0.4cm}
\]
where $\alpha$ is either $1/3$ or $1/2$  depending  on the particulars of the problem (namely, whether the cost functions in $\mathcal{V}$ are convex/strongly convex, and whether the feedback $\phi$ is a noisy observation of the cost/gradient); see below for more specificity, and further details in \S4 and \S5.

{\it 3. The ``price of non-stationarity.''} The minimax regret characterization allows, among other things, to contrast the stationary and non-stationary environments, where  the ``price'' of the latter relative to the former is expressed in terms of the ``radius'' (variation budget) of the temporal uncertainty set. The table below summarizes our main findings.
\begin{table}[h]
  \centering
    \begin{tabular}{|c||c||c||c|c}
    \hline
\multicolumn{2}{|c||}{Setting} & \multicolumn{2}{|c|}{Order of regret}\\
\hline
\multicolumn{1}{|c||}{Class of functions} & \multicolumn{1}{|c||}{Feedback} & \multicolumn{1}{|c||}{Stationary} & \multicolumn{1}{|c|}{Non-stationary}\\
\hline
convex  & noisy gradient  &  $\sqrt{T}$     & $V_T^{1/3} T^{2/3}$ \\
strongly convex & noisy gradient &  $\log{T}$     & $\sqrt{V_T T}$ \\
strongly convex & noisy function &   $\sqrt{T}$     & $V_T^{1/3} T^{2/3}$ \\
\hline
    \end{tabular}
    \caption{\small \textbf{The price of non-stationarity.} The rate of growth of the minimax regret in the stationary and non-stationary settings under different assumptions on the cost functions and feedback signal.}\vspace{-0.3cm}
  \label{tab:summary}
\end{table}
Note  that even in the most ``forgiving'' non-stationary environment, where the variation budget $V_T$ is a constant and independent of $T$, there is a marked degradation in performance between the stationary and non-stationary settings. (The table omits the general convex case with noisy cost observations; this will be explained later in the paper.)

{\it 4. A meta principle for constructing optimal policies.} One of the key insights we wish to communicate in this paper pertains to the construction of well performing policies, either long-run-average, or rate optimal. The main idea is a result of bridging two relatively disconnected streams of literature that deal with dynamic optimization under uncertainty from very different perspectives: the so-called {\it adversarial} and the {\it stochastic}  frameworks.
The former, which in our context is often referred to as online convex optimization (OCO), allows nature to select the worst possible function at {\it each} point in time depending on the actions of the decision maker, and with little constraints on nature's choices. \label{pg:contibut-OCO-v-SA}
This constitutes a more pessimistic environment compared with the traditional stochastic setting where the function is picked a priori at $t=0$ and held fixed thereafter, or the setting we propose here, where the {\it sequence} of functions is chosen by nature subject to a variation constraint.
Because of the  freedom awarded to nature in OCO settings, a policy's performance is typically measured relative to a rather coarse benchmark, known as the {\it single best action in hindsight}; the best static action that would have been picked ex post, namely, after having observed all of nature's choices of functions.
While typically a policy that is designed to compete with the single best action benchmark in an adversarial OCO setting does not admit performance guarantees in our non-stationary stochastic problem setting (relative to a dynamic oracle), we establish an important connection between performance in the former and the latter environments, given roughly by the following ``meta principle":
\begin{quote}
If a policy has ``good" performance with respect to the  single best action in the adversarial framework, it can be adapted in a manner that guarantees ``good" performance in the stochastic non-stationary environment subject to the variation budget constraint.
\end{quote}
In particular, according to this principle, a policy with sublinear regret in an OCO setting can be adapted to achieve sublinear regret in the non-stationary stochastic setting, and in a similar manner we can port over the property of rate-optimality. It is important to emphasize that while policies that admit these properties have, by and large, been identified in the OCO literature\footnote{For the sake of completeness, to establish the connection between the adversarial and the stochastic literature streams, we adapt, where needed,  results in the former setting to the case of noisy feedback.}, to the best of our knowledge there are no counterparts to date in a non-stationary stochastic setting, including the one considered in this paper. (It is worthwhile noting that the construction of said policies is mostly done with the intent of providing a relatively simple and unified way to highlight key tradeoffs at play.)\label{pg:theory-v-practice}

\textbf{Relation to literature.} The use of the cumulative performance criterion and regret, while mostly absent from the traditional SA stream of literature, has been adapted in several occasions. Examples include the work of \cite{Cope2009}, which is couched in an environment where the feedback structure is noisy observations of the cost and the target function is strongly convex. That paper shows that the  estimation scheme of \cite{Kie-Vol1952} is rate optimal and the minimax regret in such a setting is of order $\sqrt{T}$. Considering a convex (and differentiable) cost function, \cite{Aga-Fos-Hsu-Kak-Rak2011} showed that the minimax regret is of the same order, building on estimation methods presented in \cite{Nem-Yud1983}. In the context of gradient-type feedback and strongly convex cost, it is straightforward to verify that the scheme of \cite{Rob-Mon1951} is rate optimal, and the minimax regret is of order $\log T$.

While temporal changes in the cost function are typically not discussed within the traditional stationary SA literature (see chapter 3 in \cite{Kus-Yin2003}, and chapter 4 in \cite{Ben-Pri-Met1990} for exceptions), the literature on OCO, which has mostly evolved in the machine learning community starting with \cite{Zin2003},
allows the cost function to be selected at any point in time by an \emph{adversary}.
As discussed above, the performance of a policy in this setting is compared against a relatively weak benchmark, namely, the single best action in hindsight; or, a {\it static} oracle. These ideas have their origin in game theory with the work of \cite{Bla1956} and \cite{Han1957}, and have since seen significant development in several sequential decision making settings; cf. \cite{Cesa-Bia-Lugosi} for an overview. The OCO literature largely focuses on a class of either convex or strongly convex cost functions, and sub-linearity and rate optimality of policies have been studied for a variety of feedback structures. The original work of  \cite{Zin2003} considered the class of convex functions, and focused on a feedback structure in which the function $f_{t}$ is \emph{entirely revealed} after the selection of $X_{t}$, providing an \emph{online gradient descent} algorithm with regret of order $\sqrt{T}$; see also \cite{FLA2005}. \cite{Haz-Aga-Kal2007} achieve regret of order $\log T$ for a class of strongly convex cost functions, when the gradient of $f_{t}$, evaluated at $X_{t}$ is observed. Additional algorithms were shown to be rate optimal under further assumptions on the function class (see, e.g., \citealt{Kal-Vem2003}, \citealt{Haz-Aga-Kal2007}), or other feedback structures such as multi-point access (\citealt{Aga-Dek-Xia2010}).
A closer paper, at least in spirit, is that of \cite{HazKal2010}. It derives upper bounds on the regret with respect to the static single best action, in terms of a measure of dispersion of the cost functions chosen by nature, akin to variance. The cost functions in their setting are restricted to be linear and are revealed to the decision maker after each action.


\label{pg:lit-adversarial} It is important to draw attention to a significant distinction between the framework we pursue in this paper and the adversarial setting, concerning the quality of the benchmark that is used in each of the two formulations. Recall, in the adversarial setting the performance of a policy is compared to the ex post best static feasible solution, while in our setting  the benchmark is given by a dynamic oracle (where ``dynamic'' refers to the sequence of minima $\{f_{t}(x_t^*)\}$ and minimizers $\{x_t^*\}$ that is changing throughout the time horizon). It is fairly straightforward that the gap between the performance of the static oracle that uses the single best action, and that of the dynamic oracle can be significant, in particular, these quantities may differ by order $T$; for an illustrative example see \S2, Example 1.
Therefore, even if it is possible to show that a policy has a ``small'' regret relative to the best static action, there is no guarantee on how well such a policy will perform when measured against the best dynamic sequence of decisions. A second potential limitation of the adversarial framework lies in its rather pessimistic assumption of the world in which policies are to operate in, to wit, the environment can change at any point in time in the worst possible way as a {\it reaction} to the policy's chosen actions. In most application domains, one can argue, the operating environment is not nearly as harsh.

Key to establishing the connection between the adversarial setting and the non-stationary stochastic framework proposed herein is the notion of a variation budget, and the corresponding temporal uncertainty set, that curtails nature's actions in our 
formulation.
These ideas echo, at least philosophically, concepts that have permeated the robust optimization literature, where uncertainty sets are fundamental predicates; see, e.g., \cite{BenTal-Nem1997}, and a survey 
by \cite{BerBroCar-2007}.

\label{pg:lit-kalman} Another line of research considers sequential 
stochastic optimization using Kalman filters (\citealt{Kal1960}). There, the typical objective is to minimize the mean square error when estimating a state, under zero-mean Gaussian noise. In the non-stationary variant of this problem the state may change; in such cases the aforementioned change is typically well-structured by some parameterized dynamics. An overview of this research domain is given in \cite{Hay2001}, where Chapters 3,4, and 6 include a survey of methods and applications for state-dynamic models. The focus, formulation, and analysis in this paper are different from the ones adopted in the literature on Kalman filters in the following key aspects. First, 
a main interest of the current study is in characterizing the extent of non-stationarity 
under which one may achieve sublinear regret with respect to the dynamic oracle benchmark. In particular, we show 
that whenever the variation is a sublinear function of the time horizon $T$, one may achieve sublinear regret relative to the dynamic oracle, but when variation is at least linear in $T$ sublinear regret is not achievable. While non-stationary instances that are considered in the literature on Kalman filters typically fall under the latter case (linear variation), the focus of the current paper is on characterizing the minimax regret in the former. Second, the formulation in this paper is more general than the one adopted in the literature on Kalman filters; most importantly, we consider very general classes of cost functions, and temporal changes that are constrained only by a budget of variation, and are otherwise arbitrary (and in particular, non-parametric).

\label{pg:lit-operations}
A rich line of work in the literature considers concrete sequential decision problems embedded in an SA setting (namely, noisy observations of the cost or the gradient, where the underlying cost function is unknown).
Various studies consider dynamic pricing problems where the demand function is unknown, and noisy cost observations are obtained at each step; see recent works by \cite{Bro-Rus}, \cite{Den-Zwa}, and \cite{Keskin-Zeevi2014}, as well as the review by \cite{Ara-Cal2011} for both parametric and non-parametric approaches. Other studies consider a problem of inventory control with censored demand, where noisy observations of the gradient can be obtained in each step; see, e.g., \cite{Huh-Rus}, and \cite{Bes-Muh}. Other applications arise in queueing networks, online advertisements, wireless communications, and manufacturing systems, among other areas; see \cite{Kus-Yin2003} for an overview.

Most of the studies mentioned above focus on a setting in which the underlying environment (while unknown) is stationary. While several papers have considered settings where changes in the environment may occur, these papers typically assume a very specific structure on said changes (for example, considering dynamic pricing in the absence of capacity constrains, \cite{Kel-Rad} study a setting where demand is switching between two known demand functions according to a known Markov process; \cite{Bes-Zee-minimax} consider a similar problem in a setting where the timing of a single (known) change in the demand function is unknown). The current paper suggests a general framework to study stochastic optimization problems while allowing a broad array of changes in the underlying environment. In that sense, special cases of the formulation given in the current paper allow an extension of studies such as the ones mentioned above for a variety of non-stationary settings. 

\textbf{Structure of the paper.} \S\ref{sec:prob} contains the problem formulation. In \S\ref{sec:meta} we establish a principle that connects achievable regret of policies in the adversarial and non-stationary stochastic settings, in particular, proving that the property of sub-linearity of the regret can be 
carried over from the former to the latter.
\S\ref{sec:weakly} and \S\ref{sec:strictly} include the main rate optimality results for the convex and strongly convex settings, respectively.
\S 6 presents concluding remarks. Proofs can be found in Appendix A in the main text, and in Appendices B and C that appear in an online companion.

\vspace{-0.1cm}
\section{Problem Formulation}\label{sec:prob}
\vspace{-0.2cm}
Having already laid out in the previous section the key building blocks and ideas  behind our problem formulation, the purpose of the present section is to fill in any gaps and make that exposition more precise where needed; some repetition is expected but is kept to a minimum.

\textbf{Preliminaries and admissible polices.}
Let $\mathcal{X}$ be a convex, compact, non-empty {\it action set}, and $\mathcal{T} =~\left\{1,\ldots, T\right\}$ be the sequence of decision epochs. Let $\mathcal{F}$ be a class of sequences $f:=\{f_t:t=1,\ldots,T\}$ of convex cost functions from $\mathcal{X}$ into $\mathbb{R}$, that submit to the following two conditions:
\begin{enumerate}
  \item There is a finite number $G$ such that for any action $x \in \mathcal{X}$ and for any epoch $t\in \mathcal{T}$:\vspace{-0.2cm}
\begin{equation}\label{eq:G}
\left| f_{t}(x)\right| \leq G,\quad\quad\quad
\left\|\nabla f_{t}(x)\right\| \leq G.\vspace{-0.2cm}
\end{equation}
  \item There is some $\nu > 0$ such that\vspace{-0.2cm}
\begin{equation}\label{eq:interior}
\left\{x \in \mathbb{R}^{d} \;:\;\left\|x - x_{t}^{\ast}\right\| \leq \nu \right\}  \subset \mathcal{X} \quad\quad \mbox{for all } t\in \mathcal{T},\vspace{-0.2cm}
\end{equation}
\end{enumerate}
where $x^{\ast}_{t} := x^{\ast}_{t}\left(f_{t}\right) \in \argmin_{x\in\mathcal{X}}f_{t}(x)$.
Here $\nabla f_{t}(x)$ denotes the gradient of $f_{t}$ evaluated at point $x$, and $\|\cdot\|$ the Euclidean norm. In every epoch $t\in \mathcal{T}$ a decision maker selects a point $X_{t}\in \mathcal{X}$ and then observes a feedback $\phi_{t}:=~\phi_{t}(X_t,f_t)$ which takes one of two forms:\vspace{-0.2cm}
\begin{itemize}
\item noisy access to the cost, denoted by $\phi^{(0)}$, such that $ \mathbb{E}[\phi^{(0)}_{t}\left(X_{t},f_{t}\right) \;|X_{t}=x] = f_{t}(x)$;\vspace{-0.2cm}
\item noisy access to the gradient, denoted by $\phi^{(1)}$, such that $\mathbb{E}[\phi^{(1)}_{t}\left(X_{t},f_{t}\right) \;|X_{t}=x] = \nabla f_{t}(x)$,\vspace{-0.2cm}
\end{itemize}
For all $x\in\mathcal{X}$ and $f_{t}$, $t\in\left\{1,\ldots,T\right\}$, we will use $\phi_{t}(x,f_{t})$ to denote the feedback observed at epoch $t$, conditioned on  $X_{t}=x$, and $\phi$ will be used in reference to a generic feedback structure. The feedback signal is assumed to possess a second moment uniformly bounded over $\mathcal{F}$ and $\mathcal{X}$. \vspace{-0.1cm}
\begin{example}\textbf{\textup{(Independent noise)}} {\rm A conventional cost feedback structure is $\phi^{(0)}_{t}(x,f_{t}) = f_{t}(x) + \varepsilon_{t}$, where $\varepsilon_{t}$ are, say, independent Gaussian random variables with zero mean and variance uniformly bounded by $\sigma^{2}$. A gradient counterpart is $\phi^{(1)}_{t}(x,f_{t}) = \nabla f_{t}(x) + \varepsilon_{t}$, where $\varepsilon_{t}$ are independent Gaussian random vectors with zero mean and covariance matrices with entries uniformly bounded by $\sigma^{2}$.}
\qed\vspace{-0.1cm}
\end{example}
We next describe the class of admissible policies. Let $U$ be a random variable defined over a probability space $\left(\mathbb{U}, \mathcal{U},\mathbf{P}_{u}\right)$. Let $\pi_{1}:\mathbb{U}\rightarrow \mathbb{R}^{d}$ and $\pi_{t}:\mathbb{R}^{(t-1)k}\times\mathbb{U}\rightarrow \mathbb{R}^{d}$ for $t=2,3,\ldots$ be measurable functions, such that $X_{t}$, the action at time $t$, is given by\vspace{-0.3cm}
\begin{displaymath}
   X_{t} = \left\{
     \begin{array}{lr}
       \pi_{1}\left(U\right) & t=1,\quad\quad\quad\quad\\
       \pi_{t}\left(\phi_{t-1}\left(X_{t-1},f_{t-1}\right), \ldots, \phi_{1}\left(X_{1},f_{1}\right),U\right) & t=2,3,\ldots,\;\;
     \end{array}
   \right.\vspace{-0.3cm}
\end{displaymath}
where $k=1$ if $\phi = \phi^{(0)}$, namely, the feedback is noisy observations of the cost, and $k=d$ if $\phi = \phi^{(1)}$, namely, the feedback is noisy observations of the gradient. The mappings $\left\{\pi_{t}:\;t=1,\ldots,T\right\}$ together with the distribution $\mathbf{P}_{u}$ define the class of admissible policies with respect to feedback $\phi$. We denote this class by $\mathcal{P}_{\phi}$. We further denote by $\left\{\mathcal{H}_{t},\;t=1,\ldots,T\right\}$ the \emph{filtration} associated with a policy $\pi\in\mathcal{P}_{\phi}$, such that $\mathcal{H}_{1} = \sigma\left(U\right)$ and $\mathcal{H}_{t} = \sigma\left(\left\{\phi_{j}(X_{j},f_{j})\right\}_{j=1}^{t-1},U\right)$ for all $t\in~\left\{2,3,\ldots\right\}$. Note that policies in $\mathcal{P}_{\phi}$ are non-anticipating, i.e., depend only on the past history of actions and observations, and allow for randomized strategies via their dependence on $U$.

\textbf{Temporal uncertainty and regret.} As indicated already in the previous section, the class of sequences $\mathcal{F}$ is too ``rich,'' insofar as the latitude it affords nature. With that in mind, we further restrict the set of admissible cost function sequences, in particular, the manner in which its elements can change from one period to the other. Define the following notion of \emph{variation} based on the sup-norm:\vspace{-0.3cm}
\begin{equation}\label{eq:TVF}
\mbox{Var}(f_1,\ldots, f_T) \;:=\; \sum_{t=2}^{T}\|f_{t} - f_{t-1}\|,\vspace{-0.3cm}
\end{equation}
where for any bounded functions $g$ and $h$ from $\mathcal{X}$ into $\mathbb{R}$ we denote $\|g - h\| := \sup_{x\in \mathcal{X}}\left|g(x) - h(x)\right|$.
Let $\{V_t: t=1,2,\ldots\}$ be a non-decreasing sequence of real numbers such that $V_t \leq t$ for all $t$, $V_{1}=0$, and for normalization purposes set $V_2\geq 1$.
We refer to $V_T$ as the {\it variation budget} over $\mathcal{T}$. Using this as a primitive, define the corresponding
\emph{temporal uncertainty set}, as the set of admissible cost function sequences that are subject to the variation budget $V_T$ over the set of decision epochs  $\{1,\ldots,T\}$:\vspace{-0.3cm}
\begin{equation}\label{eq:TUS}
\mathcal{V} = \left\{\left\{f_{1},\ldots,f_{T}\right\} \subset \mathcal{F} \;:\; \sum_{t=2}^{T}\|f_{t} - f_{t-1}\| \leq V_{T} \right\}.\vspace{-0.3cm}
\end{equation}
While the variation budget places some restrictions on the possible evolution of the cost functions,   it still allows for many different temporal patterns: continuous change;  discrete shocks;  and  a non-constant rate of change. Two possible variations instances are illustrated in Figure \ref{fig:variation}; other variation patterns are considered in the numerical analysis described in Appendix~\ref{subsec:numerics}.
\begin{figure}[!ht]
\centering
\includegraphics[height=1.7in]{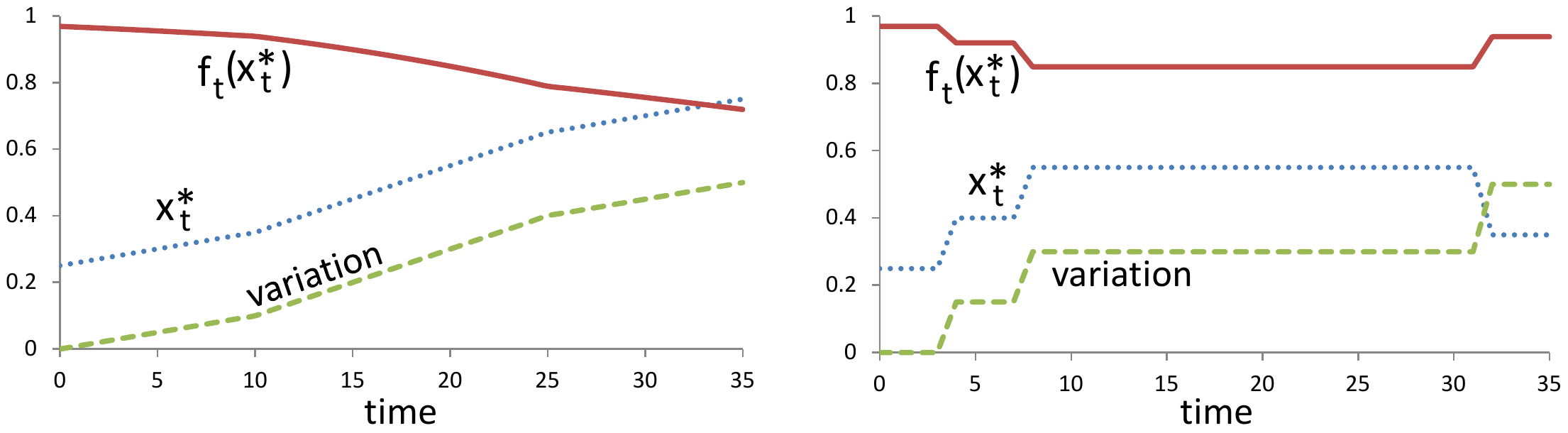}
\vspace{-0.2cm}\caption{\small \textbf{Variation instances within a temporal uncertainty set.}  Assume $\mathcal{X} = \left[0,1\right]$ and consider a sequence of quadratic cost functions of the form $f_{t}(x) = \frac{1}{2}x^{2} - b_{t}x + 1$. The change in the minimizer $x_{t}^{\ast} = b_{t}$, the optimal performance $f_{t}(x_{t}^{\ast}) = 1 - \frac{1}{2}b_{t}^{2}$, and the variation measured by (\ref{eq:TVF}), is illustrated for cases characterized by continuous changes (left), and ``jump" changes (right) in $b_{t}$. In both instances the variation budget is $V_{T}= 1/2$.
} \label{fig:variation}\vspace{-0.2cm}
\end{figure}

As described in \S1, the performance metric we adopt pits a policy $\pi$ against a dynamic oracle:\vspace{-0.2cm}
\begin{equation}\label{eq:regret}
\mathcal{R}^{\pi}_{\phi}(\mathcal{V},T) = \sup_{f\in \mathcal{V}}\left\{\mathbb{E}^{\pi}\left[\sum_{t=1}^{T}f_{t}(X_{t})\right] -\sum_{t=1}^{T}f_{t}(x_{t}^{\ast})\right\},\vspace{-0.2cm}
\end{equation}
where the expectation $\mathbb{E}^{\pi}\left[\cdot\right]$ is taken with respect to any randomness in the feedback, as well as in the policy's actions. Assuming a setup in which first a policy $\pi$ is chosen and then nature selects $f\in \mathcal{V}$ to maximize the regret, our formulation allows nature to select the worst possible sequence of cost functions for that policy, subject to the variation budget\footnote{In particular, while for the sake of simplicity and concreteness we use the above notation, our analysis applies to the case of sequences in which in every step only the next cost function is selected, in a fully adversarial manner that takes into account the realized trajectory of the policy and is subjected only to the bounded variation constraint.}. Recall that a policy $\pi$ is said to have \emph{sublinear} regret if $\mathcal{R}^{\pi}_{\phi}\left(\mathcal{V},T\right) = o\left(T\right)$, where for sequences $\{a_t\}$ and $\{b_t\}$ we write $a_t = o(b_t)$ if $a_t/b_t \to 0$ as $t \to \infty$.  Recall also that the {\it minimax regret}, being the minimal worst-case regret that can be guaranteed by an  admissible policy $\pi\in\mathcal{P}_{\phi}$, is given by:\vspace{-0.2cm}
\[
\mathcal{R}^{\ast}_{\phi}\left(\mathcal{V},T\right)
 = \inf_{\pi\in \mathcal{P}_{\phi}} \mathcal{R}^{\pi}_{\phi}\left(\mathcal{V},T\right).\vspace{-0.2cm}
\]
We refer to a policy $\pi$ as \emph{rate optimal} if it achieves the lowest possible growth rate of regret: there exists a constant $\bar{C}\geq 1$, independent of $V_{T}$ and $T$, such that for any $T \geq 1$,\vspace{-0.1cm}
\[
\mathcal{R}^{\pi}_{\phi}\left(\mathcal{V},T\right) \leq \bar{C}\cdot \mathcal{R}^{\ast}_{\phi}\left(\mathcal{V},T\right).\vspace{-0.1cm}
\]

\textbf{Contrasting with the adversarial online convex optimization paradigm.}
An OCO problem consists of a convex  set $\mathcal{X} \subset~\mathbb{R}^{d}$ and an a-priori unknown sequence $f =~\left\{f_{1},\ldots, f_{T}\right\}\in \mathcal{F}$ of convex cost functions. At any epoch $t$ the decision maker selects a point $X_{t}\in \mathcal{X}$,   and observes some feedback $\phi_{t}$. The efficacy of a policy over a given time horizon $T$ is typically measured relative to a benchmark which is defined by the \emph{single best action in hindsight}: the best {\it static} action fixed throughout the horizon, and chosen with benefit of having observed the sequence of cost functions. We use the notions of admissible, long-run-average optimal, and rate optimal policies in the adversarial OCO context as defined in the stochastic non-stationary context laid out before. Under the single best action benchmark, the objective is to minimize the regret incurred by an admissible online optimization algorithm  $\mathcal{A}$:\vspace{-0.1cm}
\begin{equation}\label{eq:advRegret}
\mathcal{G}^{\mathcal{A}}_{\phi}(\mathcal{F},T) = \sup_{f\in \mathcal{F}}\left\{\mathbb{E}^{\pi}\left[\sum_{t=1}^{T}f_{t}(X_{t})\right]
 - \min_{x\in \mathcal{X}}\left\{\sum_{t=1}^{T}f_{t}(x)\right\}\right\},\vspace{-0.1cm}
\end{equation}
where the expectation is taken with respect to possible randomness in the feedback and in the actions of the policy (We use the term ``algorithm'' to distinguish this from what we have defined as a ``policy,'' and this distinction will be important in what follows)\footnote{OCO settings typically allow sequences of cost functions that can adjust adversarially at each epoch. For the sake of consistency with the definition of (\ref{eq:regret}), in the above regret measure nature commits to a sequence of functions in advance.}. Interchanging the sum and $\min\left\{\cdot\right\}$ operators in the right-hand-side of (\ref{eq:advRegret}) we obtain the definition of regret in the non-stationary stochastic setting, as in (\ref{eq:regret}). As the next example shows, the dynamic oracle used as benchmark in the latter can be a significantly harder target than the single best action defining the static oracle in (\ref{eq:advRegret}).\vspace{-0.1cm}

\begin{example}\textbf{\textup{(Contrasting the static and dynamic oracles)}} {\rm Assume an action set $\mathcal{X} = \left[-1,2\right]$, and variation budget $V_{T} = 1$. Set\vspace{-0.3cm}
\begin{displaymath}
   f_{t}(x) = \left\{
     \begin{array}{lr}
       x^{2} & \text{ if } t \leq T/2\\
       x^{2} - 2x & \text{ otherwise},
     \end{array}
   \right.\vspace{-0.2cm}
\end{displaymath}
for any $x\in\mathcal{X}$. Then, the single best action is sub-optimal at each decision epoch, and \vspace{-0.2cm}
\[
\min_{x\in\mathcal{X}}\left\{\sum_{t=1}^{T}f_{t}(x)\right\} - \sum_{t=1}^{T}\min_{x\in\mathcal{X}}\left\{f_{t}(x)\right\}
\;=\; \frac{T}{4}.\quad\qed\vspace{-0.1cm}
\]}
\end{example}

Hence, algorithms that achieve performance that is ``close" to the static oracle in the adversarial OCO setting may perform quite poorly in the non-stationary stochastic setting (in particular they may, as the example above suggests, incur linear regret in that setting). Nonetheless, as the next section unravels, we will see that algorithms designed in the adversarial online convex optimization context can in fact  be adapted to perform well in the non-stationary stochastic setting laid out in this paper.

\section{A General Principle for Designing Efficient Policies}\label{sec:meta}
\vspace{-0.2cm}
In this section we will develop policies that operate well in non-stationary environments with given budget of variation $V_T$. Before exploring the question of what performance one may aspire to in the non-stationary variation constrained world, we first formalize what cannot be achieved.\vspace{-0.1cm}
\begin{prop} \textbf{\textup{(Linear variation budget implies linear regret)}} \label{prop:lin}
Assume a feedback structure $\phi \in~\left\{ \phi^{(0)},\phi^{(1)} \right\}$. If there exists a positive constant $C_{1}$ such that $V_{T}\geq C_{1}T$ for any $T\geq 1$, then there exists a positive constant $C_{2}$, such that for any admissible policy $\pi \in \mathcal{P}_{\phi}$,\vspace{-0.4cm}
  \[
  \mathcal{R}^{\pi}_{\phi}\left(\mathcal{V},T\right) \geq C_{2}T.\vspace{-0.2cm}
  \]
\end{prop}

\label{pg:sublinear variation} The proposition states that whenever the variation budget is at least of order $T$, \emph{any} policy which is admissible (with respect to the feedback) must incur a regret of order $T$, so under such circumstances it is not possible to have long-run-average optimality relative to the dynamic oracle benchmark. With that in mind, hereon we will focus on the case in which the variation budget is sublinear in $T$. We will show that in this case sublinear regret is achievable, \label{pg:variation-patterns} and study the behaviour of the minimax regret as a function of $V_T$ and $T$, when $T$ is large. We note that when $V_T$ is sublinear in $T$, the set $\cal{V}$ defined in \eqref{eq:TUS} is still very rich and includes many general patterns, such as sequences of functions $\left\{f_{t}\right\}$ that might change significantly from one period to the next but only do so rarely (a special case of which is a single change point), or sequences in which functions change often (even infinitely many times) but do so only locally. For example, consider the setting described in Figure~\ref{fig:variation} with a sequence of coefficients $\left\{b_{t}\right\}$ that does not converge, yet satisfies $\left|b_{t-1}-b_{t}\right|=t^{-1/2}$. Then, the variation budget $V_{T}$ is of order $\sqrt{T}$.

\textbf{A class of candidate policies.} We introduce a class of policies that leverages existing algorithms designed for fully adversarial environments. We denote by $\mathcal{A}$ an online optimization algorithm that given a feedback structure $\phi$ achieves a regret $\mathcal{G}^{\mathcal{A}}_{\phi}\left(\mathcal{F},T\right)$ (see (\ref{eq:advRegret})) with respect to the \emph{static} benchmark of the single best action. Consider the following generic ``restarting" procedure, which takes as input $\mathcal{A}$ and a batch size $\Delta_T$, with $1\leq \Delta_{T}\leq T$, and consists of restarting $\mathcal{A}$ every $\Delta_T$ periods. \label{pg:adjusted-formulation} To formalize this idea we first refine our definition of history-adapted policies and the actions they generate. Given a feedback $\phi$ and epochs $t'\geq1$, and $t>t'$ we define the history from $t'$ to $t$ by $\mathcal{H}_{t',t} = \sigma\left(\left\{\phi_{j}\left(X_{j},f_{j}\right)\right\}_{j=t'}^{t-1},U\right)$.
Then, for each restarting epoch $\tau\geq 1$, we have $X_{t} =~\mathcal{A}_{t-\tau}\left(\mathcal{H}_{\tau+1,t}\right)$ for each $\tau+1< t \leq \min\left\{\tau+\Delta_{T},T\right\}$, and $X_{\tau+1} =~\mathcal{A}_{1}\left(\mathcal{H}_{\tau-\Delta_{T}+1,\tau}\right)$.
Indeed, $X_{t}$ is $\mathcal{H}_{\tau+1,t}$-measurable for each $\tau+1< t \leq \min\left\{\tau+\Delta_{T}+1,T\right\}$, and $X_{\tau+1}$ is $\mathcal{H}_{\tau-\Delta_{T}+1,\tau}$-measurable.
The sequence of measurable mappings $\mathcal{A}_{t}$, $t=1,2,\ldots$ is prescribed by the algorithm $\mathcal{A}$, where we allow the initial action $\mathcal{A}_{1}$ to be based on information from the previous batch (if such exists).
The following procedure restarts $\mathcal{A}$ every $\Delta_{T}$ epochs. In what follows, let $\lceil\cdot\rceil$ denote the ceiling function (rounding its argument to the nearest larger integer).
\vspace{-0.5cm}\begin{center}
\line(1,0){490}
\end{center}
\vspace{-5mm}\textbf{Restarting procedure.} Inputs: an algorithm $\mathcal{A}$, and a batch size $\Delta_{T}$.\vspace{-0.2cm}
\begin{enumerate}
  \item Set $j = 1$\vspace{-0.2cm}
  \item Repeat while $j \leq \left\lceil T/\Delta_{T}\right\rceil$:\vspace{-0.2cm}
  \begin{enumerate}
    \item Set $\tau = \left(j-1\right)\Delta_{T}.$\vspace{-0.1cm}
    \item If $\tau = 0$ set $X_{\tau+1}=\mathcal{A}_{1}(U)$, otherwise set $X_{\tau+1}=\mathcal{A}_{1}\left(\mathcal{H}_{\tau-\Delta_{T}+1,\tau}\right)$;
        \\For any $t = \tau + 2,\ldots,\min\left\{T,\tau + \Delta_{T}\right\}$, select $X_{t} =~\mathcal{A}_{t-\tau}\left(\mathcal{H}_{\tau+1,t}\right)$.\vspace{-0.1cm}
    \item Set $j=j+1$.
    \vspace{-8mm}
  \end{enumerate}
\end{enumerate}
\begin{center}
\line(1,0){490}
\end{center}\vspace{-4mm}
Clearly $\pi\in \mathcal{P}_{\phi}$. Next we analyze the performance of policies defined via the restarting procedure, with suitable subroutine $\mathcal{A}$.

\textbf{First order performance.} The next result establishes a close connection between $\mathcal{G}^{\mathcal{A}}_{\phi}\left(\mathcal{F},T\right)$, the performance that is achievable in the adversarial environment by ${\cal A}$, and $\mathcal{R}^{\pi}_{\phi}\left(\mathcal{V},T\right)$, the performance in the non-stationary stochastic environment under temporal uncertainty set $\mathcal{V}$ of the restarting procedure that uses $\mathcal{A}$ as input.
\begin{theorem} \textbf{\textup{(Long-run-average optimality)}} \label{thm:meta}
Set a feedback structure $\phi \in\left\{ \phi^{(0)},\phi^{(1)} \right\}$. Let $\mathcal{A}$ be an OCO algorithm with $\mathcal{G}^{\mathcal{A}}_{\phi}\left(\mathcal{F},T\right) = o(T)$.
Let $\pi$ be the policy defined by the restarting procedure that uses $\mathcal{A}$ as a subroutine, with batch size $\Delta_{T}$. If $V_{T} = o(T)$, then for any $\Delta_{T}$ such that $\Delta_{T}= o(T/V_{T})$ and $\Delta_{T}\rightarrow\infty$ as $T\rightarrow\infty$,\vspace{-0.1cm}
\[
\mathcal{R}^{\pi}_{\phi}\left(\mathcal{V},T\right) \;=\; o(T).\vspace{-0.0cm}
\]
\end{theorem}
In other words, the theorem establishes the following meta-principle: whenever the variation budget is a sublinear function of the horizon length $T$, it is possible to construct a long-run-average optimal policy in the stochastic non-stationary SA environment by a suitable adaptation of an algorithm that achieves sublinear regret in the adversarial OCO environment. For a given structure of a function class and feedback signal, Theorem \ref{thm:meta} is meaningless unless there exists an algorithm with sublinear regret with respect to the single best action in the adversarial setting, under such structure. To that end, for the structures $\left(\mathcal{F}, \phi^{(0)}\right)$ and $\left(\mathcal{F}, \phi^{(1)}\right)$ an online gradient descent policy was shown to achieve sublinear regret in \cite{FLA2005}. We will see in the next sections that, surprisingly, the simple restarting mechanism introduced above allows to carry over not only first order optimality but also rate optimality from the OCO paradigm to the non-stationary SA setting.

\textbf{Key ideas behind the proof.} Theorem \ref{thm:meta} is driven directly by the next proposition that connects the performance of the restarting procedure with respect to the dynamic benchmark in the stochastic non-stationary environment, and the performance of the input subroutine algorithm $\mathcal{A}$ with respect to the single best action in the adversarial setting.\vspace{-0.1cm}
\begin{prop}{\textbf{\textup{(Connecting performance in OCO and non-stationary SA)}}} \label{prop:meta}
Set $\phi\in~\left\{\phi^{(0)}, \phi^{(1)}\right\}$. Let $\pi$ be the policy defined by the restarting procedure that uses $\mathcal{A}$ as a subroutine, with batch size $\Delta_{T}$. Then, for any $T\geq 1$,\vspace{-0.1cm}
\begin{equation}\label{eq:connection}
\mathcal{R}^{\pi}_{\phi}\left(\mathcal{V},T\right) \;\leq\;
\left\lceil\frac{T}{\Delta_{T}}\right\rceil\cdot
\mathcal{G}^{\mathcal{A}}_{\phi}(\mathcal{F},\Delta_{T}) \;+\; 2 \Delta_{T} V_{T}.
\end{equation}
\end{prop}
We next describe the high-level arguments. The main idea of the proof lies in analyzing the difference between the dynamic oracle and the static oracle benchmarks, used respectively in the OCO and the non-stationary SA contexts. We define a partition of the decision horizon 
into batches $\mathcal{T}_{1},\ldots,\mathcal{T}_{m}$ of size $\Delta_{T}$ each (except, possibly the last batch):\vspace{-0.2cm}
\begin{equation}\label{eq:batches}
\mathcal{T}_{j} = \left\{ t\; :\; (j-1)\Delta_{T} + 1 \leq t \leq \min\left\{j\Delta_{T},T\right\} \right\}, \text{ for all } \; j = 1,\ldots,m,\vspace{-0.2cm}
\end{equation}
where $m = \left\lceil T/\Delta_{T}\right\rceil$ is the number of batches. Then, one may write:\vspace{-0.3cm}
\[
\mathcal{R}_{\phi}^{\pi}(\mathcal{V},T)
 \;=\; \sup_{f\in \mathcal{V}}\left\{\sum_{j=1}^{m}\underbrace{\left(\mathbb{E}^{\pi}\left[\sum_{t \in \mathcal{T}_{j}}f_{t}(X_{t})\right]
 - \min_{x\in \mathcal{X}}\left\{\sum_{t \in \mathcal{T}_{j}}f_{t}(x)\right\}\right)}_{J_{1,j}}
 \;+\; \sum_{j=1}^{m}\underbrace{\left(\min_{x\in \mathcal{X}}\left\{\sum_{t \in \mathcal{T}_{j}}f_{t}(x)\right\}
  -  \sum_{t \in \mathcal{T}_{j}}f_{t}(x^{\ast}_{t})\right)}_{J_{2,j}}\right\}.\vspace{-0.3cm}
\]
The regret with respect to the dynamic benchmark is represented as two sums. The first, $\sum_{j=1}^{m}J_{1,j}$, sums the regret terms with respect to the single best action within each batch $\mathcal{T}_{j}$, which are each bounded by $\mathcal{G}^{\mathcal{A}}_{\phi}(\mathcal{F},\Delta_{T})$. Noting that there are $\left\lceil T/\Delta_{T}\right\rceil$ batches, this gives rise to the first term on the right-hand-side of (\ref{eq:connection}). The second sum, $\sum_{j=1}^{m}J_{2,j}$, is the sum of differences between the performances of the single best action benchmark and the dynamic benchmark within each batch. The latter is driven by the rate of functional change in the batch. While locally this gap can be large, we show that given the variation budget the second sum is at most of order $\Delta_{T} V_{T}$. This leads to the result of the proposition.  Intuitively, Proposition~\ref{prop:meta} highlights the following tradeoff. When $\Delta_{T}$ is large, the performance of a ``good" subroutine policy approaches the one of the static oracle; when $\Delta_{T}$ is small, the sequence of static oracles approaches the dynamic oracle. Theorem \ref{thm:meta} follows by balancing this tradeoff. \qed

\textbf{Remark} \textbf{(Alternative forms of feedback)} The principle laid out in Theorem \ref{thm:meta} can also be derived for other forms of feedback using Proposition~\ref{prop:meta}. For example, the proof of Theorem~\ref{thm:meta} holds for settings with richer feedback structures, such as noiseless access to the full cost function (\citealt{Zin2003}), or a multi-point access (\citealt{Aga-Dek-Xia2010}).

\vspace{-0.1cm}
\section{Rate Optimality: The General Convex Case}\label{sec:weakly}
\vspace{-0.2cm}

A natural question arising from the analysis of \S\ref{sec:meta} is what type of performance one may achieve in non-stationary environments and how does such performance depend on the variation one may face. We first focus on the feedback structure $\phi^{(1)}$, for which rate optimal polices are known in the OCO setting (as these will serve as inputs for the restarting procedure). To answer such a question, we first develop a lower bound for a subclass of problems and then establish that such a lower bound is achievable.

\textbf{A lower bound on achievable performance.}  We will establish a fundamental bound on the performance of any admissible policy under the following technical assumption on the structure of the gradient feedback signal (a cost feedback counterpart will be provided in the next section).\vspace{-0.2cm}
\begin{assumption}{\textbf{\textup{(Gradient feedback structure)}}}\label{ass:noise}
\begin{enumerate}\vspace{-0.2cm}
  \item $\phi_{t}^{(1)}(x,f_{t}) = \nabla f_{t}(x) + \varepsilon_{t}$ for any $f\in\mathcal{F}$, $x\in\mathcal{X}$, and $t\in \mathcal{T}$, where $\varepsilon_{t}$, $t\geq 1$, are iid random vectors with zero mean and covariance matrix with bounded entries. \vspace{-0.1cm}
  \item Let $G(\cdot)$ be the cumulative distribution function of $\varepsilon_{t}$. There exists a constant $\tilde{C}$ such that for any $a\in\mathbb{R}^{d}$, $\int\log\left(\frac{dG(y)}{dG(y+a)}  \right)dG(y)\leq \tilde{C} \left\|a\right\|^{2}$. \vspace{-0.1cm}
\end{enumerate}
\end{assumption}
For the sake of concreteness we impose an additive noise feedback structure, given in the first part of the assumption. This simplifies notation and streamlines proofs, but otherwise is not essential. The key properties that are needed are: $\mathbb{P}\left(\phi_{t}^{(1)}(x,f_{t})\in A\right) > 0$ for any $f\in \mathcal{F}$, $t\in\mathcal{T}$, $x\in\mathcal{X}$, and $A\subset \mathbb{R}^{d}$; and that the feedback observed at any epoch $t$, conditioned on the action $X_{t}$, is independent of the history that is available at that epoch. Given the structure imposed in the first part of the assumption, the second part implies that if gradients of two cost functions are ``close" to each other, the probability measures of the observed feedbacks are also ``close". The structure imposed by Assumption~1 is satisfied in many settings. For instance, it applies to Example~1 (with $\mathcal{X}\subset \mathbb{R}$), with $\tilde{C} =~1/2\sigma^{2}$.

\begin{theorem}\textbf{\textup{(Lower bound on achievable performance)}}\label{thm:weakGrad-low}
Let Assumption~1 hold. Then, there exists a constant $C>0$, independent of $T$ and $V_{T}$, such that for any policy $\pi\in\mathcal{P}_{\phi^{(1)}}$ and for all $T\geq 1$:\vspace{-0.3cm}
\[
\mathcal{R}^{\pi}_{\phi^{(1)}}(\mathcal{V},T) \;\geq\; C\cdot V_{T}^{1/3}T^{2/3}.\vspace{-0.2cm}
\]
\end{theorem}

\paragraph*{}\vspace{-0.15cm}
\textbf{Key ideas in the proof of Theorem \ref{thm:weakGrad-low}.} For two probability measures $\mathbb{P}$ and $\mathbb{Q}$ on a probability space $\mathcal{Y}$, let\vspace{-0.1cm}
\begin{equation}\label{eq:KL}
\mathcal{K}\left(\mathbb{P}\|\mathbb{Q}\right)
 \;=\; \mathbb{E}\left[\log\left( \frac{d\mathbb{P}\left\{Y\right\}}{d\mathbb{Q}\left\{Y\right\}}  \right)   \right],\vspace{-0.0cm}
\end{equation}
where $\mathbb{E}\left[\cdot\right]$ is the expectation with respect to $\mathbb{P}$, and $Y$ is a random variable defined over $\mathcal{Y}$. This quantity is is known as the Kullback-Leibler divergence. To establish the result, we consider sequences from a subset of $\mathcal{V}$ defined in the following way: in the beginning of each batch of size $\tilde{\Delta}_{T}$ (nature's decision variable), one of two ``almost-flat" functions is independently drawn according to a uniform distribution, and set as the cost function throughout the next $\tilde{\Delta}_{T}$ epochs. Then, the distance between these functions, and the batch size $\tilde{\Delta}_{T}$ are tuned such that: $(a)$ any drawn sequence must maintain the variation constraint; and $(b)$ the functions are chosen to be ``close" enough while the batches are sufficiently short, such that distinguishing between the two functions over the batch is subject to a significant error probability, yet the two functions are sufficiently ``separated" to maximize the incurred regret. (Formally, the KL divergence is bounded throughout each batch, and hence any admissible policy trying to identify the current cost function 
can only do so with a strictly positive error probability.)


 \textbf{Upper bound on performance.} In order to establish that the lower bound is achievable, we show that the restarting procedure introduced in \S\ref{sec:meta} enables to carry over the property of rate optimality from the adversarial setting to the non-stationary stochastic setting.
As a subroutine algorithm, we will use an adaptation of the \emph{online gradient descent} (OGD) algorithm introduced by \cite{Zin2003}:
\vspace{-5mm}\begin{center}
\line(1,0){490}
\end{center}
\vspace{-5mm}\textbf{OGD algorithm.} Input: a decreasing sequence of non-negative real numbers $\left\{\eta_{t}\right\}_{t=2}^{T}$.
\begin{enumerate}\vspace{-0.2cm}
  \item Select some $X_{1}\in \mathcal{X}$\vspace{-0.2cm}
  \item For any $t =1,\ldots,T-1$, set $X_{t+1}\;=\; P_{\mathcal{X}}\left(X_{t} - \eta_{t+1}\phi^{(1)}_{t}(X_{t},f_{t})\right)$, where $P_{\mathcal{X}}\left(y\right) = \arg\min_{x\in \mathcal{X}}\left\|x-y\right\|$ is the Euclidean projection operator on $\mathcal{X}$.\vspace{-6mm}
\end{enumerate}
\begin{center}
\line(1,0){490}
\end{center}\vspace{-5mm}
For any value of $\tau$ that is dictated by the restarting procedure, the OGD algorithm can be defined via the sequence of mappings $\left\{\mathcal{A}_{t-\tau}\right\}$, $t\geq\tau+1$, as follows:\vspace{-0.2cm}

\begin{displaymath}
   \mathcal{A}_{t-\tau}\left(\mathcal{H}_{\tau,t}\right) = \left\{
     \begin{array}{lr}
       \text{some } X_{1} & \text{ if } t=\tau+1\\
       P_{\mathcal{X}}\left(X_{t-1} - \eta_{t-\tau}\phi^{(1)}_{t-1}\right) & \text{ if } t>\tau+1,
     \end{array}
   \right.\vspace{-0.2cm}
\end{displaymath}
for any epoch $t\geq \tau+1$. 
For the structure $\left(\mathcal{F},\phi^{(1)}\right)$ of convex cost functions and noisy gradient access, \cite{FLA2005} consider the OGD algorithm with $X_{1}=0$ and the selection $\eta_{t} = r/G\sqrt{T}$, $t=2,\ldots, T$. Here $r$ denotes the radius of the action set: $r = \inf\left\{y>0:\;\mathcal{X}\subseteq \mathbf{B}_{y}(x) \text{ for some } x\in \mathbb{R}^{d}\right\}$, where $\mathbf{B}_{y}(x)$ is a ball with radius $y$, centered at point $x$, and show that this algorithm achieves a regret of order $\sqrt{T}$ in the adversarial setting. For completeness, we prove in Lemma \ref{lem:lowerFlax} (given in Appendix~\ref{app:sba}) that under Assumption~1 this performance cannot be improved upon in the adversarial OCO setting.

We next characterize the regret of the restarting procedure that uses the OGD policy as an input.
\begin{theorem}\textbf{\textup{(Performance of restarted OGD under noisy gradient access)}}\label{thm:weakGrad-up}
 Consider the feedback setting $\phi = \phi^{(1)}$, and let $\pi$ be the policy defined by the restarting procedure with a batch size $\Delta_{T} =~\left\lceil\left(T/V_{T}\right)^{2/3}\right\rceil$, and the OGD algorithm parameterized by $\eta_{t} = \frac{r}{G\sqrt{\Delta_{T}}}$, $t=2,\ldots, \Delta_{T}$ as a subroutine. Then, there is some finite constant $\bar{C}$, independent of $T$ and $V_{T}$, such that for all $T\geq 2$:\vspace{-0.3cm}
\[
\mathcal{R}^{\pi}_{\phi}(\mathcal{V},T) \;\leq\; \bar{C}\cdot V_{T}^{1/3}T^{2/3}.\vspace{-0.3cm}
\]
\end{theorem}
Recalling the connection between the regret in the adversarial setting and the one in the non-stationary SA setting (Proposition \ref{prop:meta}), the result of the theorem is essentially a direct consequence of bounds in the OCO literature. In particular, Flaxman et al. (\citeyear[Lemma 3.1]{FLA2005}) provide a bound on $\mathcal{G}^{\mathcal{A}}_{\phi^{(1)}}(\mathcal{F},\Delta_{T})$ of order $\sqrt{\Delta_{T}}$, and the result follows by balancing the terms in (\ref{eq:connection}) by a proper selection of $\Delta_{T}$.

When selecting a large batch size, the ability to track the single best action within each batch improves, but the single best action within a certain batch may have substantially worse performance than that of the dynamic oracle. On the other hand, when selecting a small batch size, the performance of tracking the single best action within each batch gets worse, but over the whole horizon the series of single best actions (one for each batch) achieves a performance that approaches the dynamic oracle.

We note that Theorem \ref{thm:weakGrad-up} holds for any (deterministic or random) initial action of the subroutine OGD algorithm; a practical special case is one in which the initial action of any batch $j>1$ is determined by taking one further gradient step from the last action of batch $j-1$.
 \label{pg:OGD-random}

Recalling the lower bound in Theorem \ref{thm:weakGrad-low}, Theorem \ref{thm:weakGrad-up} implies that the performance of restarted OGD is rate optimal, and the minimax regret under structure $\left(\mathcal{V},\phi^{(1)}\right)$ is:\vspace{-0.4cm}
\[
\mathcal{R}^{\ast}_{\phi^{(1)}}(\mathcal{V},T) \;\asymp\; V_{T}^{1/3}T^{2/3}.\vspace{-0.4cm}
\]
Roughly speaking, this characterization provides a mapping between the variation budget $V_T$ and the minimax regret under noisy gradient observations. For example, when $V_{T} = T^{\alpha}$ for some $0\leq \alpha \leq 1$, the minimax regret is of order $T^{(2+\alpha)/3}$, hence we obtain the minimax regret in a full spectrum of variation scales, from order $T^{2/3}$ when the variation is a constant (independent of the horizon length), up to order $T$ that corresponds to the case where $V_{T}$ scales linearly with $T$ (consistent with Proposition~\ref{prop:lin}).

\vspace{-0.5cm}
\paragraph*{\quad Alternative algorithms.} \label{pg:continuous}
 While the restarting procedure (together with suitable balancing of the batch size) can be used as a template for deriving ``good" policies in non-stationary stochastic settings, it serves mainly as a tool to articulate a general and unified principle for designing rate optimal policies. Indeed, rate optimal performance may also be achieved by taking alternative paths that may be considered as more appealing from practical points of view. One of these may rely on attempting to directly re-tune the parameters of the subroutine OCO algorithm. While, not surprisingly, OGD-type policies with classical step size selections (such as $1/t$ or $1/\sqrt{t}$) may perform poorly in non-stationary environments (see Example \ref{ex:LinearRegretUnderOGD} in Appendix \ref{app:add}), we establish next that one may fine tune such a policy to achieve rate optimality, matching the lower bound given in Theorem \ref{thm:weakGrad-low}.
\begin{prop}\label{thm:ebe}
\textbf{\textup{(Optimal tuning of OGD)}} Assume $\phi = \phi^{(1)}$, and let $\pi$ be the OGD algorithm with $\eta_{t} = \frac{r}{G}\left(V_{T}/T\right)^{1/3}$, $t=2,\ldots, T$. Then, there exists a finite constant $\bar{C}$, independent of $T$ and $V_{T}$, such that for all $T\geq 2$:\vspace{-0.2cm}
\[
\mathcal{R}^{\pi}_{\phi}(\mathcal{V},T) \;\leq\; \bar{C}\cdot V_{T}^{1/3}T^{2/3}.\vspace{-0.1cm}
\]
\end{prop}
The key in tuning the OGD algorithm to achieves rate optimal performance in the non-stationary SA setting is a suitable adjustment of the step size sequence as a function of the variation budget $V_{T}$. A sequence of ``larger" steps that converge ``slower" to zero allows the the policy to respond efficiently to potential changes in the environment; the larger the variation budget is (relative to the horizon length $T$), the larger the step sizes that are required in order to ``keep up" with the potential changes.

\paragraph*{}\vspace{-0.5cm}
\textbf{Noisy access to the function value.} Considering the feedback structure $\phi^{(0)}$ and the class $\mathcal{F}$, \cite{FLA2005} show that in the adversarial OCO setting, a modification of the OGD algorithm can be tuned to achieve regret of order $T^{3/4}$. There is no indication that this regret rate is the best possible, and to the best of our knowledge, under cost observations and general convex cost functions, the question of rate optimality is an open problem in the adversarial OCO setting. By Proposition \ref{prop:meta}, the regret of order $T^{3/4}$ that is achievable in the OCO setting implies that a regret of order $V_{T}^{1/5}T^{4/5}$ is achievable in the non-stationary SA setting, by applying the restarting procedure. While at present, we are not aware of any algorithm that guarantees a lower regret rate for arbitrary action spaces of dimension $d$, we conjecture that a rate optimal algorithm in the OCO setting can be lifted to a rate optimal procedure in the non-stationary stochastic setting by applying the restarting procedure.\footnote{Considering the special case of linear cost functions, there are known policies that guarantee regret of order $\sqrt{T}$ relative to the single best action in adversarial settings (see, e.g., \cite{Kal-Vem2003} for full access to the function, and \cite{McM-Blu2004} for point feedback). Using an adaptation of such policy as a subroutine of the restarting procedure would guarantee regret of order $V_{T}^{1/3}T^{2/3}$ relative to the dynamic oracle in 
our setting; a matching lower bound may be obtained by a rather straightforward adaptation of the proof of Theorem \ref{thm:weakGrad-low}.}
The next section 
supports this conjecture examining the case of strongly convex cost functions.

\vspace{-0.2cm}
\section{Rate Optimality: The Strongly Convex Case}\label{sec:strictly}
\vspace{-0.25cm}
\textbf{Preliminaries.} We now focus on the class of strongly convex functions $\mathcal{F}_{s}\subseteq \mathcal{F}$, defined such that in addition to the conditions that are stipulated by membership in $\mathcal{F}$, for a finite number $H>0$, the sequence $\left\{f_{t}\right\}$ satisfies \vspace{-0.2cm}
\begin{equation}\label{eq:H}
H\mathbf{I}_{d} \;\preceq\; \nabla^{2} f_{t}(x) \;\preceq\; G\mathbf{I}_{d} \quad \text{ for all } x\in\mathcal{X}, \text{ and all } t\in \mathcal{T},\vspace{-0.2cm}
\end{equation}
where $\mathbf{I}_{d}$ denotes the $d$-dimensional identity matrix. Here for two square matrices of the same dimension $A$ and $B$, we write $A \preceq B$ to denote that $B-A$ is positive semi-definite, and $\nabla^{2} f(x)$ denotes the Hessian of $f(\cdot)$, evaluated at point $x\in\mathcal{X}$; for the sake of simplicity we assume that $G$ is a unified bound that also appears in (\ref{eq:G}). In the presence of strongly convex cost functions, it is well known that local properties of the functions around their minimum play a key role in the performance of sequential optimization procedures. To localize the analysis, we adapt the functional variation definition so that it is measured by the uniform norm over the \emph{convex hull of the minimizers}, denoted by:\vspace{-0.2cm}
\[
\mathcal{X}^{\ast} = \left\{x\in \mathbb{R}^{d} \;: x= \sum_{t=1}^{T}\lambda_{t}x_{t}^{\ast},\; \sum_{t=1}^{T}\lambda_{t}=1,\; \lambda_{t}\geq 0 \; \text{ for all } t\in\mathcal{T}\right\}.\vspace{-0.3cm}
\]
Using the above, we measure variation by:\vspace{-0.2cm}
\begin{equation}\label{eq:TVFconv}
\mbox{Var}_{s}(f_1,\ldots,f_T) := \sum_{t=2}^{T}\sup_{x\in \mathcal{X}^{\ast}}\left|f_{t}(x) - f_{t-1}(x)\right|.\vspace{-0.3cm}
\end{equation}
Given the class $\mathcal{F}_s$ and a variation budget $V_{T}$, we define the temporal uncertainty set as follows:\vspace{-0.3cm}
\[
\mathcal{V}_{s} = \left\{f = \left\{f_{1},\ldots,f_{T}\right\} \subset \mathcal{F}_{s} \;:\; \mbox{Var}_{s}(f_1,\ldots,f_T) \leq V_{T} \right\}.\vspace{-0.3cm}
\]
We note that the proof of Proposition \ref{prop:meta} effectively holds without change under the above structure. Hence first order optimality is carried over from the OCO setting, as long as $V_{T}$ is sublinear. We next examine rate-optimality results.

\subsection{Noisy access to the gradient}\label{subsec:strictGrad}\vspace{-0.2cm}
For the class $\mathcal{F}_s$ and gradient feedback $\phi_{t}\left(x,f_{t}\right) = \nabla f_{t}(x)$, \cite{Haz-Aga-Kal2007} consider the OGD algorithm with a tuned selection of $\eta_{t} = 1/Ht$ for $t=2,\ldots T$, and provide in the OCO framework a regret guarantee of order $\log T$ (relative to the single best action benchmark). For completeness, we provide in Appendix \ref{app:sba} (Lemma \ref{lem:sbaGrad}) a simple adaptation of this result to the case of \emph{noisy} gradient access and an arbitrary, random $X_{1}$. \cite{HazKal2011} show that this algorithm is rate optimal in the OCO setting under strongly convex functions and a class of unbiased gradient feedback.\footnote{In fact, \cite{HazKal2011} show that even in a stationary stochastic setting with strongly convex cost function and a class of unbiased gradient access, any policy must incur regret of at least order $\log T$ compared to a static benchmark.}

\begin{theorem}\textbf{\textup{(Rate optimality for strongly convex functions and noisy gradient access)}}\label{thm:strictGrad}
\begin{enumerate}\vspace{-0.2cm}
  \item Consider the feedback structure $\phi = \phi^{(1)}$, and let $\pi$ be the policy defined by the restarting procedure with a batch size $\Delta_{T} =~\left\lceil\sqrt{T \log T/V_{T}}\;\right\rceil$, and the OGD algorithm parameterized by $\eta_{t} = \left(Ht\right)^{-1}$, $t=2,\ldots, \Delta_{T}$ as a subroutine. Then, there exists a finite positive constant $\bar{C}$, independent of $T$ and $V_{T}$, such that for all $T\geq 2$:\vspace{-0.4cm}
\[
\mathcal{R}^{\pi}_{\phi}(\mathcal{V}_{s},T) \leq \bar{C}\cdot \log\left(\frac{T}{V_{T}} + 1\right)\sqrt{V_{T}T}.\vspace{-0.3cm}
\]
  \item Let Assumption~1 hold. Then, there exists a constant $C>0$, independent of $T$ and $V_{T}$, such that for any policy $\pi\in\mathcal{P}_{\phi^{(1)}}$ and for all $T\geq 1$:\vspace{-0.3cm}
\[
\mathcal{R}^{\pi}_{\phi}(\mathcal{V}_{s},T) \;\geq\; C\cdot \sqrt{V_{T}T}.\vspace{-0.3cm}
\]
\end{enumerate}
\end{theorem}
Up to a logarithmic term, Theorem \ref{thm:strictGrad} establishes rate optimality in the non-stationary SA setting of the policy defined by the restarting procedure with the tuned $OGD$ algorithm as a subroutine. In \S\ref{sec:rem} we show that one may achieve a performance of $O(\sqrt{V_T T})$ through a slightly modified procedure, and hence the \emph{minimax regret} under structure $\left(\mathcal{F}_{s},\phi^{(1)}\right)$ is:\vspace{-0.3cm}
\[
\mathcal{R}^{\ast}_{\phi^{(1)}}(\mathcal{V}_{s},T) \;\asymp\; \sqrt{V_{T}T}.\vspace{-0.2cm}
\]
Theorem \ref{thm:strictGrad} further validates the ``meta-principle" in the case of strongly convex functions and noisy gradient feedback: rate optimality in the adversarial setting (relative to the single best action benchmark) can be adapted by the restarting procedure to guarantee an essentially optimal regret rate in the non-stationary stochastic  setting (relative to the dynamic benchmark).

The first part of Theorem \ref{thm:strictGrad} is derived directly from Proposition \ref{prop:meta}, by plugging in a bound on $\mathcal{G}^{\mathcal{A}}_{\phi^{(1)}}(\mathcal{F}_{s},\Delta_{T})$ of order $\log T$ (given by Lemma \ref{lem:sbaGrad} in the case of noisy gradient access), and a tuned selection of $\Delta_{T}$. The proof of the second part follows by arguments similar to the ones used in the proof of Theorem \ref{thm:weakGrad-low}, adjusting for strongly convex cost functions.

\vspace{-0.1cm}
\subsection{Noisy access to the cost}\label{subsec:strictCost}
\vspace{-0.2cm}
We now consider the structure $\left(\mathcal{V}_{s},\phi^{(0)}\right)$, in which the cost functions are strongly convex and the decision maker has noisy access to the cost. In order to show that rate optimality is carried over from the adversarial setting to the non-stationary stochastic setting, we first need to introduce an algorithm that is rate optimal in the adversarial setting under the structure $\left(\mathcal{F}_{s},\phi^{(0)}\right)$.

\textbf{Estimated gradient step.} For a small $\delta$, we denote by $\mathcal{X}_{\delta}$ the $\delta$-interior of the action set $\mathcal{X}$:\vspace{-0.3cm}
\[
\mathcal{X}_{\delta} = \left\{x\in \mathcal{X}\;:\; \mathbf{B}_{\delta}(x) \subseteq \mathcal{X}\right\}.\vspace{-0.3cm}
\]
We assume access to the projection operator $P_{\mathcal{X}_{\delta}}\left(y\right) = \arg\min_{x\in \mathcal{X}_{\delta}}\left\|x-y\right\|$ on the set $\mathcal{X}_{\delta}$.

For $k=1,...,d$, let $e^{(k)}$ denote the unit vector with $1$ at the $k^{\text{th}}$ coordinate. The \emph{estimated gradient step} (EGS) algorithm is defined through three sequences of real numbers $\left\{h_{t}\right\}$, $\left\{a_{t}\right\}$, and $\left\{\delta_{t}\right\}$, where\footnote{For any $t$ such that $\nu < \delta_{t}$, one may use the numbers $h'_{t} = \delta'_{t} =\min\left\{\nu, \delta_{t}\right\}$ instead, with the rate optimality obtained in Lemma \ref{lem:sbaRegret} remaining unchanged.} $\nu \geq \delta_{t}\geq h_{t}$ for all $t\in \mathcal{T}$:
\vspace{-6mm}\begin{center}
\line(1,0){490}
\end{center}
\vspace{-4mm}\textbf{EGS algorithm.} Inputs: decreasing sequences of real numbers $\left\{a_{t}\right\}_{t=1}^{T-1}, \left\{h_{t}\right\}_{t=1}^{T-1}, \left\{\delta_{t}\right\}_{t=1}^{T-1}$.
\begin{enumerate}\vspace{-0.2cm}
  \item Select some initial point $X_{1} = Z_{1}$ in $\mathcal{X}$.\vspace{-0.2cm}
  \item For each $t = 1,\ldots,T-1$:\vspace{-0.2cm}
  \begin{enumerate}
    \item Draw $\psi_{t}$ uniformly over the set $\left\{\pm e^{(1)},\ldots, \pm e^{(d)}\right\}$.\vspace{-0.1cm}
    \item Compute stochastic gradient estimate \quad $\hat{\nabla}_{h_{t}}f_{t}(Z_t) \:=\: h_{t}^{-1}\phi_{t}^{(0)}(Z_t + h_t\psi_{t})\psi_{t}$.\vspace{-0.1cm}
    \item Update \quad $Z_{t+1} \:=\: P_{\mathcal{X}_{\delta_{t}}}\left(Z_{t} - a_{t}\hat{\nabla}_{h_{t}}f_{t}(Z_{t})\right)$.\vspace{-0.1cm}
    \item Select the action \quad $X_{t+1} \:=\: Z_{t+1} + h_{t+1}\psi_{t}$.
  \end{enumerate}
\end{enumerate}\vspace{-1.1cm}
\begin{center}
\line(1,0){490}
\end{center}\vspace{-3mm}
For any value of $\tau$ dictated by the restarting procedure, the EGS policy can be formally defined by\vspace{-0.2cm}

\begin{displaymath}
   \mathcal{A}_{t-\tau} = \left\{
     \begin{array}{lr}
       \text{ some } Z_{1} & \text{ if } t=\tau+1\\
       Z_{t-\tau} + h_{t-\tau}\psi_{t-\tau-1} & \text{ if } t>\tau+1.
     \end{array}
   \right.\vspace{-0.2cm}
\end{displaymath}

Note that $\Expect[\hat{\nabla}_{h}f_{t}(Z_t)|X_t] = \nabla f_{t}(Z_t)$ (cf. \citealt{Nem-Yud1983}, chapter 7),\label{pg:EGS} and that the EGS algorithm essentially consists of estimating a stochastic direction of improvement and following this direction. In Lemma \ref{lem:sbaRegret} (Appendix \ref{app:sba}) we show that when tuned by $a_{t} = 2d/Ht$ and $\delta_{t} = h_{t} = a_{t}^{1/4}$ for all $t\in \left\{1,\ldots,T-1\right\}$, the EGS algorithm achieves a regret of order $\sqrt{T}$ compared to a single best action in the adversarial setting under structure $\left(\mathcal{F}_{s},\phi^{(0)}\right)$. In Lemma~\ref{lem:lowerEGS} (Appendix \ref{app:sba}) we establish that under Assumption~2 (given below) this performance is rate optimal in the adversarial setting.

Before analyzing the minimax regret in the non-stationary SA setting, let us introduce a counterpart to Assumption~1 for the case of cost feedback, that will be used in deriving a lower bound on the regret.
\begin{assumption}{\textbf{\textup{(Cost feedback structure)}}}\label{ass:noise2}
\begin{enumerate}\vspace{-0.2cm}
  \item $\phi_{t}^{(0)}(x,f_{t}) = f_{t}(x) + \varepsilon_{t}$ for any $f\in\mathcal{F}$, $x\in\mathcal{X}$, and $t\in \mathcal{T}$, where $\varepsilon_{t}$, $t\geq 1$, are iid random variables with zero mean and bounded variance.
  \item Let $G(\cdot)$ be the cumulative distribution function of $\varepsilon_{t}$. Then, there exists a constant $\tilde{C}$ such that for any $a\in\mathbb{R}$, $\int\log\left(\frac{dG(y)}{dG(y+a)}  \right)dG(y) \leq\tilde{C}\cdot a^{2}$. \vspace{-0.2cm}
\end{enumerate}
\end{assumption}
\begin{theorem}\textbf{\textup{(Rate optimality for strongly convex functions and noisy cost access)}}\label{thm:strictCost}
\begin{enumerate}\vspace{-0.2cm}
  \item Consider the feedback structure $\phi = \phi^{(0)}$, and let $\pi$ be the policy defined by the restarting procedure with EGS parameterized by $a_{t} = 2d/Ht$,
$h_{t} = \delta_{t} =  \left(2d/Ht\right)^{1/4}$, $t = 1,\ldots, T-1$, as subroutine, and a batch size $\Delta_{T} =~\left\lceil\left(T/V_{T}\right)^{2/3}\right\rceil$. Then, there exists a finite constant $\bar{C}>0$, independent of $T$ and $V_{T}$, such that for all $T\geq 2$:\vspace{-0.2cm}
\[
\mathcal{R}^{\pi}_{\phi}(\mathcal{V}_{s},T) \;\leq\; \bar{C}\cdot V_{T}^{1/3}T^{2/3}.\vspace{-0.4cm}
\]
  \item Let Assumption~2 hold. Then, there exists a constant $C>0$, independent of $T$ and $V_{T}$, such that for any policy $\pi\in\mathcal{P}_{\phi^{(0)}}$ and for all $T\geq 1$:\vspace{-0.2cm}
\[
\mathcal{R}^{\pi}_{\phi}(\mathcal{V}_{s},T) \;\geq\; C\cdot V_{T}^{1/3}T^{2/3}.\vspace{-0.3cm}
\]
\end{enumerate}
\end{theorem}
Theorem \ref{thm:strictCost} again establishes the ability to ``port over" rate optimality from the adversarial OCO setting to the non-stationary stochastic setting, this time under structure $\left(\mathcal{F}_{s},\phi^{(0)}\right)$. The theorem establishes a characterization of the minimax regret under structure $\left(\mathcal{V}_{s},\phi^{(0)}\right)$:\vspace{-0.3cm}
\[
\mathcal{R}^{\ast}_{\phi^{(0)}}(\mathcal{V}_{s},T) \;\asymp\; V_{T}^{1/3}T^{2/3}.\vspace{-0.3cm}
\]

\textbf{Illustrative Numerical Results.}
\vspace{-0.1cm}
\label{pg:numerics} In Appendix \ref{subsec:numerics} we illustrate the upper bounds on the regret by numerical experiments measuring the average regret that is incurred in the presence of various patterns of changing costs of fixed variation, different feedback structures and noise. Under noisy gradient access $(\phi^{(1)})$ our results support a regret of order $\sqrt{T}$ achieved by the restarted OGD, where the multiplicative factor ranges in the interval $[0.05,0.94]$. Under noisy cost access $(\phi^{(0)})$ our results support a regret of order $T^{2/3}$ achieved by the restarted EGS, where the multiplicative factor ranges in the interval $[2.09,2.88]$. In both cases when observations are more noisy, the multiplicative constant increases. While the above policies were introduced mainly as a tool to study the minimax regret rates in the non-stationary stochastic optimization problem, and in particular were not designed to optimize performance in practical settings, we note that in most of the instances that we considered these restarting policies perform at least ``on par" with policies that use fixed step-sizes; while such fixed-step policies are considered as possible heuristics in many practical instances (see, e.g., chapter~4 of \cite{Ben-Pri-Met1990}), they have no performance guarantees relative to the dynamic oracle considered here.

\section{Concluding Remarks}\label{sec:rem}\vspace{-0.2cm}
\paragraph*{\quad On the transition from stationary to non-stationary settings.} Throughout the paper we address ``significant" variation in the cost function, and for the sake of concreteness assume $V_{T}\geq 1$. Nevertheless, one may show (following the proofs of Theorems 2-5) that under each of the different cost and feedback structures, the established bounds hold for ``smaller" variation scales, and if the variation scale is sufficiently ``small," the minimax regret rates coincide with the ones in the classical stationary SA settings. We refer to the variation scales at which the stationary and the non-stationary complexities coincide as ``critical variation scales." Not surprisingly, these transition points between the stationary and the non-stationary regimes differ across cost and feedback structures. The following table summarizes the minimax regret rates for a variation budget of the form $V_{T}=T^{\alpha}$, and documents the critical variation scales in different settings.
\begin{table}[h]
  \centering
    \begin{tabular}{|c||c||c||c||c|c}

\hline
\multicolumn{2}{|c||}{Setting} & \multicolumn{2}{|c||}{Order of regret} &  \multicolumn{1}{|c|}{Critical variation scale}  \\
\hline
\multicolumn{1}{|c||}{Class of functions} &  \multicolumn{1}{|c||}{Feedback} &      \multicolumn{1}{|c||}{Stationary}  &  \multicolumn{1}{|c||}{Non-stationary} &\multicolumn{1}{|c|}{\cellcolor{blue!0} }  \\
\hline
convex  & noisy gradient & $T^{1/2}$ & $\max\left\{T^{1/2}, T^{(2+\alpha)/3}\right\}$ & $T^{-1/2}$\\
strongly convex & noisy gradient & $\log{T}$        &  $\max\left\{\log T, T^{(1+\alpha)/2}\right\}$  & $\left(\log T\right)^{2}T^{-1}$ \\
strongly convex & noisy function & $T^{1/2}$   & $\max\left\{T^{1/2}, T^{(2+\alpha)/3}\right\}$  & $T^{-1/2}$\\
\hline
    \end{tabular}
    \caption{\small \textbf{Critical variation scales.} The growth rates of the minimax regret in different settings for $V_{T}=T^{\alpha}$ (where $\alpha\leq 1$) and the variation scales that separate the stationary and the non-stationary regimes.}
  \label{tab:summary}
\end{table}
In all cases highlighted in the table, the transition point occurs for variation scales that diminish with $T$; this critical quantity therefore measures how ``small" should the temporal variation be, relative to the horizon length, to make non-stationarity effects insignificant relative to other problem primitives insofar as the regret measure goes.

\vspace{-0.3cm}
\paragraph*{\quad Inaccurate or no information on the variation budget.}
The policies introduced in this paper rely on prior knowledge of the variation budget $V_{T}$, but predictions of $V_{T}$ may underestimate or overestimate it. Denoting the ``real" variation budget by $V_{T}$ and the estimate that is used by the agent when tuning the restarting procedure by $\hat{V}_{T}$, one may observe that Proposition~\ref{prop:meta} holds with $V_{T}$ (and the respective class $\mathcal{V}$), but $\Delta_{T}$ is tuned (e.g., in Theorems~2, 4, and~5) using the estimate $\hat{V}_{T}$. This implies that in all the settings that have been considered here, when the ``real" budget is close enough to the estimate $\hat{V}_{T}$, the restarting procedure still guarantees long-run average optimality (naturally, the respective performance is dominated by the one achieved with an accurate knowledge of $V_T$).

For example, consider the case of cost observation with strongly convex cost functions, and suppose that the restarting procedure is tuned by $\hat{V}_T = T^{\alpha}$, but the variation is $V_{T} = T^{\alpha + \delta}$. Then sublinear regret (of order $T^{2/3 + \alpha/3 + \delta}$) can be guaranteed as long as $\delta < 1 - \alpha/3 - 2/3$, and otherwise sublinear regret cannot be guaranteed; e.g., if $\alpha=0$ and $\delta = 1/4$, the restarting procedure may guarantee order $T^{11/12}$ (accurate tuning of the restarting procedure would have guarantee order $T^{3/4}$).

Since there are essentially no restrictions on the rate at which the variation budget can be consumed (in particular, nature is not constrained to sequences with epoch-homogenous variation), an interesting and potentially challenging open problem is to delineate to what extent it is possible to design adaptive policies that do not have a-priori knowledge of the variation budget, yet have performance ``close" to the order of the minimax regret characterized in this paper.Moreover, for both known or unknown variation budgets, characterizing the minimax regret more finely, including the multiplicative constants, remains an important open research avenue of clear practical importance.

\renewcommand{\theequation}{\thesection-\arabic{equation}}
\setcounter{equation}{0}
\renewcommand{\thelemma}{\thesection-\arabic{lemma}}
\setcounter{lemma}{0}

\vspace{-0.2cm}
\appendix
\section{Proofs of main results}\label{app:proofs}
\vspace{-0.2cm}

\paragraph*{Proof of Proposition \ref{prop:lin}.} See Appendix \ref{app:add}.\vspace{-0.4cm}

\paragraph*{Proof of Theorem \ref{thm:meta}.} Fix $\phi \in\left\{ \phi^{(0)},\phi^{(1)} \right\}$, and assume $V_{T}=o(T)$. Let $\mathcal{A}$ be a policy such that $\mathcal{G}^{\mathcal{A}}_{\phi}\left(\mathcal{F},T\right)=~o(T)$, and let $\Delta_{T}\in \left\{1,\ldots,T\right\}$. Let $\pi$ be the policy defined by the restarting procedure that uses $\mathcal{A}$ as a subroutine with batch size $\Delta_{T}$. Then, by Proposition \ref{prop:meta},\vspace{-0.3cm}
\[
\frac{\mathcal{R}^{\pi}_{\phi}(\mathcal{V},T)}{T}
 \;\leq\; \frac{\mathcal{G}^{\mathcal{A}}_{\phi}(\mathcal{F},\Delta_{T})}{\Delta_{T}}
  + \frac{\mathcal{G}^{\mathcal{A}}_{\phi}(\mathcal{F},\Delta_{T})}{T}
  + 2 \Delta_{T} \cdot \frac{V_{T}}{T},\vspace{-0.3cm}
\]
for any $1\leq \Delta_{T} \leq T$. Since $V_{T} = o(T)$, for any selection of $\Delta_{T}$ such that $\Delta_{T}= o(T/V_{T})$ and $\Delta_{T}\rightarrow \infty$ as $T\rightarrow \infty$, the right-hand-side of the above converges to zero as $T\rightarrow \infty$, concluding the  proof.\qed\vspace{-0.2cm}

\paragraph*{Proof of Proposition \ref{prop:meta}.} Fix $\phi\in\left\{\phi^{(0)}, \phi^{(1)}\right\}$, $T\geq1$, and $1\leq V_{T} \leq T$. For $\Delta_{T}\in \left\{1,\ldots,T\right\}$, we break the horizon~$\mathcal{T}$ into a sequence of batches $\mathcal{T}_{1},\ldots,\mathcal{T}_{m}$ of size $\Delta_{T}$ each (except, possibly the last batch) according to (\ref{eq:batches}). Fix $\mathcal{A}\in\mathcal{P}_{\phi}$, and let $\pi$ be the policy defined by the restarting procedure that uses $\mathcal{A}$ as a subroutine with batch size $\Delta_{T}$. Let $f\in\mathcal{V}$. We decompose the regret in the following way: $R^{\pi}(f,T) = \sum_{j=1}^{m}R^{\pi}_{j}$, where\vspace{-0.5cm}
\begin{eqnarray}
R^{\pi}_{j}
& := & \mathbb{E}^{\pi}\left[\sum_{t \in \mathcal{T}_{j}}\left(f_{t}(X_{t}) - f_{t}(x^{\ast}_{t})\right)\right]\notag\\
\label{eq:decompose}
&=& \underbrace{\mathbb{E}^{\pi}\left[\sum_{t \in \mathcal{T}_{j}}f_{t}(X_{t})\right]
 - \min_{x\in \mathcal{X}}\left\{\sum_{t \in \mathcal{T}_{j}}f_{t}(x)\right\}}_{J_{1,j}}
 \;+\; \underbrace{\min_{x\in \mathcal{X}}\left\{\sum_{t \in \mathcal{T}_{j}}f_{t}(x)\right\}
  -  \sum_{t \in \mathcal{T}_{j}}f_{t}(x^{\ast}_{t})}_{J_{2,j}}.
\end{eqnarray}
The first component, $J_{1,j}$, is the regret with respect to the single-best-action of batch $j$, and the second component, $J_{2,j}$, is the difference in performance along batch $j$ between the single-best-action of the batch and the dynamic benchmark. We next analyze $J_{1,j}$, $J_{2,j}$, and the regret throughout the horizon.

\textbf{Step 1} (\textbf{Analysis of $J_{1,j}$}). By taking the sup over all sequences in $\mathcal{F}$ (recall that $\mathcal{V}\subseteq \mathcal{F}$) and using the regret with respect to the single best action in the adversarial setting, one has:\vspace{-0.4cm}
\begin{equation}\label{eq:J1}
J_{1,j}
\;\leq\; \sup_{f\in \mathcal{F}}\left\{\mathbb{E}^{\pi}\left[\sum_{t \in \mathcal{T}_{j}}f_{t}(X_{t})\right] - \min_{x\in \mathcal{X}}\left\{\sum_{t \in \mathcal{T}_{j}}f_{t}(x)\right\}\right\}
\;\leq\; \mathcal{G}_{\phi}^{\mathcal{A}}\left(\mathcal{F}, \Delta_{T}\right),\vspace{-0.3cm}
\end{equation}
where the last inequality holds using (\ref{eq:advRegret}), and since in each batch decisions are dictated by $\mathcal{A}$, and since in each batch there are at most $\Delta_{T}$ epochs (recall that $\mathcal{G}^{\mathcal{A}}_{\phi}$ is non-decreasing in the number of epochs).

\textbf{Step 2} (\textbf{Analysis of $J_{2,j}$}). Defining $f_{0}(x) = f_{1}(x)$, we denote by $V_{j}=\sum_{t\in\mathcal{T}_{j}}\left\|f_{t}-f_{t-1}\right\|$ the variation along batch $\mathcal{T}_{j}$.
By the variation constraint (\ref{eq:TVF}), one has:\vspace{-0.3cm}
\begin{equation}\label{eq:meta}
\sum_{j= 1}^{m}V_{j}
\; =\; \sum_{j= 1}^{m} \sum_{t \in \mathcal{T}_{j}}\sup_{x\in \mathcal{X}}\left|f_{t}(x) - f_{t-1}(x)\right|
\; \leq\; V_{T}.\vspace{-0.3cm}
\end{equation}
Let $\tilde{t}$ be the first epoch of batch $\mathcal{T}_{j}$. Then,\vspace{-0.3cm}
\begin{equation}\label{eq:upp4a}
\min_{x\in \mathcal{X}}\left\{\sum_{t \in \mathcal{T}_{j}}f_{t}(x)\right\}
 - \sum_{t \in \mathcal{T}_{j}}f_{t}(x^{\ast}_{t})
 \;\leq\; \sum_{t \in \mathcal{T}_{j}}\left(f_{t}(x_{\tilde{t}}^{\ast}) - f_{t}(x^{\ast}_{t})\right)
\;\leq\; \Delta_{T} \cdot \max_{t\in \mathcal{T}_{j}}\left\{ f_{t}(x_{\tilde{t}}^{\ast}) - f_{t}(x^{\ast}_{t}) \right\}.\vspace{-0.3cm}
\end{equation}
We next show that $\max_{t\in \mathcal{T}_{j}}\left\{ f_{t}(x_{\tilde{t}}^{\ast}) - f_{t}(x^{\ast}_{t}) \right\} \leq 2V_{j}$. Suppose otherwise. Then, there is some epoch $t_0\in\mathcal{T}_{j}$ at which $f_{t_0}(x_{\tilde{t}}^{\ast}) - f_{t_0}(x^{\ast}_{t_{0}}) > 2V_{j}$, implying\vspace{-0.3cm}
\[
f_{t}(x^{\ast}_{t_{0}})
\;\stackrel{(a)}\leq\; f_{t_{0}}(x^{\ast}_{t_{0}}) + V_{j}
\;<\; f_{t_0}(x_{\tilde{t}}^{\ast}) - V_{j}
\;\stackrel{(b)}\leq\; f_{t}(x_{\tilde{t}}^{\ast}),\quad\quad \text{ for all } t\in \mathcal{T}_{j},\vspace{-0.3cm}
\]
where $(a)$ and $(b)$ follows from the fact that $V_{j}$ is the maximal variation along batch $\mathcal{T}_{j}$. 
In particular, the above holds for $t=\tilde{t}$, contradicting the optimality of $x_{\tilde{t}}^{\ast}$ at epoch $\tilde{t}$. Therefore, one has from (\ref{eq:upp4a}):\vspace{-0.3cm}
\begin{equation}\label{eq:upp4}
\min_{x\in \mathcal{X}}\left\{\sum_{t \in \mathcal{T}_{j}}f_{t}(x)\right\}
 - \sum_{t \in \mathcal{T}_{j}}f_{t}(x^{\ast}_{t})
\;\leq\; 2\Delta_{T} V_{j}.\vspace{-0.3cm}
\end{equation}

\textbf{Step 3} (\textbf{Analysis of the regret over $T$ periods}). Summing (\ref{eq:upp4}) over batches and using (\ref{eq:meta}), one has\vspace{-0.3cm}
\begin{equation}\label{eq:Fbound}
\sum_{j= 1}^{m}\left(\min_{x\in \mathcal{X}}\left\{\sum_{t \in \mathcal{T}_{j}}f_{t}(x)\right\}
 - \sum_{t \in \mathcal{T}_{j}}f_{t}(x^{\ast}_{t})\right)
 \;\leq\; \sum_{j= 1}^{m}2\Delta_{T} V_{j}
 \;\leq\; 2\Delta_{T} V_{T}.\vspace{-0.1cm}
\end{equation}
Therefore, by the regret decomposition in (\ref{eq:decompose}), and following (\ref{eq:J1}) and (\ref{eq:Fbound}), one has:\vspace{-0.3cm}
\[
R^{\pi}\left(f,T\right)
\;\leq\; \sum_{j= 1}^{m}\mathcal{G}^{\mathcal{A}}_{\phi}(\mathcal{F},\Delta_{T})
+ 2 \Delta_{T} V_{T}.\vspace{-0.4cm}
\]
Since the above holds for any $f\in \mathcal{V}$, and recalling that $m = \left\lceil \frac{T}{\Delta_{T}} \right\rceil$, we have\vspace{-0.3cm}
\[
\mathcal{R}^{\pi}_{\phi}(\mathcal{V},T)
\; =\; \sup_{f \in \mathcal{V}}R^{\pi}\left(f,T\right)
\;\leq\; \left\lceil\frac{T}{\Delta_{T}}\right\rceil\cdot\mathcal{G}^{\mathcal{A}}_{\phi}(\mathcal{F},\Delta_{T}) + 2 \Delta_{T} V_{T}.\quad\quad\vspace{-0.3cm}
\]
This concludes the proof. \qed\vspace{-0.2cm}

\paragraph*{Proof of Theorem \ref{thm:weakGrad-low}.} Fix $T\geq 1$ and $1\leq V_{T} \leq T$. We will restrict nature to a specific class of function sequences $\mathcal{V}' \subset \mathcal{V}$. In any element of $\mathcal{V}'$ the cost function is limited to be one of two known quadratic functions, selected by nature in the beginning of every batch of $\tilde{\Delta}_{T}$ epochs, and applied for the following $\tilde{\Delta}_{T}$ epochs. Then we will show that any policy in $\mathcal{P}_{\phi^{(1)}}$ must incur regret of order $V_{T}^{1/3}T^{2/3}$.

\textbf{Step 1 (Preliminaries).} Let $\mathcal{X} = [0,1]$ and consider the following two functions:\vspace{-0.4cm}
\begin{equation}\label{eq:convFuns}
   f^{1}(x) = \left\{
     \begin{array}{lr}
       \frac{1}{2} + \delta - 2\delta x + \left(x-\frac{1}{4}\right)^{2}&  x<\frac{1}{4} \\
       \frac{1}{2} + \delta - 2\delta x &  \frac{1}{4}\leq x \leq\frac{3}{4} \\
       \frac{1}{2} + \delta - 2\delta x + \left(x-\frac{3}{4}\right)^{2}&  x>\frac{3}{4}
     \end{array}
   \right.\;\;;\;\;
   f^{2}(x) = \left\{
     \begin{array}{lr}
     \frac{1}{2} - \delta + 2\delta x + \left(x-\frac{1}{4}\right)^{2}&  x<\frac{1}{4} \\
     \frac{1}{2} - \delta + 2\delta x &  \frac{1}{4} \leq x \leq \frac{3}{4} \\
       \frac{1}{2} - \delta + 2\delta x + \left(x-\frac{3}{4}\right)^{2}&  x>\frac{3}{4},
     \end{array}
   \right.\vspace{-0.1cm}
\end{equation}
for some $\delta>0$ that will be specified shortly. Denoting $x_{k}^{\ast} = \argmin_{x\in[0,1]}f^{k}(x)$, one has $x_{1}^{\ast} = \frac{3}{4} + \delta$, and $x_{2}^{\ast} = \frac{1}{4} - \delta$. It is immediate that 
$f^{1}$ and $f^{2}$ are convex and for any $\delta\in(0,1/4)$ obtain a global minimum in an interior point in $\mathcal{X}$. For some $\tilde{\Delta}_{T}\in\left\{1,\ldots,T\right\}$ that will be specified below, define a partition of the 
horizon $\mathcal{T}$ to $m =~\left\lceil T/\tilde{\Delta}_{T}\right\rceil$ batches $\mathcal{T}_{1},\ldots,\mathcal{T}_{m}$ of size $\tilde{\Delta}_{T}$ each (except perhaps 
$\mathcal{T}_{m}$), according to (\ref{eq:batches}). Define:\vspace{-0.4cm}
\begin{equation}\label{eq:Vprime}
\mathcal{V}'
\;=\; \left\{f\;:\; f_{t}\in \left\{f^{1},f^{2}\right\} \;\text{and}\; f_{t} = f_{t+1}\;\text{for}\; (j-1)\tilde{\Delta}_{T}+1\leq t \leq \min\left\{j\tilde{\Delta}_{T},T\right\}-1,\; j=1,\ldots,m \right\}.\vspace{-0.1cm}
\end{equation}
In every sequence in $\mathcal{V}'$ the cost function is restricted to the set $\left\{f^{1},f^{2}\right\}$, and cannot change throughout a batch. Let $\delta =  V_{T}\tilde{\Delta}_{T}/2T$. Any sequence in $\mathcal{V}'$ consists of convex functions, with minimizers that are interior points in $\mathcal{X}$. In addition, one has:\vspace{-0.3cm}
\[
\sum_{t=2}^{T}\left\|f_{t} - f_{t-1}\right\|
 \;\leq\; \sum_{j=2}^{m}\sup_{x\in \mathcal{X}}\left|f^{1}(x) - f^{2}(x)\right|
 \;=\; \left(\left\lceil\frac{T}{\tilde{\Delta}_{T}}\right\rceil-1 \right)\cdot 2\delta
 \;\leq\; \frac{2T\delta}{\tilde{\Delta}_{T}}
 \;\leq\; V_{T},\vspace{-0.3cm}
 \]
where the first inequality holds since the function can only change between batches. Therefore, $\mathcal{V}'\subset \mathcal{V}$.

\textbf{Step 2 (Bounding the relative entropy within a batch).} Fix any policy $\pi\in\mathcal{P}_{\phi^{(1)}}$. At each $t\in\mathcal{T}_{j}$, the decision maker selects $X_{t}\in \mathcal{X}$ and observes a noisy feedback $\phi^{(1)}_{t}(X_{t},f_{t})$. For any $f\in\mathcal{F}$: denote by $\mathbb{P}_{f}^{\pi}$ the probability measure under policy $\pi$ when $f$ is the sequence of cost functions that is selected by nature, and by $\mathbb{E}_{f}^{\pi}$ the associated expectation operator; For any $\tau \geq 1$, $A\subset\mathbb{R}^{d\times \tau}$ and $B\subset \mathcal{U}$, denote $\mathbb{P}_{f}^{\pi,\tau}(A,B):= \mathbb{P}_{f}^{\pi}\left\{\left\{\phi^{(1)}_{t}(X_{t},f_{t})\right\}_{t=1}^{\tau}\in A, U\in B\right\}$. In what follows we make use of the Kullback-Leibler divergence defined in (\ref{eq:KL}).\vspace{-0.3cm}

\begin{lemma}\label{lem:KL}
\textbf{\textup{(Bound on KL divergence for noisy gradient observations)}} Consider the feedback structure $\phi=\phi^{(1)}$ and let Assumption~1 holds. Then, for any 
$\tau \geq 1$ and 
$f,g\in \mathcal{F}$:\vspace{-0.3cm}
\[
\mathcal{K}\left(\mathbb{P}^{\pi,\tau}_{f}\|\mathbb{P}^{\pi,\tau}_{g}\right)
 \;\leq\; \tilde{C}\mathbb{E}^{\pi}_{f}\left[ \sum_{t=1}^{\tau}\left\|\nabla f_{t}(X_{t}) - \nabla g_{t}(X_{t})\right\|^{2} \right],\vspace{-0.3cm}
\]
where $\tilde{C}$ is the constant that appears in the second part of Assumption 1.
\end{lemma}

The proof of Lemma \ref{lem:KL} appears in Appendix \ref{app:add}. We also use the following result for the minimal error probability in distinguishing between two distributions:\vspace{-0.2cm}
\begin{lemma}\label{lem:tsy}
\textbf{\textup{(Theorem 2.2 in \cite{Tsy})}} Let $\mathbb{P}$ and $\mathbb{Q}$ be two probability distributions on $\mathcal{H}$, such that $\mathcal{K}(\mathbb{P}\|\mathbb{Q})\leq \beta < \infty$. Then, for any $\mathcal{H}$-measurable real function $\varphi:\mathcal{H}\rightarrow \left\{0,1\right\}$,\vspace{-0.3cm}
\[
\max\left\{\mathbb{P}(\varphi = 1), \mathbb{Q}(\varphi = 0)\right\} \;\geq\; \frac{1}{4}\exp\left\{-\beta\right\}.\vspace{-0.2cm}
\]
\end{lemma}
Set $\tilde{\Delta}_{T} = \max\left\{\left\lfloor\left(\frac{1}{4\tilde{C}}\right)^{1/3}\left(\frac{T}{V_{T}}\right)^{2/3}\right\rfloor,1\right\}$, (where $\tilde{C}$ is the constant that appears in part~2 of Assumption~1). We next show that for each batch $\mathcal{T}_{j}$, 
$\mathcal{K}\left(\mathbb{P}^{\pi,\tau}_{f^{1}}\|\mathbb{P}^{\pi,\tau}_{f^{2}}\right)$ is bounded for any $1\leq \tau \leq \left|\mathcal{T}_{j}\right|$. Fix $j\in\left\{1,\ldots,m\right\}$. Then:\vspace{-0.5cm}
\begin{eqnarray}
\mathcal{K}\left(\mathbb{P}^{\pi,|\mathcal{T}_{j}|}_{f^{1}}\|\mathbb{P}^{\pi,|\mathcal{T}_{j}|}_{f^{2}}\right)
&\stackrel{(a)} \leq& \tilde{C}\mathbb{E}^{\pi}_{f^{1}}\left[\sum_{t\in\mathcal{T}_{j}}\left(\nabla f^{1}_{t}(X_{t}) - \nabla f^{2}_{t}(X_{t})\right)^{2}\right]\notag\\
 &=& \tilde{C}\mathbb{E}^{\pi}_{f^{1}}\left[\sum_{t\in\mathcal{T}_{j}}16\delta^{2}X_{t}^{2}\right]
\;\leq\; 16\tilde{C}\tilde{\Delta}_{T}\delta^{2}\notag\\
&\stackrel{(b)} =& \frac{4\tilde{C}V_{T}^{2}\tilde{\Delta}_{T}^{3}}{T^{2}}
\;\stackrel{(c)}\leq\; \max\left\{1,\frac{2\tilde{C}V_{T}}{T}\right\} 
\;\stackrel{(d)}\leq\; \max\left\{1,2\tilde{C}\right\},\notag
\end{eqnarray}
where: $(a)$ follows from Lemma \ref{lem:KL}; $(b)$ and $(c)$ hold given the respective values of $\delta$ and $\tilde{\Delta}_{T}$; and $(d)$ holds by $V_{T}\leq T$. Set $\beta = \max\left\{1,2\tilde{C}\right\}$. Since $\mathcal{K}\left(\mathbb{P}^{\pi,\tau}_{f^{1}}\|\mathbb{P}^{\pi,\tau}_{f^{2}}\right)$ is non-decreasing in $\tau$ throughout a batch, we deduce that $\mathcal{K}\left(\mathbb{P}^{\pi,\tau}_{f^{1}}\|\mathbb{P}^{\pi,\tau}_{f^{2}}\right)$ is bounded by $\beta$ throughout each batch. Then, for any $x_0\in \mathcal{X}$, using Lemma \ref{lem:tsy} with $\varphi_{t} = \mathbbm{1}\{X_{t} \leq x_0\}$, one has:\vspace{-0.3cm}
\begin{equation}\label{eq:phierrorGrad}
\max\left\{\mathbb{P}_{f^{1}}\left\{X_{t} \leq x_0\right\}, \mathbb{P}_{f^{2}}\left\{X_{t} > x_0\right\}\right\}
\;\geq\; \frac{1}{4e^{\beta}} \quad\quad \text{ for all } t\in \mathcal{T}.\vspace{-0.3cm}
\end{equation}
\textbf{Step 3 (A lower bound on the incurred regret for $\boldsymbol{f} \boldsymbol{\in} \boldsymbol{\mathcal{V}'}$).} Set $x_{0} =\frac{1}{2}\left(x_{1}^{\ast}+x_{2}^{\ast}\right) =\frac{1}{2}$. Let $\tilde{f}$ be a random sequence in which in the beginning of each batch $\mathcal{T}_{j}$ a cost function is independently drawn according to a discrete uniform distribution over $\left\{f^{1},f^{2}\right\}$, and applied throughout the whole batch. In particular, note that for any $1\leq j \leq m$, for any epoch $t\in \mathcal{T}_{j}$, $f_{t}$ is independent of $\mathcal{H}_{(j-1)\tilde{\Delta}_{T}+1}$ (the history that is available at the beginning of the batch). Clearly any realization of $\tilde{f}$ is in $\mathcal{V}'$. In particular, taking expectation over $\tilde{f}$, one has:\vspace{-0.3cm}
\begin{eqnarray}
\mathcal{R}^{\pi}_{\phi^{(1)}}\left(\mathcal{V}',T\right)
&\geq& \mathbb{E}^{\pi,\tilde{f}}\left[\sum_{t=1}^{T}\tilde{f}_{t}(X_{t}) - \sum_{t=1}^{T}\tilde{f}_{t}(x_{t}^{\ast})\right]
\;=\; \mathbb{E}^{\pi,\tilde{f}}\left[\sum_{j=1}^{m}\sum_{t\in \mathcal{T}_{j}}\left(\tilde{f}_{t}(X_{t}) - \tilde{f}_{t}(x_{t}^{\ast})\right)\right]\notag\\
 &=& \sum_{j=1}^{m}\left(\frac{1}{2}\cdot\mathbb{E}_{f^{1}}^{\pi}\left[\sum_{t\in \mathcal{T}_{j}}\left(f^{1}(X_{t}) - f^{1}(x_{1}^{\ast})\right)\right]
+ \frac{1}{2}\cdot\mathbb{E}_{f^{2}}^{\pi}\left[\sum_{t\in \mathcal{T}_{j}}\left(f^{2}(X_{t}) - f^{2}(x_{2}^{\ast})\right)\right]\right)\notag\\
&\stackrel{(a)}\geq& \sum_{j=1}^{m}\frac{1}{2}\left(\sum_{t\in \mathcal{T}_{j}}\left(f^{1}(x_{0}) - f^{1}(x_{1}^{\ast})\right)\mathbb{P}_{f^{1}}^{\pi}\left\{X_{t}> x_{0}\right\}
+ \sum_{t\in \mathcal{T}_{j}}\left(f^{2}(x_{0}) - f^{2}(x_{2}^{\ast})\right)\mathbb{P}_{f^{2}}^{\pi}\left\{X_{t}\leq x_{0}\right\}  \right)\notag\\
&\geq& \sum_{j=1}^{m}\frac{\delta}{4} \sum_{t\in \mathcal{T}_{j}}\left(\mathbb{P}_{f^{1}}^{\pi}\left\{X_{t}> x_{0}\right\} + \mathbb{P}_{f^{2}}^{\pi}\left\{X_{t}\leq x_{0}\right\}  \right)\notag\\
&\geq& \sum_{j=1}^{m}\frac{\delta}{4} \sum_{t\in \mathcal{T}_{j}}\max\left\{\mathbb{P}^{\pi}_{f^{1}}\left\{X_{t} > x_{0}\right\}, \mathbb{P}^{\pi}_{f^{2}}\left\{X_{t} \leq x_{0}\right\}\right\}\notag\\
&\stackrel{(b)}\geq& \sum_{j=1}^{m} \frac{\delta}{4} \sum_{t\in \mathcal{T}_{j}}\frac{1}{4e^{\beta}}
\; = \; \sum_{j=1}^{m}\frac{\delta\tilde{\Delta}_{T}}{16e^{\beta}}\notag\\
&\stackrel{(c)} =& \sum_{j=1}^{m}\frac{V_{T}\tilde{\Delta}_{T}^{2}}{32e^{\beta}T}
\;\geq \frac{T}{\tilde{\Delta}_{T}}\cdot\frac{V_{T}\tilde{\Delta}_{T}^{2}}{32e^{\beta}T}
\;=\; \frac{V_{T}\tilde{\Delta}_{T}}{32e^{\beta}},\notag\vspace{-0.1cm}
\end{eqnarray}
where $(a)$ holds since for any function $g:[0,1]\rightarrow\mathbb{R}^{+}$ and $x_{0}\in[0,1]$ such that $g(x) \geq g(x_{0})$ for all $x>x_{0}$, one has that
$\mathbb{E}\left[g(X_{t})\right]
=\mathbb{E}\left[g(X_{t})|X_{t}>x_{0}\right]\mathbb{P}\left\{X_{t}>x_{0}\right\}
+ \mathbb{E}\left[g(X_{t})|X_{t}\leq x_{0}\right]\mathbb{P}\left\{X_{t}\leq x_{0}\right\}
\geq g(x_{0})\mathbb{P}\left\{X_{t}>x_{0}\right\}$
for any $t\in\mathcal{T}$, and similarly for any $x_{0}\in[0,1]$ such that $g(x) \geq g(x_{0})$ for all $x\leq x_{0}$, one obtains $\mathbb{E}\left[g(X_{t})\right] \geq g(x_{0})\mathbb{P}\left\{X_{t}\leq x_{0}\right\}$. In addition, $(b)$ holds by (\ref{eq:phierrorGrad}) and $(c)$ holds by $\delta = V_{T}\tilde{\Delta}_{T}/2T$. Suppose that $T\geq 2^{5/2}\sqrt{\tilde{C}}\cdot V_{T}$. Applying the selected $\tilde{\Delta}_{T}$, one has: \vspace{-0.1cm}
\begin{eqnarray}
\mathcal{R}^{\pi}_{\phi^{(1)}}\left(\mathcal{V}',T\right)
&\geq&\frac{V_{T}}{32e^{\beta}}\cdot\left\lfloor\left(\frac{1}{4\tilde{C}}\right)^{1/3}\left(\frac{T}{V_{T}}\right)^{2/3}\right\rfloor\notag\\
&\geq&\frac{V_{T}}{32e^{\beta}}\cdot\left(\left(\frac{1}{4\tilde{C}}\right)^{1/3}\left(\frac{T}{V_{T}}\right)^{2/3} - 1\right)\notag\\
&=&\frac{V_{T}}{32e^{\beta}}\cdot\left( \frac{T^{2/3} - \left(4\tilde{C}\right)^{1/3}V_{T}^{2/3}}{\left(4\tilde{C}\right)^{1/3}V_{T}^{2/3}}\right)
\;\geq\; \frac{1}{64e^{\beta}\left(4\tilde{C}\right)^{1/3}}\cdot V_{T}^{1/3} T^{2/3},\notag
\end{eqnarray}
where the last inequality follows from $T\geq 2^{5/2}\sqrt{\tilde{C}}\cdot V_{T}$. If $T< 2^{5/2}\sqrt{\tilde{C}}\cdot V_{T}$, by Proposition~\ref{prop:lin} there exists a constant $C$ such that $\mathcal{R}^{\pi}_{\phi^{(1)}}\left(\mathcal{V},T\right)\geq~C\cdot T \geq~C\cdot V_{T}^{1/3}T^{2/3}$. Recalling that $\mathcal{V}'\subseteq \mathcal{V}$, we have:\vspace{-0.2cm}
\[
\mathcal{R}^{\pi}_{\phi^{(1)}}\left(\mathcal{V},T\right)
\;\geq\;\mathcal{R}^{\pi}_{\phi^{(1)}}\left(\mathcal{V}',T\right)
\;\geq\;\frac{1}{64e^{\beta}\left(4\tilde{C}\right)^{1/3}}\cdot V_{T}^{1/3} T^{2/3}.\vspace{-0.3cm}
\]
This concludes the proof.\qed\vspace{-0.2cm}

\paragraph*{Proof of Theorem \ref{thm:weakGrad-up}.} Fix $T\geq 1$, and $1\leq V_{T}\leq T$. For any $\Delta_{T}\in \left\{1,\ldots,T\right\}$, let $\mathcal{A}$ be the OGD algorithm with $\eta_{t} = \eta = \frac{r}{G\sqrt{\Delta_{T}}}$ for any $t=2,\ldots,\Delta_{T}$ (where $r$ denotes the radius of the action set $\mathcal{X}$), and let $\pi$ be the policy defined by the restarting procedure with subroutine $\mathcal{A}$ and batch size $\Delta_{T}$. \cite{FLA2005} consider the performance of the OGD algorithm (with a specific, deterministic $x_{1}=0$) relative to the single best action in the adversarial setting, and show (\citealt{FLA2005}, Lemma 3.1) that $\mathcal{G}^{\mathcal{A}}_{\phi^{(1)}}(\mathcal{F},\Delta_{T}) \leq rG\sqrt{\Delta_{T}}$. Following their analysis one obtains that for an arbitrary (potentially random) initial action $\mathcal{G}^{\mathcal{A}}_{\phi^{(1)}}(\mathcal{F},\Delta_{T}) \leq 2rG\sqrt{\Delta_{T}}$\label{pg:adjusttorandomx1}. Therefore, by Proposition \ref{prop:meta},
\vspace{-0.3cm}
\[
\mathcal{R}^{\pi}_{\phi^{(1)}}(\mathcal{V},T)
\;\leq \; \left(\frac{T}{\Delta_{T}}+1\right)\cdot \mathcal{G}^{\mathcal{A}}_{\phi^{(1)}}(\mathcal{F},\Delta_{T})
 + 2 V_{T}\Delta_{T}
\;\leq \;\frac{2rG\cdot  T}{\sqrt{\Delta_{T}}}
 + 2rG\sqrt{\Delta_{T}}
 + 2V_{T}\Delta_{T}.\vspace{-0.3cm}
\]
Selecting $\Delta_{T} =  \left\lceil\left(T / V_{T}\right)^{2/3}\right\rceil$, one has\vspace{-0.3cm}
\begin{eqnarray}
\mathcal{R}^{\pi}_{\phi^{(1)}}(\mathcal{V},T)
&\leq& \frac{2rG\cdot  T}{\left(T/V_{T}\right)^{1/3}}
 + 2rG\left(\left(\frac{T}{V_{T}}\right)^{1/3}+ 1\right)
 + 2V_{T}\left(\left(\frac{T}{V_{T}}\right)^{2/3} + 1\right)\notag\\
&\stackrel{(a)} \leq& \left(2rG+4\right)\cdot V_{T}^{1/3}T^{2/3}
 + 2rG\cdot\left(\frac{T}{V_{T}}\right)^{1/3} + 2rG \notag\\
&\stackrel{(b)}\leq& \left(6rG+4\right)\cdot V_{T}^{1/3}T^{2/3},
\end{eqnarray}
where $(a)$ and $(b)$ follows since $1\leq V_{T}\leq T$. This concludes the proof. \qed\vspace{-0.2cm}

\paragraph*{Proof of Theorem \ref{thm:strictGrad}.} \textbf{Part 1.} We begin with the first part of the Theorem. Fix $T\geq 1$, and $1\leq V_{T} \leq T$. For any $\Delta_{T}\in \left\{1,\ldots,T\right\}$ let $\mathcal{A}$ be the OGD algorithm with $\eta_{t} = 1/Ht$ for any $t=2,\ldots,\Delta_{T}$, and let $\pi$ be the policy defined by the restarting procedure with subroutine $\mathcal{A}$ and batch size $\Delta_{T}$. By Lemma \ref{lem:sbaGrad} (see Appendix \ref{app:sba}), one has:\vspace{-0.4cm}
\begin{equation}\label{eq:Ggrad}
\mathcal{G}^{\mathcal{A}}_{\phi^{(1)}}(\mathcal{F}_{s},\Delta_{T})
\;\leq\; \frac{\left(G^{2} + \sigma^{2}\right)}{2H}\left(1 + \log \Delta_{T}\right).\vspace{-0.3cm}
\end{equation}
 Therefore, by Proposition \ref{prop:meta},\vspace{-0.3cm}
\[
\mathcal{R}^{\pi}_{\phi^{(1)}}(\mathcal{V}_{s},T)
\;\leq \; \left(\frac{T}{\Delta_{T}}+1\right) \cdot \mathcal{G}^{\mathcal{A}}_{\phi^{(1)}}(\mathcal{F}_{s},\Delta_{T}) + 2 V_{T}\Delta_{T}
\;\leq \; \left(\frac{T}{\Delta_{T}}+ 1\right)\frac{\left(G^{2}+ \sigma^{2}\right)}{2H}\left(1 + \log \Delta_{T}\right)
 + 2V_{T}\Delta_{T}.\vspace{-0.3cm}
\]
Selecting $\Delta_{T} =  \left\lceil\sqrt{T / V_{T}}\right\rceil$, one has:\vspace{-0.4cm}
\begin{eqnarray}
\mathcal{R}^{\pi}_{\phi^{(1)}}(\mathcal{V}_{s},T)
&\leq& \left(\frac{T}{\sqrt{T/V_{T}}}+ 1\right)\frac{\left(G^{2}+ \sigma^{2}\right)}{2H}\left(1 + \log \left(\sqrt{\frac{T}{V_{T}}}+1\right)\right)
 + 2V_{T}\left(\sqrt{\frac{T}{V_{T}}}+1\right)\notag\\
&\stackrel{(a)}\leq& \left(4 + \frac{\left(G^{2}+ \sigma^{2}\right)}{2H}\left(1 + \log \left(\sqrt{\frac{T}{V_{T}}}+1\right)\right)\right)\cdot\sqrt{V_{T}T}
 + \frac{\left(G^{2}+ \sigma^{2}\right)}{2H}\left(1 + \log \left(\sqrt{\frac{T}{V_{T}}}+1\right)\right)\notag\\
&\stackrel{(b)}\leq& \left(4 + \frac{\left(2G^{2}+ 2\sigma^{2}\right)}{H} \right)\cdot\log\left(\sqrt{\frac{T}{V_{T}}}+1\right)\sqrt{V_{T}T},\notag
\end{eqnarray}
where $(a)$ and $(b)$ hold since $1\leq V_{T} \leq T$.\vspace{-0.3cm}

\paragraph*{Part 2.} We next prove the second part of the Theorem. The proof follows steps and notation appearing in the proof of Theorem \ref{thm:weakGrad-low}. For strongly convex cost functions a different choice of $\delta$ is used in step 2 and $\tilde{\Delta}_{T}$ is modified accordingly in step 3. The regret analysis in step 4 is adjusted as well.

\textbf{Step 1.} Let $\mathcal{X} = [0,1]$, and consider the following two quadratic functions:\vspace{-0.3cm}
\begin{equation}\label{eq:sconvFuns}
   f^{1}(x) = x^{2} - x + \frac{3}{4},
    \quad\quad\quad\quad
   f^{2}(x) = x^{2} - \left(1+\delta\right)x + \frac{3}{4} + \frac{\delta}{2}\vspace{-0.3cm}
\end{equation}
for some small $\delta>0$. Note that $x_{1}^{\ast} = \frac{1}{2}$, and $x_{2}^{\ast} = \frac{1+\delta}{2}$. We define a partition of $\mathcal{T}$ into batches $\mathcal{T}_{1},\ldots,\mathcal{T}_{m}$ of size $\tilde{\Delta}_{T}$ each (perhaps except $\mathcal{T}_{m}$), according to (\ref{eq:batches}), where $\tilde{\Delta}_{T}$ will be specified below. Define the class $\mathcal{V}'_{s}$ according to (\ref{eq:Vprime}), such that in every $f\in\mathcal{V}'_{s}$ the cost function is restricted to the set $\left\{f^{1},f^{2}\right\}$, and cannot change throughout a batch. The sequences in $\mathcal{V}'_{s}$ consist of strongly convex functions ((\ref{eq:H}) holds for any $H \leq 1$), with minimizers that are interior points in $\mathcal{X}$. Set $\delta =  \sqrt{2V_{T}\tilde{\Delta}_{T}/T}$. Then:\vspace{-0.3cm}
\[
\sum_{t=2}^{T}\sup_{x\in \mathcal{X}^{\ast}}\left|f_{t}(x) - f_{t-1}(x)\right|
\;\leq\; \sum_{j=2}^{m}\sup_{x\in \mathcal{X}^{\ast}}\left|f^{1}(x) - f^{2}(x)\right|
 \;\leq\; \frac{T}{\tilde{\Delta}_{T}}\cdot \frac{\delta^{2}}{2}
 \;=\; V_{T},\vspace{-0.3cm}
 \]
where the first inequality holds since the function can change only between batches. Therefore, $\mathcal{V}_{s}'\subset \mathcal{V}_{s}$.

\textbf{Step 2.} Fix $\pi\in\mathcal{P}_{\phi^{(1)}}$, and let $\tilde{\Delta}_{T} = \max\left\{\left\lfloor\frac{1}{\sqrt{2\tilde{C}}}\cdot \sqrt{\frac{T}{V_{T}}}\right\rfloor,1\right\}$ ($\tilde{C}$ appears in part~2 of Assumption~1). Fix $j\in\left\{1,\ldots,m\right\}$. Then:\vspace{-0.4cm}
\begin{eqnarray}
\mathcal{K}\left(\mathbb{P}^{\pi,|\mathcal{T}_{j}|}_{f^{1}}\|\mathbb{P}^{\pi,|\mathcal{T}_{j}|}_{f^{2}}\right)
&\stackrel{(a)} \leq& \tilde{C}\mathbb{E}^{\pi}_{f^{1}}\left[\sum_{t\in\mathcal{T}_{j}}\left(\nabla f^{1}(X_{t}) - \nabla f^{2}(X_{t})\right)^{2}\right]\notag\\
&\leq& \tilde{C}\tilde{\Delta}_{T}\delta^{2}
\;\stackrel{(b)} =\; \frac{2\tilde{C}V_{T}\tilde{\Delta}_{T}^{2}}{T}\notag\\
\label{eq:boundOnKLGrad}
&\stackrel{(c)}\leq& \max\left\{1,\frac{2\tilde{C}V_{T}}{T}\right\}
\;\stackrel{(d)}\leq\; \max\left\{1,2\tilde{C}\right\},
\end{eqnarray}
where: $(a)$ follows from Lemma \ref{lem:KL}; $(b)$ and $(c)$ hold by the selected values of $\delta$ and $\tilde{\Delta}_{T}$ respectively; and $(d)$ holds by $V_{T}\leq T$. Set $\beta = \max\left\{1,2\tilde{C}\right\}$. Then, for any $x_0\in \mathcal{X}$, using Lemma \ref{lem:tsy} with $\varphi_{t} = \mathbbm{1}\{X_{t} > x_0\}$, one has:\vspace{-0.3cm}
\begin{equation}\label{eq:phierrorGrad2}
\max\left\{\mathbb{P}_{f^{1}}\left\{X_{t} > x_0\right\}, \mathbb{P}_{f^{2}}\left\{X_{t} \leq x_0\right\}\right\} \geq \frac{1}{4e^{\beta}}, \quad\quad \forall t\in \mathcal{T}.\vspace{-0.3cm}
\end{equation}
\textbf{Step 3.} Set $x_{0} = \frac{1}{2}\left(x_{1}^{\ast}+x_{2}^{\ast}\right) = 1/2 + \delta/4$. Let $\tilde{f}$ be a random sequence in which in the beginning of each batch $\mathcal{T}_{j}$ a cost function is independently drawn according to a discrete uniform distribution over $\left\{f^{1},f^{2}\right\}$, and applied throughout the batch. Taking expectation over $\tilde{f}$ one has:\vspace{-0.4cm}
\begin{eqnarray}
\mathcal{R}^{\pi}_{\phi^{(1)}}\left(\mathcal{V}'_{s},T\right)
&\geq&  \sum_{j=1}^{m}\left(\frac{1}{2}\cdot\mathbb{E}_{f^{1}}^{\pi}\left[\sum_{t\in \mathcal{T}_{j}}\left(f^{1}(X_{t}) - f^{1}(x_{1}^{\ast})\right)\right]
+ \frac{1}{2}\cdot\mathbb{E}_{f^{2}}^{\pi}\left[\sum_{t\in \mathcal{T}_{j}}\left(f^{2}(X_{t}) - f^{2}(x_{2}^{\ast})\right)\right]\right)\notag\\
 &\geq& \sum_{j=1}^{m}\frac{1}{2}\left(\sum_{t\in \mathcal{T}_{j}}\left(f^{1}(x_{0}) - f^{1}(x_{1}^{\ast})\right)\mathbb{P}_{f^{1}}^{\pi}\left\{X_{t}> x_{0}\right\}
+ \sum_{t\in \mathcal{T}_{j}}\left(f^{2}(x_{0}) - f^{2}(x_{2}^{\ast})\right)\mathbb{P}_{f^{2}}^{\pi}\left\{X_{t}\leq x_{0}\right\}  \right)\notag\\
&\geq& \sum_{j=1}^{m}\frac{\delta^{2}}{16} \sum_{t\in \mathcal{T}_{j}}\left(\mathbb{P}_{f^{1}}^{\pi}\left\{X_{t}> x_{0}\right\} + \mathbb{P}_{f^{2}}^{\pi}\left\{X_{t}\leq x_{0}\right\}  \right)\notag\\
&\geq& \sum_{j=1}^{m}\frac{\delta^{2}}{16} \sum_{t\in \mathcal{T}_{j}}\max\left\{\mathbb{P}^{\pi}_{f^{1}}\left\{X_{t} > x_{0}\right\}, \mathbb{P}^{\pi}_{f^{2}}\left\{X_{t} \leq x_{0}\right\}\right\}\notag\\
&\;\stackrel{(a)}\geq&  \sum_{j=1}^{m}\frac{\delta^{2}}{16}\sum_{t\in \mathcal{T}_{j}}\frac{1}{4e^{\beta}}
\; = \;  \sum_{j=1}^{m}\frac{\delta^{2}\tilde{\Delta}_{T}}{64e^{\beta}}
\;\stackrel{(b)} =\; \sum_{j=1}^{m}\frac{V_{T}\tilde{\Delta}_{T}^{2}}{32e^{\beta}T}
\;\geq\; \frac{V_{T}\tilde{\Delta}_{T}}{32e^{\beta}},\notag
\end{eqnarray}
where: the first four inequalities follow from arguments given in step 3 in the proof of Theorem \ref{thm:weakGrad-low}; $(a)$ holds by (\ref{eq:phierrorGrad2}); and $(b)$ holds by $\delta = \sqrt{2V_{T}\tilde{\Delta}_{T}/T}$. Given the selection of $\tilde{\Delta}_{T}$, one has:\vspace{-0.3cm}
\[
\mathcal{R}^{\pi}_{\phi^{(1)}}\left(\mathcal{V}'_{s},T\right)
\;\geq\;\frac{V_{T}}{32e^{\beta}}\cdot\left\lfloor\frac{1}{\sqrt{2\tilde{C}}}\cdot \sqrt{\frac{T}{V_{T}}}\right\rfloor
\;\geq\; \frac{V_{T}}{32e^{\beta}}\cdot\left( \frac{\sqrt{T} - \sqrt{2\tilde{C}V_{T}}}{\sqrt{2\tilde{C}V_{T}}}   \right)
\geq\; \frac{1}{64e^{\beta}\sqrt{2\tilde{C}}}\cdot \sqrt{V_{T}T},\vspace{-0.3cm}
\]
where the last inequality holds if $T\geq 8\tilde{C}V_{T}$. If $T < 8\tilde{C} V_{T}$, by Proposition~\ref{prop:lin} there exists a constant $C$ such that $\mathcal{R}^{\pi}_{\phi^{(1)}}\left(\mathcal{V}_{s},T\right)\geq C T \geq C \sqrt{V_{T}T}$. Then, recalling that $\mathcal{V}'_{s}\subseteq \mathcal{V}_{s}$, we have established that\vspace{-0.4cm}
\[
\mathcal{R}^{\pi}_{\phi^{(1)}}\left(\mathcal{V}_{s},T\right)
\;\geq\;\mathcal{R}^{\pi}_{\phi^{(1)}}\left(\mathcal{V}'_{s},T\right)
\;\geq\;\frac{1}{64e^{\beta}\sqrt{2\tilde{C}}}\cdot \sqrt{V_{T}T}.\vspace{-0.4cm}
\]
This concludes the proof.\qed\vspace{-0.2cm}

\paragraph*{Proof of Theorem \ref{thm:strictCost}.} \textbf{Part 1.} Fix $T\geq 1$, and $1 \leq V_{T} \leq T$. For any $\Delta_{T}\in \left\{1,\ldots,T\right\}$, consider the EGS algorithm $\mathcal{A}$ given in \S5.2 with $a_{t} = 2/Ht$ and $\delta_{t} = h_{t} = a_{t}^{1/4}$ for $t=1,\ldots,\Delta_{T}$, and let $\pi$ be the policy defined by the restarting procedure with subroutine $\mathcal{A}$ and batch size $\Delta_{T}$. By Lemma \ref{lem:sbaRegret} (see Appendix \ref{app:sba}), we have: \vspace{-0.4cm}
\begin{equation}\label{eq:Gcost}
\mathcal{G}^{\mathcal{A}}_{\phi}(\mathcal{F}_{s},\Delta_{T}) \leq C_{1}\cdot \sqrt{\Delta_{T}},\vspace{-0.1cm}
\end{equation}
with $C_{1} = 2G +\left(G^{2}+\sigma^{2} + H\right)d^{3/2}/\sqrt{2H}$. Therefore, by Proposition~\ref{prop:meta},\vspace{-0.3cm}
\[
\mathcal{R}^{\pi}_{\phi^{(0)}}(\mathcal{V}_{s},T)
\;\leq \; \left(\frac{T}{\Delta_{T}}+1\right) \cdot \mathcal{G}^{\mathcal{A}}_{\phi^{(0)}}(\mathcal{F},\Delta_{T}) + 2 V_{T}\Delta_{T}
\;\stackrel{(a)} \leq\;  C_{1}\cdot \frac{T}{\sqrt{\Delta_{T}}}
 + C_{1}\cdot \sqrt{\Delta_{T}}
 + 2 V_{T}\Delta_{T},\vspace{-0.3cm}
\]
where $(a)$ holds by (\ref{eq:Gcost}). By selecting $\Delta_{T} =~\left\lceil\left(T / V_{T}\right)^{2/3}\right\rceil$, one obtains\vspace{-0.5cm}
\begin{eqnarray}
\mathcal{R}^{\pi}_{\phi^{(0)}}(\mathcal{V}_{s},T)
&\leq& C_{1}\cdot \frac{T}{\left(T/V_{T}\right)^{1/3}}
 + C_{1}\cdot \left(\left(\frac{T}{V_{T}}\right)^{1/3}+1\right)
 + 2 V_{T}\left(\left(\frac{T}{V_{T}}\right)^{2/3}+1\right)\notag\\
&\stackrel{(b)} \leq& \left(C_{1}+4\right)V_{T}^{1/3}T^{2/3}
 + C_{1}\cdot \left(\frac{T}{V_{T}}\right)^{1/3}
 + C_{1}\notag\\
&\stackrel{(c)} \leq& \left(3C_{1} + 4\right)V_{T}^{1/3}T^{2/3},\notag
\end{eqnarray}
where $(b)$ and $(c)$ hold since $1 \leq V_{T} \leq T$.\vspace{-0.2cm}

\paragraph*{Part 2.} The proof of this part of the theorem follows the steps and uses notation introduced in the proof of Theorem~\ref{thm:weakGrad-low}. The different feedback structure affects the bound on the KL divergence and the selected value of $\tilde{\Delta}_{T}$ in step 2 as well as the resulting regret analysis in step 3. Details are given below.

\textbf{Step 1.} We define a class $\mathcal{V}_{s}'$ as it is defined in the proof of Theorem \ref{thm:strictGrad}, using the quadratic functions $f^{1}$ and $f^{2}$ that are given in (\ref{eq:sconvFuns}), and the partition of $\mathcal{T}$ to batches in (\ref{eq:batches}). Again, selecting $\delta = \sqrt{2V_{T}\tilde{\Delta}_{T}/T}$, we have $\mathcal{V}_{s}'\subset \mathcal{V}_{s}$.

\textbf{Step 2.} Fix some policy $\pi\in\mathcal{P}_{\phi^{(0)}}$. At each $t\in\mathcal{T}_{j}$, $j=1,\ldots,m$, the decision maker selects $X_{t}\in \mathcal{X}$ and observes a noisy feedback $\phi^{(0)}_{t}(X_{t},f^{k})$. For any $f\in \mathcal{F}$, $\tau \geq 1$, $A\subset\mathbb{R}^{\tau}$ and $B\subset \mathcal{U}$, denote $\mathbb{P}_{f}^{\pi,\tau}(A,B):= \mathbb{P}_{f}\left\{\left\{\phi^{(0)}_{t}(X_{t},f_{t})\right\}_{t=1}^{\tau}\in A, U\in B\right\}$. In this part of the proof we use the following counterpart of Lemma~\ref{lem:KL} for the case of noisy cost feedback structure.\vspace{-0.3cm}

\begin{lemma}\label{lem:KL2}
\textbf{\textup{(Bound on KL divergence for noisy cost observations)}} Consider the feedback structure $\phi=\phi^{(0)}$ and let Assumption~2 holds. Then, for any $\tau \geq 1$ and $f,g\in \mathcal{F}$:\vspace{-0.3cm}
\[
\mathcal{K}\left(\mathbb{P}^{\pi,\tau}_{f}\|\mathbb{P}^{\pi,\tau}_{g}\right)
 \;\leq\; \tilde{C}\mathbb{E}^{\pi}_{f}\left[ \sum_{t=1}^{\tau}\left(f_{t}(X_{t}) - g_{t}(X_{t})\right)^{2} \right],\vspace{-0.3cm}
\]
where $\tilde{C}$ is the constant that appears in the second part of Assumption 2.\vspace{-0.1cm}
\end{lemma}

The proof of the lemma appears in Appendix \ref{app:add}. We next bound $\mathcal{K}\left(\mathbb{P}^{\pi,|\mathcal{T}_{j}|}_{f^{1}}\|\mathbb{P}^{\pi,|\mathcal{T}_{j}|}_{f^{2}}\right)$ throughout an arbitrary batch $\mathcal{T}_{j}$, $j\in\left\{1,\ldots,m\right\}$, for a given batch size $\tilde{\Delta}_{T}$. Define:\vspace{-0.4cm}
\[
R_{j}^{\pi} \;=\; \frac{1}{2}\mathbb{E}_{f^{1}}^{\pi}\left[\sum_{t\in \mathcal{T}_{j}}\left(f^{1}(X_{t}) - f^{1}(x_{1}^{\ast})\right)\right]
+ \frac{1}{2}\mathbb{E}_{f^{2}}^{\pi}\left[\sum_{t\in \mathcal{T}_{j}}\left(f^{2}(X_{t}) - f^{2}(x_{2}^{\ast})\right)\right].\vspace{-0.4cm}
\]
Then, one has:\vspace{-0.5cm}
\begin{eqnarray}
\mathcal{K}\left(\mathbb{P}^{\pi,|\mathcal{T}_{j}|}_{f^{1}}\|\mathbb{P}^{\pi,|\mathcal{T}_{j}|}_{f^{2}}\right)
&\stackrel{(a)} \leq& \tilde{C}\mathbb{E}^{\pi}_{f^{1}}\left[\sum_{t\in\mathcal{T}_{j}}\left(f^{1}(X_{t}) - f^{2}(X_{t})\right)^{2}\right]
\;=\; \tilde{C}\mathbb{E}^{\pi}_{f^{1}}\left[\sum_{t\in\mathcal{T}_{j}}\left(\delta X_{t} - \frac{\delta}{2}\right)^{2}\right]\notag\\
&=& \tilde{C}\mathbb{E}^{\pi}_{f^{1}}\left[\delta^{2}\sum_{t\in\mathcal{T}_{j}}\left(X_{t} - x_{1}^{\ast}\right)^{2}\right]
\;\stackrel{(b)} =\; 2\tilde{C}\delta^{2}\mathbb{E}^{\pi}_{f^{1}}\left[\sum_{t\in\mathcal{T}_{j}}\left(f^{1}(X_{t}) - f^{1}(x_{1}^{\ast})\right)\right]\notag\\
\label{eq:boundOnKL}
&\stackrel{(c)} \leq& \frac{8\tilde{C}\tilde{\Delta}_{T}V_{T}}{T}\cdot R_{j}^{\pi}
\end{eqnarray}
where: $(a)$ follows from Lemma \ref{lem:KL2}; $(b)$ holds since\vspace{-0.3cm}
\[
f^{1}(x) - f^{1}(x_{1}^{\ast}) = \nabla f^{1}(x_{1}^{\ast})(x-x_{1}^{\ast}) + \frac{1}{2}\cdot \nabla f^{1}(x_{1}^{\ast})(x-x_{1}^{\ast})^{2} = \frac{1}{2}(x-x_{1}^{\ast})^{2},\vspace{-0.3cm}
\]
for any $x\in \mathcal{X}$; and $(c)$ holds since $\delta = \sqrt{2V_{T}\tilde{\Delta}_{T}/T}$, and $R_{j}^{\pi} \geq \frac{1}{2}\mathbb{E}_{f^{1}}^{\pi}\left[\sum_{t\in \mathcal{T}_{j}}\left(f^{1}(X_{t}) - f^{1}(x_{1}^{\ast})\right)\right]$. Thus, for any $x_{0}\in \mathcal{X}$, using Lemma \ref{lem:tsy} with $\varphi_{t} = \mathbbm{1}\{X_{t} > x_{0}\}$, we have:\vspace{-0.2cm}
\begin{equation}\label{eq:phierror}
\max\left\{\mathbb{P}^{\pi}_{f^{1}}\left\{X_{t} > x_{0}\right\}, \mathbb{P}^{\pi}_{f^{2}}\left\{X_{t} \leq x_{0}\right\}\right\}
 \geq \frac{1}{4}\exp\left\{-\frac{8\tilde{C}\tilde{\Delta}_{T}V_{T}}{T}\cdot R_{j}^{\pi}\right\}, \quad\quad \text{ for all } t\in \mathcal{T}_{j},\quad 1\leq j \leq m.\vspace{-0.2cm}
\end{equation}

\textbf{Step 3.} Set $x_{0} = \frac{1}{2}\left(x_{1}^{\ast}+x_{2}^{\ast}\right) = 1/2 + \delta/4$. Let $\tilde{f}$ be the random sequence of functions that is described in Step 3 in the proof of Theorem \ref{thm:strictGrad}. Taking expectation over $\tilde{f}$, one has:\vspace{-0.3cm}
\[
\mathcal{R}^{\pi}_{\phi^{(0)}}\left(\mathcal{V}'_{s},T\right)
\;\geq\;\sum_{j=1}^{m}\left(\frac{1}{2}\cdot\mathbb{E}_{f^{1}}^{\pi}\left[\sum_{t\in \mathcal{T}_{j}}\left(f^{1}(X_{t}) - f^{1}(x_{1}^{\ast})\right)\right]
+ \frac{1}{2}\cdot\mathbb{E}_{f^{2}}^{\pi}\left[\sum_{t\in \mathcal{T}_{j}}\left(f^{2}(X_{t}) - f^{2}(x_{2}^{\ast})\right)\right]\right)
\;=:\; \sum_{j=1}^{m}R_{j}^{\pi}.\vspace{-0.3cm}
\]
In addition, for each $1\leq j \leq m$ one has:\vspace{-0.5cm}
\begin{eqnarray}
R_{j}^{\pi}
&\geq& \frac{1}{2}\left(\sum_{t\in \mathcal{T}_{j}}\left(f^{1}(x_{0}) - f^{1}(x_{1}^{\ast})\right)\mathbb{P}_{f^{1}}^{\pi}\left\{X_{t}> x_{0}\right\}
+ \sum_{t\in \mathcal{T}_{j}}\left(f^{2}(x_{0}) - f^{2}(x_{2}^{\ast})\right)\mathbb{P}_{f^{2}}^{\pi}\left\{X_{t}\leq x_{0}\right\}  \right)\notag\\
&\geq& \frac{\delta^{2}}{16} \sum_{t\in \mathcal{T}_{j}}\left(\mathbb{P}_{f^{1}}^{\pi}\left\{X_{t}> x_{0}\right\} + \mathbb{P}_{f^{2}}^{\pi}\left\{X_{t}\leq x_{0}\right\}  \right)\notag\\
&\geq& \frac{\delta^{2}}{16} \sum_{t\in \mathcal{T}_{j}}\max\left\{\mathbb{P}^{\pi}_{f^{1}}\left\{X_{t} > x_{0}\right\}, \mathbb{P}^{\pi}_{f^{2}}\left\{X_{t} \leq x_{0}\right\}\right\}\notag \quad\quad\quad\quad\quad\quad\quad\quad\quad\quad\quad\quad\quad\quad\quad  \\
&\stackrel{(a)} \geq& \frac{\delta^{2}}{16}\sum_{t\in \mathcal{T}_{j}}\frac{1}{4}\exp\left\{-\frac{8\tilde{C}\tilde{\Delta}_{T}V_{T}}{T}\cdot R_{j}^{\pi}\right\}
 \;=\; \frac{\delta^{2}\tilde{\Delta}_{T}}{64}\exp\left\{-\frac{8\tilde{C}\tilde{\Delta}_{T}V_{T}}{T}\cdot R_{j}^{\pi}\right\}\notag\\
&\stackrel{(b)}=& \frac{\tilde{\Delta}_{T}^{2}V_{T}}{32T}\exp\left\{-\frac{8\tilde{C}\tilde{\Delta}_{T}V_{T}}{T}\cdot R_{j}^{\pi}\right\},\notag
\end{eqnarray}
where: the first three inequalities follow arguments given in step 3 in the proof of Theorem 3; $(a)$ holds by (\ref{eq:phierror}); and $(b)$ holds by $\delta = \sqrt{2V_{T}\tilde{\Delta}_{T}/T}$. Assume that $\sqrt{\tilde{C}}\cdot V_{T} \leq 2T$. Then, taking $\tilde{\Delta}_{T} =~ \left\lceil\left(\frac{4}{\tilde{C}}\right)^{1/3} \left(\frac{T}{V_{T}}\right)^{2/3}\right\rceil$, one has:\vspace{-0.3cm}
\begin{eqnarray}
R^{\pi}_{j}
&\geq& \frac{1}{32}\cdot \left(\frac{4}{\tilde{C}}\right)^{2/3}\left(\frac{T}{V_{T}}\right)^{1/3} \exp\left\{-\frac{8\tilde{C}V_{T}}{T}\cdot\left(\left(\frac{4}{\tilde{C}}\right)^{1/3}\left(\frac{T}{V_{T}}\right)^{2/3} +1 \right)\cdot R_{j}^{\pi}\right\}\notag\\
&\geq& \frac{1}{32}\cdot \left(\frac{4}{\tilde{C}}\right)^{2/3}\left(\frac{T}{V_{T}}\right)^{1/3} \exp\left\{-16\tilde{C}^{2/3}\cdot 4^{1/3}\cdot \left(\frac{V_{T}}{T}\right)^{1/3}R_{j}^{\pi}\right\},\notag
\end{eqnarray}
where the last inequality follows from $\sqrt{\tilde{C}}\cdot V_{T} \leq 2T$. Then, for $\beta = 16\left(4\tilde{C}^{2}\cdot \frac{V_{T}}{T}\right)^{1/3}$, one has:\vspace{-0.3cm}
\begin{equation}\label{eq:inequality}
\beta R_{j}^{\pi}
\;\geq\; \frac{32T}{\tilde{\Delta}_{T}^{2}V_{T}}
\;\geq\; \exp\left\{-\beta R_{j}^{\pi}\right\}.\vspace{-0.3cm}
\end{equation}
Let $y_{0}$ be the unique solution to the equation $y = \exp\left\{-y\right\}$. Then, (\ref{eq:inequality}) implies $\beta R_{j}^{\pi} \geq y_{0}$. In particular, since $y_{0}>1/2$ this implies $R_{j}^{\pi} \geq 1/\left(2\beta\right) =\frac{1}{32\left(2\tilde{C}\right)^{2/3}}\left(\frac{T}{V_{T}}\right)^{1/3}$ for all $1\leq j \leq m$. Hence:\vspace{-0.3cm}
\begin{eqnarray}
\mathcal{R}^{\pi}_{\phi^{(0)}}\left(\mathcal{V}_{s},T\right)
&\geq& \sum_{j=1}^{m}R_{j}^{\pi}
\;\geq\; \frac{T}{\tilde{\Delta}_{T}}\cdot\frac{1}{32\left(2\tilde{C}\right)^{2/3}}\left(\frac{T}{V_{T}}\right)^{1/3}\notag\\
&\stackrel{(a)}\geq& \frac{1}{64\cdot 2^{4/3}\tilde{C}^{1/3}}  \cdot V_{T}^{1/3}T^{2/3},\notag
\end{eqnarray}
where $(a)$ holds if $\sqrt{\tilde{C}}\cdot V_{T} \leq 2T$. If $\sqrt{\tilde{C}}\cdot V_{T} > 2T$, by Proposition~\ref{prop:lin} there is a constant $C$ such that $\mathcal{R}^{\pi}_{\phi^{(0)}}\left(\mathcal{V}_{s},T\right)\geq C T \geq C V_{T}^{1/3}T^{2/3}$; the last inequality holds by $T\geq V_{T}$. This concludes the proof.~\qed




\small
\vspace{-0.1cm}
\setstretch{1}
\bibliographystyle{chicago}
\bibliography{regret_bib2}
\normalsize

\newpage
\setstretch{1.3}
\noindent
\begin{center}
{\Large \textbf{Online Companion:\\  Non-stationary Stochastic Optimization}}
\end{center}

{\large{
\begin{center}
\hspace{0mm}{\sf Omar Besbes} \hspace{26mm}
{\sf Yonatan Gur} \hspace{27mm}
{\sf Assaf Zeevi}\symbolfootnote[1]{Correspondence: {\tt ob2105@columbia.edu}, {\tt ygur@stanford.edu}, {\tt assaf@gsb.columbia.edu}}
\\  Columbia University \hspace{13mm}Stanford University \hspace{13mm}Columbia University

\end{center}}}

\def\spacen{1.05}  
\def\spaces{1.0}  


\setcounter{footnote}{0}
\setcounter{page}{1}
\vspace*{-0.7cm}

\appendix
\setcounter{section}{1}

\renewcommand{\theequation}{\thesection-\arabic{equation}}
\setcounter{equation}{0}


\setstretch{1.4}
\section{Proofs of additional results}  \label{app:add}\vspace{-0.3cm}

\paragraph*{Proof of Proposition \ref{prop:lin}.} The proof of the proposition is established in two steps. In the first step, we limit nature to a class of function sequences $\mathcal{V}'$ where in every epoch nature is limited to one of two specific cost functions, and show that $\mathcal{V}'\subset \mathcal{V}$. In the second step, we show that whenever $\phi \in \left\{\phi^{(0)},\phi^{(1)}\right\}$, any admissible policy must incur regret of at least order $T$, even when nature is limited to the set $\mathcal{V}'$.

\textbf{Step 1.} Let $\mathcal{X} = [0,1]$ and fix $T\geq 1$. Let $V_{T}\in\left\{1,\ldots,T\right\}$ and assume that $C_{1}$ is a constant such that $V_{T}\geq C_{1}T$. Let $C =~\min\left\{C_{1}, \left(\frac{1}{2}-\nu\right)^{2}\right\}$ where $\nu$ appears in (\ref{eq:interior}), and we assume $\nu< 1/2$. Consider the following two quadratic functions:\vspace{-0.3cm}
\[
   f^{1}(x) = x^{2} - x + \frac{3}{4},
    \quad\quad\quad\quad
   f^{2}(x) = x^{2} - \left(1+ 2C \right)x + \frac{3}{4} + C.\vspace{-0.3cm}
\]
Denoting $x_{k}^{\ast} = \arg\min_{x\in[0,1]}f^{k}(x)$, we have $x_{1}^{\ast} = \frac{1}{2}$, and $x_{2}^{\ast} = \frac{1}{2}+C$. Define $\mathcal{V}' = \left\{f\;; f_{t}\in\left\{f^{1},f^{2}\right\}\;\forall t\in \mathcal{T}\right\}$. Then, for any sequence in $\mathcal{V}'$ the total functional variation is:\vspace{-0.3cm}
\[
\sum_{t=2}^{T}\sup_{x\in\mathcal{X}}\left|f_{t} - f_{t-1}\right|
\;\leq\; \sum_{t=2}^{T}\sup_{x\in\mathcal{X}}\left|2Cx - C\right|
 \;\leq\;  C T
 \;\leq\; C_{1}T
 \;\leq\; V_{T}.\vspace{-0.3cm}
\]
For any sequence in $\mathcal{V}'$ the total functional variation (\ref{eq:TVF}) is bounded by $V_{T}$, and therefore $\mathcal{V}'\subset \mathcal{V}$.

\textbf{Step 2.} Fix $\phi \in \left\{\phi^{(0)},\phi^{(1)}\right\}$, and let $\pi\in\mathcal{P}_{\phi}$. Let $\tilde{f}$ to be a random sequence in which in each epoch $f_{t}$ is drawn according to a discrete uniform distribution over $\left\{f^{1},f^{2}\right\}$ ($\tilde{f}_{t}$ is independent of $\mathcal{H}_{t}$ for any $t\in \mathcal{T}$). Any realization of $\tilde{f}$ is a sequence in $\mathcal{V}'$. In particular, taking expectation over $\tilde{f}$, one has:\vspace{-0.5cm}
\begin{eqnarray}
\mathcal{R}^{\pi}_{\phi}(\mathcal{V}',T) &\geq& \mathbb{E}^{\pi,\tilde{f}}\left[\sum_{t=1}^{T}\tilde{f}_{t}(X_{t}) - \sum_{t=1}^{T}\tilde{f}_{t}(x_{t}^{\ast})\right] \notag \\
 &=& \mathbb{E}^{\pi}\left[\sum_{t=1}^{T}\left(\frac{1}{2}\left(f^{1}(X_{t})+ f^{2}(X_{t})\right) - \frac{1}{2}\left(f^{1}(x_{1}^{\ast})+ f^{2}(x_{2}^{\ast})\right)\right)\right]  \notag \\ 
 &\geq& \sum_{t=1}^{T}\min_{x\in\left[0,1\right]}\left\{x^{2} - \left(1+C\right)x + \frac{1}{4} + \frac{C}{2} + \frac{C^{2}}{2}\right\}
 \;=\;  T\cdot \frac{C^{2}}{4},\notag
\end{eqnarray}
where the minimum is obtained at $x^{\ast} = \frac{1+C}{2}$. Since $\mathcal{V}'\subseteq \mathcal{V}$, we have established that\vspace{-0.4cm}
\[
\mathcal{R}^{\pi}_{\phi}(\mathcal{V},T)
\;\geq\; \mathcal{R}^{\pi}_{\phi}(\mathcal{V}',T)
\;\geq\; \frac{C^{2}}{4}\cdot T,\quad\vspace{-0.5cm}
\]
which concludes the proof.\qed\vspace{-0.2cm}

\paragraph*{Proof of Proposition \ref{thm:ebe}.} Fix $T\geq 1$, and $1\leq V_{T}\leq T$. Let $\pi$ be the OGD algorithm with $\eta_{t+1} = \eta$ 
for any $t=1,\ldots,T-1$. Fix $\Delta_{T}\in \left\{1,\ldots,T\right\}$ (to be specified below), and define a partition of $\mathcal{T}$ into batches $\mathcal{T}_{1},\ldots,\mathcal{T}_{m}$ of size $\Delta_{T}$ each (except perhaps $\mathcal{T}_{m}$) according to (\ref{eq:batches}); this partition is only for analysis purposes. Fix $f\in \mathcal{V}$. By \cite{FLA2005} we have that (see analysis in their Lemma~3.1):\footnote{The expression adjusts the analysis in \cite{FLA2005} to allow an arbitrary (and apotentially random)$X_{0}$.}\vspace{-0.2cm}
\begin{equation}\label{eq:fromlem3.1}
\mathbb{E}^{\pi}\left[\sum_{t\in \mathcal{T}_{j}}f_{t}(X_{t})\right] - \inf_{x\in\mathcal{X}}\left\{\sum_{t\in \mathcal{T}_{j}}f_{t}(x)\right\}
\;\leq\; \frac{4r^{2}}{\eta} + \Delta_{T}\cdot\frac{\eta G^{2}}{2},\vspace{-0.2cm}
\end{equation}
for any $j=1,\ldots,m$, where $r$ is the radius of the set $\mathcal{X}$. Following the proof of Proposition \ref{prop:meta}, we have:
\vspace{-0.2cm}
\begin{eqnarray}
\mathbb{E}^{\pi}\left[\sum_{t=1}^{T}f_{t}(X_{t})\right] - \sum_{t=1}^{T}f_{t}(x^{\ast}_{t})
 &\leq& \sum_{j=1}^{m}\left(\mathbb{E}^{\pi}\left[\sum_{t\in \mathcal{T}_{j}}f_{t}(X_{t})\right] - \inf_{x\in\mathcal{X}}\left\{\sum_{t\in \mathcal{T}_{j}}f_{t}(x)\right\}\right) + 2\cdot \Delta_{T}V_{T}\notag\\
&\stackrel{(a)} \leq& \frac{8T}{\Delta_{T}}\cdot \frac{r^{2}}{\eta} + T\eta G^{2}  + 2\cdot \Delta_{T}V_{T}\notag,
\end{eqnarray}
for any $1\leq\Delta_{T}\leq T$, where $(a)$ follows (\ref{eq:fromlem3.1}). Taking $\Delta_{T} = \left\lceil\left(T/V_{T}\right)^{2/3}\right\rceil$ and $\eta = \frac{r}{G}\left(V_{T}/T\right)^{1/3}$ we get:\vspace{-0.3cm}
\[
\mathbb{E}^{\pi}\left[\sum_{t=1}^{T}f_{t}(X_{t})\right] - \sum_{t=1}^{T}f_{t}(x^{\ast}_{t})
 \;\leq\; \left(9rG + 4\right) \cdot V^{1/3}T^{2/3}.\vspace{-0.2cm}
\]
Since the above holds for any $f\in \mathcal{V}$, we established:\vspace{-0.2cm}
\[
\mathcal{R}^{\pi}_{\phi^{(1)}}(\mathcal{V},T)
 \;\leq\; \left(9rG + 4\right) \cdot V^{1/3}T^{2/3}.\vspace{-0.2cm}
\]
This concludes the proof.\qed\vspace{-0.2cm}

\paragraph{Proofs of Lemma \ref{lem:KL} and Lemma \ref{lem:KL2}.} We start by proving Lemma \ref{lem:KL}. Suppose that 
$\phi = \phi^{(1)}$. In the proof we use the notation defined in \S4 and in the proof of Theorem~3. For any $t\in \mathcal{T}$ denote $Y_{t} = \phi^{(1)}(X_{t},\cdot)$, and denote by $y_{t}\in \mathbb{R}^{d}$ the realized feedback observation at epoch $t$. For convenience, for any $t\geq 1$ we further denote $y^{t} = \left(y_{1},\ldots,y_{t}\right)$. Fix $\pi\in\mathcal{P}_{\phi}$. Letting $u\in \mathcal{U}$, we denote $x_{1} = \pi_{1}(u)$, and $x_{t} := \pi_{t}\left(y^{t-1},u\right)$ for $t\in\left\{2,\ldots,T\right\}$. For any $f\in\mathcal{F}$ and $\tau\geq 2$, one has:\vspace{-0.2cm}
\begin{eqnarray}
d\mathbb{P}_{f}^{\pi,\tau}\left\{y^{\tau},u\right\}
&=&
d\mathbb{P}_{f}\left\{y^{\tau-1},u\right\}
d\mathbb{P}_{f}^{\pi,\tau-1}\left\{y^{\tau-1},u\right\}\notag\\
&\stackrel{(a)}=&
d\mathbb{P}_{f}\left\{y_{\tau}\;|x_{\tau}\right\}
d\mathbb{P}_{f}^{\pi,\tau-1}\left\{y^{\tau-1},u\right\}\notag\\
\label{eq:product}
&\stackrel{(b)}=&
dG\left(y_{\tau} - \nabla f(x_{\tau})\right)
d\mathbb{P}_{f}^{\pi,\tau-1}\left\{y^{\tau-1},u\right\},
\end{eqnarray}
where: $(a)$ holds since by the first part of Assumption~1 the feedback at epoch $\tau$ depends on the history only through $x_{\tau} =~ \pi_{\tau}\left(y^{\tau-1},u\right)$; and $(b)$ follows from the feedback structure given in the first part of Assumption~1. Fix $f,g \in \mathcal{F}$ and $\tau\geq 2$. One has:\vspace{-0.2cm}
\begin{eqnarray}
\mathcal{K}\left(\mathbb{P}^{\pi,\tau}_{f}\|\mathbb{P}^{\pi,\tau}_{g}\right)
 &=& \int_{u,y^{\tau}} \log\left(
\frac{d\mathbb{P}^{\pi,\tau}_{f}\left\{y^{\tau},u\right\}}
{d\mathbb{P}^{\pi,\tau}_{g}\left\{y^{\tau},u\right\}}\right)
d\mathbb{P}^{\pi,\tau}_{f}\left\{y^{\tau},u\right\}\notag\\
&\stackrel{(a)}=& \int_{u,y^{\tau}} \log\left(
\frac{dG\left(y_{\tau} - \nabla f(x_{\tau})\right)
d\mathbb{P}_{f}^{\pi,\tau-1}\left\{y^{\tau-1},u\right\}}
{dG\left(y_{\tau} - \nabla g(x_{\tau})\right)
d\mathbb{P}_{g}^{\pi,\tau-1}\left\{y^{\tau-1},u\right\}}\right)
dG\left(y_{\tau} - \nabla f(x_{\tau})\right)
d\mathbb{P}_{f}^{\pi,\tau-1}\left\{y^{\tau-1},u\right\}\notag
\end{eqnarray}
where $(a)$ holds by (\ref{eq:product}). We have that $\mathcal{K}\left(\mathbb{P}^{\pi,\tau}_{f}\|\mathbb{P}^{\pi,\tau}_{g}\right) = A_{\tau} + B_{\tau}$, where:\vspace{-0.4cm}
\begin{eqnarray}
A_{\tau}
&:=& \int_{u,y^{\tau}} \log\left(
\frac{d\mathbb{P}_{f}^{\pi,\tau-1}\left\{y^{\tau-1},u\right\}}
{d\mathbb{P}_{g}^{\pi,\tau-1}\left\{y^{\tau-1},u\right\}}\right)
dG\left(y_{\tau} - \nabla f(x_{\tau})\right)
d\mathbb{P}_{f}^{\pi,\tau-1}\left\{y^{\tau-1},u\right\}\notag\\
&=& \int_{u,y^{\tau-1}}\log\left(
\frac{d\mathbb{P}_{f}^{\pi,\tau-1}\left\{y^{\tau-1},u\right\}}
{d\mathbb{P}_{g}^{\pi,\tau-1}\left\{y^{\tau-1},u\right\}}\right)\left[\int_{y_{\tau}}
dG\left(y_{\tau} - \nabla f(x_{\tau})\right)\right]
d\mathbb{P}_{f}^{\pi,\tau-1}\left\{y^{\tau-1},u\right\}\notag\\
&=& \int_{u,y^{\tau-1}}\log\left(
\frac{d\mathbb{P}_{f}^{\pi,\tau-1}\left\{y^{\tau-1},u\right\}}
{d\mathbb{P}_{g}^{\pi,\tau-1}\left\{y^{\tau-1},u\right\}}\right)
d\mathbb{P}_{f}^{\pi,\tau-1}\left\{y^{\tau-1},u\right\}
\;=\; \mathcal{K}\left(\mathbb{P}^{\pi,\tau-1}_{f}\|\mathbb{P}^{\pi,\tau-1}_{g}\right),\notag
\end{eqnarray}
and\vspace{-0.6cm}
\begin{eqnarray}
B_{\tau}
&:=& \int_{u,y^{\tau}} \log\left(
\frac{dG\left(y_{\tau} - \nabla f(x_{\tau})\right)}
{dG\left(y_{\tau} - \nabla g(x_{\tau})\right)}\right)
dG\left(y_{\tau} - \nabla f(x_{\tau})\right)
d\mathbb{P}_{f}^{\pi,\tau-1}\left\{y^{\tau-1},u\right\}\notag\\
&=&\int_{u,y^{\tau-1}}\int_{y_{\tau}}\left[ \log\left(
\frac{dG\left(y_{\tau} - \nabla f(x_{\tau})\right)}
{dG\left(y_{\tau} - \nabla g(x_{\tau})\right)}\right)
dG\left(y_{\tau} - \nabla f(x_{\tau})\right)\right]
d\mathbb{P}_{f}^{\pi,\tau-1}\left\{y^{\tau-1},u\right\}\notag\\
&\stackrel{(b)}\leq& \tilde{C}\int_{u,y^{\tau-1}}\left\|\nabla f_{\tau}(x_{\tau}) - g_{\tau}(x_{\tau})\right\|^{2}
d\mathbb{P}_{f}^{\pi,\tau-1}\left\{y^{\tau-1},u\right\}
\;=\; \tilde{C}\mathbb{E}^{\pi}_{f}\left\|\nabla f_{\tau}(x_{\tau}) - g_{\tau}(x_{\tau})\right\|^{2},\notag
\end{eqnarray}
where $(b)$ follows the second part of Assumption~1. Repeating the above arguments, one has:\vspace{-0.4cm}
\[
\mathcal{K}\left(\mathbb{P}^{\pi,\tau}_{f}\|\mathbb{P}^{\pi,\tau}_{g}\right)
\;\leq\;\mathcal{K}\left(\mathbb{P}^{\pi,1}_{f}\|\mathbb{P}^{\pi,1}_{g}\right)
+ \tilde{C}\mathbb{E}^{\pi}_{f}\left[\sum_{t=2}^{\tau}\left\|\nabla f_{t}(x_{t}) - g_{t}(x_{t})\right\|^{2} \right].\vspace{-0.4cm}
\]
From the above it is also clear that:\vspace{-0.4cm}
\begin{eqnarray}
\mathcal{K}\left(\mathbb{P}^{\pi,1}_{f}\|\mathbb{P}^{\pi,1}_{g}\right)
&=&\int_{u,y_{1}} \log\left(
\frac{d\mathbb{P}^{\pi,1}_{f}\left\{y_{1},u\right\}}
{d\mathbb{P}^{\pi,1}_{g}\left\{y_{1},u\right\}}\right)
d\mathbb{P}^{\pi,1}_{f}\left\{y_{1},u\right\}\notag\\
&=& \int_{u}\left[\int_{y_{1}} \log\left(
\frac{dG\left(y_{1} - \nabla f(x_{1})\right)}
{dG\left(y_{1} - \nabla g(x_{1})\right)}\right)
dG\left(y_{1} - \nabla f(x_{1})\right)
\right]d\mathbf{P}_{u}\left\{u\right\}\notag\\
&\leq&\tilde{C}\int_{u}\left\|\nabla f_{1}(x_{1}) - \nabla g_{1}(x_{1})  \right\|^{2}d\mathbf{P}_{u}\left\{u\right\}
\;=\; \tilde{C}\mathbb{E}^{\pi}_{f}\left\|\nabla f_{1}(x_{1}) - \nabla g_{1}(x_{1})  \right\|^{2}.\notag
\end{eqnarray}
Hence, we have established that for any $\tau\geq1$:\vspace{-0.4cm}
\[
\mathcal{K}\left(\mathbb{P}^{\pi,\tau}_{f}\|\mathbb{P}^{\pi,\tau}_{g}\right)
\;\leq\;  \tilde{C}\sum_{t=1}^{\tau}\mathbb{E}^{\pi}_{f}\left\|\nabla f_{t}(x_{t}) - g_{t}(x_{t})\right\|^{2}.\vspace{-0.4cm}
\]
Finally, following the steps above, the proof of Lemma \ref{lem:KL2} (for the feedback structure $\phi =~\phi^{(0)}$) is immediate, using the notation introduced in the proof of Theorem~5 for cost feedback structure, along with Assumption~2. This concludes the proof.\qed

\vspace{-0.1cm}
\paragraph*{Performance analysis of OGD algorithm without restarting.} We onsider the performance of the OGD algorithm \emph{without restarting}, relative to the dynamic benchmark. The following illustrates that this algorithm will yield linear regret for a broad set of variation budgets.
\begin{example}{\textbf{\textup{(Failure of OGD without restarting)}}} \label{ex:LinearRegretUnderOGD}
{\rm
Consider a partition of the horizon $\mathcal{T}$ into batches $\mathcal{T}_{1},\ldots \mathcal{T}_{m}$ according to (\ref{eq:batches}), with each batch of size $\Delta_{T}$. Consider the following cost functions:\vspace{-0.4cm}
\[
g_{1}(x) = (x-\alpha)^{2},\quad\quad\quad g_{2}(x) = x^{2};\quad \quad\quad x\in\left[-1,3\right].\vspace{-0.4cm}
\]
Assume that nature selects the cost function to be $g_{1}(\cdot)$ in the even batches and $g_{2}(\cdot)$ in the odd batches. Assume that at every epoch $t$, after selecting an action $x_{t}\in \mathcal{X}$, a \emph{noiseless} access to the gradient of the cost function at point $x_{t}$ is granted, that is, $\phi^{(1)}_{t}(x,f_{t}) = f'_{t}(x)$ for all $x\in\mathcal{X}$ and $t\in \mathcal{T}$. Assume that the decision maker is applying the OGD algorithm with a sequence of step sizes $\left\{\eta_{t}\right\}_{t=2}^{T}$, and $x_{1}=1$. We consider two classes of step size sequences that have been shown to be rate optimal in two instances of OCO settings (see \cite{FLA2005}, and \cite{Haz-Aga-Kal2007}).
\begin{enumerate}\vspace{-0.2cm}
  \item Suppose $\eta_{t} = \eta = C/\sqrt{T}$. Then, selecting a batch size $\Delta_{T}$ of order $\sqrt{T}$, and $\alpha = 1+\left(1+2\eta\right)^{\Delta_{T}}$, the variation budget $V_{T}$ is at most of order $\sqrt{T}$, and there is a constant $C_{1}$ such that $\mathcal{R}^{\pi}_{\phi}(\mathcal{V},T)\geq~C_{1}T$.\vspace{-0.2cm}
  \item Suppose that $\eta_{t} = C/t$. Then, selecting a batch size $\Delta_{T}$ of order $T$, and $\alpha = 1$, the variation budget $V_{T}$ is a fixed constant, and there is a constant $C_{2}$ such that $\mathcal{R}^{\pi}_{\phi}(\mathcal{V},T)\geq C_{2}T$.
      \quad\qed
\end{enumerate}}
\end{example}

\paragraph*{Proof of claims made in Example \ref{ex:LinearRegretUnderOGD}.} Fix $T\geq 1$. Let $\mathcal{X} = [-1,3]$ (we assume that $\nu$, appearing in (\ref{eq:interior}), is smaller than 1) and consider the following two functions: $g^{1}(x) = (x-\alpha)^{2}$, and $g^{2}(x) = x^{2}$. We assume that in each epoch $t$, after selecting an action $x_{t}$, there is a noiseless access to the gradient of the cost function, evaluated at point $x_{t}$. The the deterministic actions are generated by an OGD algorithm:\vspace{-0.2cm}
\[
x_{t+1} = P_{\mathcal{X}}\left(x_{t} - \eta_{t+1}\cdot f'_{t}(x_{t})\right),\quad\quad \text{ for all } t\geq1,\vspace{-0.2cm}
\]
with the initial selection $x_{1}=1$. In the first part we consider the case of $\eta_{t}= \eta = C/\sqrt{T}$, and in a second part we consider the case of $\eta_{t} = C/t$. The structure of both parts is similar: first we analyze the variation of the instance, showing it is sublinear. Then, by analyzing the sequence of decisions $\left\{x_{t}\right\}_{t=1}^{T}$ that is generated by the Online Gradient Descent policy, we show that in a linear portion of the horizon there is a constant $C_{2}$ such that $\left|x_{t}- x_{t}^{\ast}\right|> C_{2}$, and therefore a linear regret is incurred.

\textbf{Part 1.} Assume that $\eta_{t}= \eta = C/\sqrt{T}\leq 1/2$. Select $\Delta_{T}=\left\lfloor 1 + \frac{1}{2\eta} \right\rfloor$, and set $\alpha = 1 + \left(1-2\eta\right)^{\Delta_{T}}$ (note that $1\leq \alpha \leq 2$). We assume that nature selects the cost function to be $g_{1}(\cdot)$ in the even batches and $g_{2}(\cdot)$ in the odd batches. We start by analyzing the variation along the horizon:\vspace{-0.3cm}
\begin{eqnarray}
\sum_{t=2}^{T}\sup_{x\in \mathcal{X}}\left|f_{t}(x) - f_{t-1}(x)\right|
&\leq& \left(\left\lceil\frac{T}{\Delta_{T}}\right\rceil-1\right)\cdot \sup_{x\in \mathcal{X}}\left|g_{2}(x) - g_{1}(x)\right|\notag\\
&\leq& \frac{T}{\Delta_{T}}\cdot \sup_{x\in \mathcal{X}}\left|\alpha^{2} - 2\alpha x\right|\notag\\
&\stackrel{(a)} \leq& \frac{8T}{\Delta_{T}}
\;=\; \frac{8T}{\left\lfloor 2 + \frac{1}{2\eta} \right\rfloor}\notag\\
&\leq& 16T\eta
  \;=\; 16C\cdot \sqrt{T},\notag
\end{eqnarray}
where $(a)$ follows from $1\leq \alpha \leq 2$ and $-1\leq x \leq 3$. Next, we analyze the incurred regret. We start by analyzing decisions generated by the OGD algorithm throughout the first two batches. Recalling that $x_{1} = 1$ and that $g_{2}(\cdot)$ is the cost function throughout the first batch, one has for any $2 \leq t \leq \Delta_{T}+1$:\vspace{-0.3cm}
\begin{eqnarray}
x_{t} &=& x_{t-1} - \eta\cdot f'(x_{t-1})\notag\\
 &=& x_{t-1} - \eta\cdot 2x_{t-1}
 \;=\; x_{t-1}\left(1-2\eta\right)\notag\\
 &=& x_{1}\left(1-2\eta\right)^{t-1}
 \;=\; \left(1-2\eta\right)^{t-1}\notag\\
&=& \exp\left\{(t-1)\ln\left(1-2\eta  \right)   \right\}\notag\\
&\stackrel{(a)} \geq& \exp\left\{(t-1)(-2\eta - 2\eta^{2})   \right\}\notag\\
&\stackrel{(b)} \geq& \exp\left\{-1-\eta   \right\}\notag\\
&\stackrel{(c)} >& \frac{1}{e^{2}},\notag
\end{eqnarray}
where: $(a)$ follows since for any $-1< x \leq 1$ one has $\ln(1+x)\geq x - \frac{x^{2}}{2}$; $(b)$ follows from $t \leq \Delta_{T}\leq 1 + \frac{1}{2\eta}$; and $(c)$ follows from $\eta \leq \frac{1}{2} < 1$. Since $x_{t}^{\ast} = 0$ for any $1\leq t \leq \Delta_{T}$, one has:\vspace{-0.3cm}
\[
x_{t} - x_{t}^{\ast} \;>\; \frac{1}{e^{2}},\vspace{-0.3cm}
\]
for any $1\leq t \leq \Delta_{T}$. At the end of the first batch the cost function changes from $f(\cdot)$ to $g(\cdot)$. Note that the first action of the second batch is $x_{\Delta_{T}+1} = \left(1-2\eta\right)^{\Delta_{T}}$. Since $g_{1}(\cdot)$ is the cost function throughout the second batch, for any $\Delta_{T}+2\leq t \leq 2\Delta_{T}+1$ one has:\vspace{-0.4cm}
\begin{eqnarray}
x_{t}&=& x_{t-1} - \eta\cdot g'(x_{t-1})\notag\\
&=& x_{t-1} - \eta\cdot 2\left(x_{t-1} - \alpha\right).\notag
\end{eqnarray}
Using the transformation $y_{t} = x_{t} - \alpha$ for all $t$, one has:\vspace{-0.3cm}
\begin{eqnarray}
y_{t} &=& y_{t-1} - \eta\cdot 2y_{t-1}
 \;=\; y_{t-1}\left(1-2\eta\right)\notag\\
 &=& y_{\Delta_{T}+1}\left(1-2\eta\right)^{t-\Delta_{T}-1}\notag\\
 &=& x_{\Delta_{T}+1}\left(1-2\eta\right)^{t-\Delta_{T}-1}
 - \alpha \left(1-2\eta\right)^{t-\Delta_{T}-1}\notag\\
 &=& \left(1-2\eta\right)^{t-1}
 - \left(1-2\eta\right)^{t-\Delta_{T}-1}
 - \left(1-2\eta\right)^{t-1}\notag\\
 &=& - \left(1-2\eta\right)^{t-\Delta_{T}-1}\notag\\
 &=&  - \exp\left\{\left(t-\Delta_{T}-1\right)\ln\left(  1-2\eta\right)    \right\}\notag\\
&\stackrel{(a)} \leq& -\exp\left\{(t-\Delta_{T}-1)(-2\eta - 2\eta^{2})   \right\}\notag\\
&\stackrel{(b)} \leq& -\exp\left\{-1-\eta   \right\}\notag\\
&\stackrel{(c)} <& -\frac{1}{e^{2}},\notag
\end{eqnarray}
where: $(a)$ holds since for any $-1< x \leq 1$ one has $\ln(1+x)\geq x - \frac{x^{2}}{2}$; $(b)$ follows from $t \leq 2\Delta_{T}\leq 1 + \frac{1}{2\eta}+\Delta_{T}$; and $(c)$ follows from $\eta \leq \frac{1}{2} < 1$. Finally, recalling that $x_{t}^{\ast} = \alpha$ and using the transformation $y_{t} = x_{t} - \alpha$, one has for any $\Delta_{T}+1 \leq t \leq 2\Delta_{T}$:\vspace{-0.2cm}
\[
x_{t}^{\ast} - x_{t}
 \;=\; y_{t}
 \;<\; -\frac{1}{e^{2}}.\vspace{-0.2cm}
\]
In the beginning of the third batch $g_{2}(\cdot)$ becomes the cost function once again. We note that the first action of the third batch is the same as the first action of the first batch:\vspace{-0.3cm}
\[
x_{2\Delta_{T} + 1}
 \;=\; \alpha + y_{2\Delta_{T} + 1}
 \;=\;  \alpha - \left(1-2\eta\right)^{2\Delta_{T} + 1-\Delta_{T}-1}
 \;=\;  \alpha - \left(1-2\eta\right)^{\Delta_{T}}
 \;=\; 1
 \;=\; x_{1},\vspace{-0.3cm}
\]
and therefore the actions taken in the first two batches are repeated throughout the horizon. We conclude that for any $1\leq t \leq T$,\vspace{-0.2cm}
\[
\left|x_{t} - x_{t}^{\ast}\right|\;>\; \frac{1}{e^{2}}.\vspace{-0.2cm}
\]
Finally, we calculate the regret incurred throughout the horizon. Using Taylor expansion, one has\vspace{-0.3cm}
\[
\sum_{t=1}^{T}\left(f_{t}(x_{t}) - f_{t}(x_{t}^{\ast})\right)
\;=\; \sum_{t=1}^{T}\left(x_{t}- x_{t}^{\ast}\right)^{2}
\;>\; \sum_{t=1}^{T}\frac{1}{e^{4}}
\;=\; \frac{T}{e^{4}}.\vspace{-0.3cm}
\]

\textbf{Part 2.} For concreteness we assume in this part that $T$ is even and larger than $2$. We show that linear regret can be incurred when $\eta_{t} = \frac{C}{t}$. Set $\alpha = 1$ and $\Delta_{T} =  T/2$ (therefore we have two batches). Assume that nature selects $g_{1}(\cdot)$ to be the cost function in the first batch, $g_{2}(\cdot)$ to be the cost function in the second batch. We start by analyzing the variation along the horizon. Recalling that there is only one change in the cost function, one has:\vspace{-0.4cm}
\begin{eqnarray}
\sum_{t=2}^{T}\sup_{x\in \mathcal{X}}\left|f_{t}(x) - f_{t-1}(x)\right|
&=& \sup_{x\in \mathcal{X}}\left|g_{2}(x) - g_{1}(x)\right|\notag\\
&=& \sup_{x\in \mathcal{X}}\left|\alpha^{2} - 2\alpha x\right|
 \;=\; \sup_{x\in \mathcal{X}}\left|1 - 2x\right|
\;\stackrel{(a)} =\; 5,\notag
\end{eqnarray}
where $(a)$ holds because $-1\leq x \leq 3$. Since $x_{1} = 1$, and $g_{1}'(1)=0$, one obtains $x_{t}=1$ for all $1\leq t \leq \left\lceil\frac{T}{2}\right\rceil + 1$. After $\left\lceil T/2 \right\rceil$ epochs, the cost function changes from $g_{1}(\cdot)$ to $g_{2}(\cdot)$, and for all $\left\lceil\frac{T}{2}\right\rceil +2 \leq t \leq T$ one has:\vspace{-0.3cm}
\begin{eqnarray}
x_{t} &=& x_{t-1} - \eta_{t}\cdot g_{2}'(x_{t-1})\notag\\
 &=& x_{t-1} - \eta_{t}\cdot 2x_{t-1}
 \;=\; x_{t-1}\left(1-2\eta_{t}\right)\notag\\
 &=& x_{\frac{T}{2}+1}\prod_{t' = \frac{T}{2}+1}^{t}\left(1-2\eta_{t'}\right)
\;=\; \prod_{t' = \frac{T}{2}+1}^{t}\left(1-2\eta_{t'}\right)\notag\\
&\stackrel{(a)} \geq& \left(1-2\eta_{\frac{T}{2} + 2}\right)^{t - \frac{T}{2} - 1}\
\;=\;  \left(1-\frac{4C}{T+4}\right)^{t - \frac{T}{2} - 1}\notag\\
&=& \exp\left\{\left(t - \frac{T}{2} - 1\right)\ln\left(1-\frac{4C}{T+4} \right)   \right\}\notag\\
&\stackrel{(b)} \geq& \exp\left\{\left(t - \frac{T}{2} - 1\right)\left(-\frac{4C}{T+4} - \frac{8C}{\left(T+4\right)^{2}}\right)   \right\}\notag\\
&\stackrel{(c)} \geq& \exp\left\{-4C -\frac{8C^{2}}{T+4}  \right\}
\;>\; \exp\left\{-4C -2C^{2}  \right\},\notag
\end{eqnarray}
where: $(a)$ holds since $\left\{\eta_{t}\right\}$ is a decreasing sequence; $(b)$ holds since $\ln(1+x)\geq x - \frac{x^{2}}{2}$ for any $-1< x \leq 1$; and $(c)$ is obtained using $t < T + \frac{T}{2} + 5$. Since $x_{t}^{\ast} = 0$ for any $\frac{T}{2} + 1 \leq t \leq T$, one has:\vspace{-0.3cm}
\[
x_{t} - x_{t}^{\ast} \;>\; \frac{1}{e^{2C\left(2+C\right)}},\vspace{-0.3cm}
\]
for all $\frac{T}{2} + 1 \leq t \leq T$. Finally, we calculate the regret incurred throughout the horizon. Recalling that throughout the first batch no regret is incurred, and using Taylor expansion, one has:\vspace{-0.3cm}
\[
\sum_{t=1}^{T}\left(f_{t}(x_{t}) - f_{t}(x_{t}^{\ast})\right)
\;=\; \sum_{t=\frac{T}{2}+1}^{T}\left(f(x_{t}) - f(x_{t}^{\ast})\right)
\;=\; \sum_{t=\frac{T}{2}+1}^{T}\left(x_{t}- x_{t}^{\ast}\right)^{2}
\;\geq\; \sum_{t=\frac{T}{2}+1}^{T}\frac{1}{e^{4C\left(2+C\right)}}
\;=\; \frac{T}{2e^{4C\left(2+C\right)}}.\vspace{-0.3cm}
\]
This concludes the proof.\qed\vspace{-0.2cm}

\section{Auxiliary results for OCO settings}  \label{app:sba}\vspace{-0.2cm}
\subsection{Preliminaries} \vspace{-0.2cm}
In this section we develop auxiliary results that provide bounds on the regret with respect to the single best action in the adversarial setting. As discussed in \S1, the OCO literature most often considers few different feedback structures; typical examples include full access to the cost/gradient after the action $X_{t}$ is selected, as well as a noiseless access to the cost/gradient evaluated at $X_{t}$. However, in this section we consider the feedback structures $\phi^{(0)}$ and $\phi^{(1)}$, where noisy access to the cost/gradient is granted.

We define admissible online algorithms exactly as admissible policies are defined in \S2.\footnote{We use the different terminology and notation only to highlight the different objectives: a policy $\pi$ is designed to minimize regret with respect to the dynamic oracle, while an online algorithm $\mathcal{A}$ is designed to minimize regret compared to the static single best action benchmark.} More precisely, letting $U$ be a random variable defined over a probability space $\left(\mathbb{U}, \mathcal{U},\mathbf{P}_{u}\right)$, we let $\mathcal{A}_{1}:\mathbb{U}\rightarrow \mathbb{R}^{d}$ and $\mathcal{A}_{t}:\mathbb{R}^{(t-1)k}\times\mathbb{U}\rightarrow \mathbb{R}^{d}$ for $t=2,3,\ldots$ be measurable functions, such that $X_{t}$, the action at time $t$, is given by\vspace{-0.2cm}
\begin{displaymath}
   X_{t} = \left\{
     \begin{array}{lr}
       \mathcal{A}_{1}\left(U\right) & t=1,\quad\quad\quad\quad\\
       \mathcal{A}_{t}\left(\phi_{t-1}\left(X_{t-1},f_{t-1}\right), \ldots, \phi_{1}\left(X_{1},f_{1}\right),U\right) & t=2,3,\ldots,\;\;
     \end{array}
   \right.\vspace{-0.1cm}
\end{displaymath}
where $k=1$ if $\phi = \phi^{(0)}$, and $k=d$ if $\phi = \phi^{(1)}$. The mappings $\left\{\mathcal{A}_{t}:\;t=1,\ldots,T\right\}$ together with the distribution $\mathbf{P}_{u}$ define the class of admissible online algorithms with respect to feedback $\phi$, which is exactly the class $\mathcal{P}_{\phi}$. The filtration $\left\{\mathcal{H}_{t},\;t=1,\ldots,T\right\}$ is defined exactly as in \S2. Given a feedback structure $\phi\in \left\{\phi^{(0)},\phi^{(1)}\right\}$, the objective is to minimize the regret compared to the single best action:\vspace{-0.1cm}
\[
\mathcal{G}^{\mathcal{A}}_{\phi}(\mathcal{F},T) = \sup_{f\in \mathcal{F}}\left\{\mathbb{E}^{\mathcal{A}}\left[\sum_{t=1}^{T}f_{t}(X_{t})\right]
 - \min_{x\in \mathcal{X}}\left\{\sum_{t=1}^{T}f_{t}(x)\right\}\right\}.\vspace{-0.1cm}
\]
We note that while most results in the OCO literature allow sequences that can adjust the cost function adversarially at each epoch, we consider the above setting where nature commits to a sequence of functions in advance. This, along with the setting of noisy cost/gradient observations, is done for the sake of consistency with the non-stationary stochastic framework we propose in this paper.\vspace{-0.2cm}

\subsection{Upper bounds} \vspace{-0.2cm}
The first two results of this section, Lemma \ref{lem:sbaRegret} and Lemma \ref{lem:sbaGrad}, analyze the performance of the EGS algorithm (given in \S5) under structure $(\mathcal{F}_{s},\phi^{(0)})$ and the OGD algorithm (given in \S4) under structure $(\mathcal{F}_{s},\phi^{(1)})$, respectively. To the best of our knowledge, the upper bound in Lemma \ref{lem:sbaRegret} is not documented in the Online Convex Optimization literature\footnote{The feasibility of an upper bound of order $\sqrt{T}$ on the regret in an adversarial setting with noisy access to the cost and with strictly convex cost functions was suggested by \cite{Aga-Dek-Xia2010} without further details or proof.}. Lemma \ref{lem:sbaGrad} adapts Theorem 1 in \cite{Haz-Aga-Kal2007} (that considered noiseless access to the gradient) to the feedback structure $\phi^{(1)}$.\vspace{-0.2cm}

\begin{lemma}\label{lem:sbaRegret}\textbf{\textup{(Performance of EGS in the adversarial setting)}}
Consider the feedback structure $\phi =~\phi^{(0)}$. Let $\mathcal{A}$ be the EGS algorithm given in \S5.2, with $a_{t} = 2d/Ht$ and $\delta_{t} = h_{t} = a_{t}^{1/4}$ for all $t\in~\left\{1,\ldots,T-1\right\}$. Then, there exists a constant $\bar{C}$, independent of $T$ such that for any $T\geq 1$,\vspace{-0.3cm}
\[
\mathcal{G}^{\mathcal{A}}_{\phi}\left(\mathcal{F}_{s},T\right)\;\leq\; \bar{C}\sqrt{T}.\vspace{-0.3cm}
\]
\end{lemma}

\textbf{Proof.} Let $\phi= \phi^{(0)}$. Fix $T\geq 1$ and $f\in\mathcal{F}_{s}$. Let $\mathcal{A}$ be the EGS algorithm, with the selection $a_{t} = 2d/Ht$ and $\delta_{t} = h_{t} = a_{t}^{1/4}$ for all $t\in \left\{1,\ldots,T-1\right\}$. We assume that $\delta_{t} \leq \nu$ for all $t\in \mathcal{T}$; in the end of the proof we discuss the case in which the former does not hold. For the sequence $\left\{\delta_{t}\right\}_{t=1}^{T}$, we denote by $\mathcal{X}_{\delta_{t}}$ the $\delta_{t}$-interior of the action set $\mathcal{X}$: $\mathcal{X}_{\delta_{t}} \;=\; \left\{x\in \mathcal{X}\;|\; \mathbf{B}_{\delta_{t}}(x) \subseteq \mathcal{X}\right\}$. We have for all $f_{t}\in\mathcal{F}_{s}$:\vspace{-0.3cm}
\begin{equation}\label{eq:sba0}
\mathbb{E}\left[\phi^{(0)}_{t}\left(X_{t},f_{t}\right)\;|X_{t} = x\right] \;=\; f_{t}(x) \quad\text{ and }\quad
\sup_{x\in \mathcal{X}}\left\{\mathbb{E}\left[\left(\phi^{(0)}_{t}\left(x,f_{t}\right)\right)^{2}\right]\right\}
\;\leq\; G^{2} + \sigma^{2},\vspace{-0.3cm}
\end{equation}
for some $\sigma \geq 0$. At any $t\in\mathcal{T}$ the gradient estimator is:\vspace{-0.3cm}
\[
\hat{\nabla}_{h_{t}}f_{t}(X_{t})
\;=\;\frac{\phi^{(0)}_{t}\left(X_{t} + h_{t}\psi_{t},f_{t}\right)\psi_{t}}{h_{t_{}}},\vspace{-0.3cm}
\]
for a fixed $h_{t}>0$, and where $\left\{\psi_{t}\right\}$ is a sequence of iid random variables, drawn uniformly over the set $\left\{\pm e^{(1)},\ldots, \pm e^{(d)}\right\}$, where $e^{(k)}$ denotes the unit vector with $1$ at the $k^{\text{th}}$ coordinate. In particular, we denote $\psi_{t} = Y_{t}W_{t}$, where $Y_{t}$ and $W_{t}$ are independent random variables, $\mathbb{P}\left\{y_{t} = 1\right\} = \mathbb{P}\left\{y_{t} = -1\right\} = 1/2$, and $W_{t} = e^{(k)}$ with probability $1/d$ for all $k\in\left\{1,\ldots, d\right\}$. The estimated gradient step is\vspace{-0.3cm}
\[
Z_{t+1} \;=\; P_{\mathcal{X}_{\delta_{t}}}\left(Z_{t} - a_{t}\hat{\nabla}_{h_{t}}f_{t}(Z_{t})\right),
\quad\quad\quad
X_{t+1} \;=\; Z_{t+1} + h_{t+1}\psi_{t},\vspace{-0.3cm}
\]
where $P_{\mathcal{X}_{\delta_{t}}}$ denotes the Euclidean projection operator over the set $\mathcal{X}_{\delta_{t}}$. Note that $Z_{t}\in\mathcal{X}$, $X_{t}\in\mathcal{X}$, and $X_{t}+ h_{t}\psi_{t}\in\mathcal{X}$ for all $t\in \mathcal{T}$. Since $\left\|\psi_{t}\right\| = 1$ for all $t\in \mathcal{T}$, one has:\vspace{-0.3cm}
\begin{equation}\label{eq:sba1}
\mathbb{E}\left[\left\|\hat{\nabla}_{h_{t}}f_{t}(Z_{t})\right\|^{2}\;|Z_{t} = z\right]
\;=\; \frac{\mathbb{E}\left[\left(\phi^{(0)}_{t}\left(z + h_{t}\psi_{t},f_{t}\right)\right)^{2}\right]}{h_{t}^{2}}
\;\leq\; \frac{G^{2}+ \sigma^{2}}{h_{t}^{2}}\quad\quad \text{ for all } z\in \mathcal{X},\vspace{-0.3cm}
\end{equation}
using (\ref{eq:sba0}). Then,\vspace{-0.3cm}
\[
\mathbb{E}\left[\hat{\nabla}_{h_{t}}f_{t}(Z_{t})\;|Z_{t} = z,\psi_{t}=\psi\right]
 \;=\; \frac{\mathbb{E}\left[\phi^{(0)}_{t}\left(Z_{t} + h_{t}\psi_{t},f_{t}\right)\psi_{t}\;|Z_{t} = z,\psi_{t} = \psi\right]}{h_{t}}
 \;=\; \frac{f_{t}\left(z + h_{t}\psi\right)\psi}{h_{t}}.\vspace{-0.3cm}
\]
Therefore, taking expectation with respect to $\psi$, one has\vspace{-0.2cm}
\begin{eqnarray}
\mathbb{E}\left[\hat{\nabla}_{h_{t}}f_{t}(Z_{t})\;|Z_{t} = z\right]
&=& \mathbb{E}_{Y,W}\left[ \frac{f_{t}\left(z + h_{t}\psi\right)\psi}{h_{t}} \right]
\;=\; \frac{1}{d}\sum_{k=1}^{d}\frac{\left(f_{t}(z+h_{t}e^{(k)})- f_{t}(z-h_{t}e^{(k)})\right)e^{(k)} }{2h_{t}}\notag\\
&\stackrel{(a)}\geq& \frac{1}{d}\sum_{k=1}^{d}\left(\nabla f_{t}(z- h_{t}e^{(k)})\cdot e^{(k)}\right)e^{(k)}\notag\\
&\stackrel{(b)}\geq& \frac{1}{d}\sum_{k=1}^{d}\left(\nabla f_{t}(z)\cdot e^{(k)} - Gh_{t}\right)e^{(k)}
\;=\; \frac{1}{d} \nabla f_{t}(z) -\frac{Gh_{t}}{d}\cdot \bar{e},\notag
\end{eqnarray}
where $\bar{e}$ denotes a vector of ones. The equalities and inequalities above hold componentwise, where $(a)$ follows from a Taylor expansion and the convexity of $f_{t}$: $f_{t}(z+h_{t}e^{(k)})- f_{t}(z-h_{t}e^{(k)})
\geq \nabla f_{t}(z-h_{t}e^{(k)})\cdot\left(2h_{t}e^{(k)}\right)$, for any $1\leq k \leq d$, and $(b)$ follows from a Taylor expansion, the convexity of $f_{t}$, and (\ref{eq:H}):\vspace{-0.3cm}
\[
\nabla f_{t}(z- h_{t}e^{(k)})\cdot e^{(k)}
\;\geq\; \nabla f_{t}(z)\cdot e^{(k)} - \left(he^{(k)}\right)\cdot\left(\nabla^{2}f_{t}\right)e^{(k)}
\;\geq\; \nabla f_{t}(z)\cdot e^{(k)} - Gh_{t},\vspace{-0.3cm}
\]
for any $1\leq k \leq d$. Therefore, for all $z\in\mathcal{X}$ and for all $t\in\mathcal{T}$:\vspace{-0.3cm}
\begin{equation}\label{eq:sba2}
\left\|\frac{1}{d} \nabla f_{t}(z) - \mathbb{E}\left[\hat{\nabla}_{h_{t}}f_{t}(Z_{t})\;|Z_{t} = z\right] \right\|
\;\leq\; \frac{Gh_{t}}{\sqrt{d}}.\vspace{-0.3cm}
\end{equation}
Define $x^{\ast}$ as the single best action: $x^{\ast} = \argmin_{x\in \mathcal{X}}\left\{\sum_{t=1}^{T}f_{t}(x)\right\}$. Then, for any $t\in \mathcal{T}$, one has\vspace{-0.3cm}
\[
f_{t}(x^{\ast})
\;\geq\; f_{t}(Z_{t}) + \nabla f_{t}(Z_{t})\cdot \left(x^{\ast} - Z_{t}\right) + \frac{1}{2}H\left\|x^{\ast} - Z_{t}\right\|^{2},\vspace{-0.3cm}
\]
and hence:\vspace{-0.3cm}
\begin{equation}\label{eq:sba3}
f_{t}(Z_{t})  - f_{t}(x^{\ast})
\;\leq\; \nabla f_{t}(Z_{t})\cdot \left(Z_{t} - x^{\ast}\right) - \frac{1}{2}H\left\|Z_{t} - x^{\ast}\right\|^{2}.\vspace{-0.3cm}
\end{equation}
Next, using the estimated gradient step, one has\vspace{-0.5cm}
\begin{eqnarray}
\left\|Z_{t+1} - x^{\ast}\right\|^{2}
&=& \left\|P_{\mathcal{X}_{\delta_{t}}}\left(Z_{t} - a_{t}\hat{\nabla}_{h_{t}}f_{t}(Z_{t})\right) - x^{\ast}\right\|^{2}\notag\\
&\stackrel{(a)}\leq& \left\|Z_{t} - a_{t}\hat{\nabla}_{h_{t}}f_{t}(Z_{t}) - x^{\ast}\right\|^{2}\notag\\
&=& \left\|Z_{t} - x^{\ast}\right\|^{2}
 - 2a_{t}\left(Z_{t} - x^{\ast}\right)\cdot\hat{\nabla}_{h_{t}}f_{t}(Z_{t})
 + a_{t}^{2}\left\|\hat{\nabla}_{h_{t}}f_{t}(Z_{t})\right\|^{2}\notag\\
&=& \left\|Z_{t} - x^{\ast}\right\|^{2}
 - \frac{2a_{t}}{d}\cdot\left(Z_{t} - x^{\ast}\right)\cdot\nabla f_{t}(Z_{t})
 + a_{t}^{2}\left\|\hat{\nabla}_{h_{t}}f_{t}(Z_{t})\right\|^{2}\notag\\
& & +\; 2a_{t}\left(Z_{t} - x^{\ast}\right)\cdot\left( \frac{1}{d}\nabla f_{t}(Z_{t}) - \hat{\nabla}_{h_{t}}f_{t}(Z_{t}) \right)\notag\\
&\leq& \left\|Z_{t} - x^{\ast}\right\|^{2}
 - \frac{2a_{t}}{d}\cdot\left(Z_{t} - x^{\ast}\right)\cdot\nabla f_{t}(Z_{t})
 + a_{t}^{2}\left\|\hat{\nabla}_{h_{t}}f_{t}(Z_{t})\right\|^{2}\notag\\
& & +\; 2a_{t}\left\|Z_{t} - x^{\ast}\right\|\cdot\left\| \frac{1}{d}\nabla f_{t}(Z_{t}) - \hat{\nabla}_{h_{t}}f_{t}(Z_{t}) \right\|,\notag
\end{eqnarray}
where $(a)$ follows from a standard contraction property of the Euclidean projection operator. Taking expectation with respect to $\psi_{t}$ and conditioning on $Z_{t}$, we follow (\ref{eq:sba1}) and (\ref{eq:sba2}) to obtain\vspace{-0.3cm}
\[
\mathbb{E}\left[\left\|Z_{t+1} - x^{\ast}\right\|^{2}\;|Z_{t}\right]
\;\leq\; \left\|Z_{t} - x^{\ast}\right\|^{2}
 - \frac{2a_{t}}{d}\cdot\left(Z_{t} - x^{\ast}\right)\cdot\nabla f_{t}(Z_{t})
 +\frac{a_{t}^{2}\left(G^{2}+ \sigma^{2}\right)}{h_{t}^{2}}
 +\frac{2Ga_{t}h_{t}}{\sqrt{d}}\cdot\left\|Z_{t} - x^{\ast}\right\|.\vspace{-0.3cm}
\]
Taking another expectation, with respect to $Z_{t}$, we get\vspace{-0.3cm}
\[
\mathbb{E}\left[\left\|Z_{t+1} - x^{\ast}\right\|^{2}\right]
 \;\leq\; \mathbb{E}\left[\left\|Z_{t} - x^{\ast}\right\|^{2}\right]
 - \frac{2a_{t}}{d}\cdot\mathbb{E}\left[\left(Z_{t} - x^{\ast}\right)\cdot\nabla f_{t}(Z_{t})\right]
 +\frac{a_{t}^{2}\left(G^{2}+ \sigma^{2}\right)}{h_{t}^{2}}
 + \frac{2Ga_{t}h_{t}}{\sqrt{d}}\cdot\mathbb{E}\left\|Z_{t} - x^{\ast}\right\|,\vspace{-0.1cm}
\]
and therefore, fixing some $\gamma >0$, we have for all $t\in\left\{1,\ldots,T-1\right\}$:\vspace{-0.5cm}
\begin{eqnarray}
\mathbb{E}\left[\left(Z_{t} - x^{\ast}\right)\cdot\nabla f_{t}(Z_{t})\right]
&\leq& \frac{d}{2a_{t}}\left(\mathbb{E}\left[\left\|Z_{t} - x^{\ast}\right\|^{2}\right] - \mathbb{E}\left[\left\|Z_{t+1} - x^{\ast}\right\|^{2}\right]\right)
+ \frac{\left(G^{2}+ \sigma^{2}\right)a_{t}d}{2h_{t}^{2}}\notag\\
& &+ \gamma \cdot \frac{1}{\gamma}\cdot Gh_{t}\sqrt{d}\cdot\mathbb{E}\left\|Z_{t} - x^{\ast}\right\|\notag\\
\label{eq:sba4}
&\stackrel{(a)} \leq& \frac{d}{2a_{t}}\left(\mathbb{E}\left[\left\|Z_{t} - x^{\ast}\right\|^{2}\right] - \mathbb{E}\left[\left\|Z_{t+1} - x^{\ast}\right\|^{2}\right]\right)
+ \frac{\left(G^{2}+ \sigma^{2}\right)a_{t}d}{2h_{t}^{2}}\notag\\
& & +\; \frac{\gamma^{2}}{2}\cdot\mathbb{E}\left[\left\|Z_{t} - x^{\ast}\right\|^{2}\right] +\frac{G^{2}h_{t}^{2}d}{2\gamma^{2}},
\end{eqnarray}
where $(a)$ holds by $ab \leq \left(a^{2} + b^{2}\right)/2$, and by Jensen's inequality. In addition, one has for any $t\in \mathcal{T}$:\vspace{-0.5cm}
\begin{eqnarray}
\mathbb{E}\left[f_{t}(X_{t})\right] &=& \mathbb{E}\left[\mathbb{E}\left[f_{t}(X_{t})|Z_{t}\right]\right]
\;=\; \mathbb{E}\left[\frac{1}{2}\left(f_{t}(Z_{t}+h_{t}) + f_{t}(Z_{t}-h_{t})\right)\right]\notag\\
&\leq& \frac{1}{2}\mathbb{E}\left[2f_{t}(Z_{t}) + h_{t}\left(\nabla f_{t}(Z_{t} + h_{t}) - \nabla f_{t}(Z_{t} - h_{t})\right) - Hh_{t}^{2}  \right]\notag\\
&\leq& \mathbb{E}\left[f_{t}(Z_{t}) + \frac{1}{2}Hh_{t}^{2}\right].
\end{eqnarray}
The regret with respect to the single best action is:\vspace{-0.1cm}
\begin{eqnarray}
\sum_{t=1}^{T}\mathbb{E}^{\mathcal{A}}\left[f_{t}(X_{t}) - f_{t}(x^{\ast})\right]
&\leq& 2G + \sum_{t=1}^{T-1}\mathbb{E}^{\pi}\left[f_{t}(Z_{t}) - f_{t}(x^{\ast}) + \frac{1}{2}Hh_{t}^{2}\right]\notag\\
&\stackrel{(a)}\leq& 2G + \sum_{t=1}^{T-1}\mathbb{E}\left[\nabla f_{t}(Z_{t})\cdot\left(Z_{t} - x^{\ast}\right) - \frac{1}{2}H\left\|Z_{t} - x^{\ast}\right\|^{2} + \frac{1}{2}Hh_{t}^{2}\right]\notag\\
&\stackrel{(b)}\leq& \mathbb{E}\left[\sum_{t=1}^{T-1}\left(\frac{d}{2a_{t}}\left(\left\|Z_{t} - x^{\ast}\right\|^{2} - \left\|Z_{t+1} - x^{\ast}\right\|^{2}\right)
+ \frac{\left(\gamma^{2} - H\right)}{2}\cdot \left\|Z_{t} - x^{\ast}\right\|^{2}\right)\right]\notag\\
& &+\; 2G + \frac{\left(G^{2}+\sigma^{2}\right)}{2}\sum_{t=1}^{T-1}\left(\frac{a_{t}d}{h_{t}^{2}}  + \frac{h_{t}^{2}d}{\gamma^{2}}\right) + \frac{H}{2}\sum_{t=1}^{T-1}h_{t}^{2}\notag\\
&\stackrel{(c)}=& \frac{1}{2}\sum_{t=2}^{T}\mathbb{E}\left[\left\|Z_{t} - x^{\ast}\right\|^{2}\right]\underbrace{\left(\frac{d}{a_{t}} - \frac{d}{a_{t-1}} + \left(\gamma^{2} - H\right) \right)}_{I_{t}}
 +\mathbb{E}\left[\left\|Z_{1}-x^{\ast}\right\|^{2}\right]\underbrace{\left(\frac{d}{2a_{1}} + \frac{\gamma^{2}- H}{2}\right)}_{I_{1}}\notag\\
 & &- \mathbb{E}\left[\left\|Z_{T}-x^{\ast}\right\|^{2}\right]\frac{d}{2a_{T-1}}
+ 2G + \frac{\left(G^{2}+\sigma^{2}\right)}{2}\sum_{t=1}^{T-1}\left(\frac{a_{t}d}{h_{t}^{2}} + \frac{h_{t}^{2}d}{\gamma^{2}}\right)
+ \frac{H}{2}\sum_{t=1}^{T-1}h_{t}^{2},\notag
\end{eqnarray}
where $(a)$ holds by (\ref{eq:sba3}), $(b)$ holds by (\ref{eq:sba4}), and $(c)$ holds by rearranging the summation. By selecting $\gamma^{2} = \frac{H}{2}$, $a_{t} = \frac{d}{(H - \gamma^{2})t}$, and $h_{t} = \delta_{t} = a_{t}^{1/4}$, we have $I_{t} = 0$ for all $t\in \mathcal{T}$, and:\vspace{-0.2cm}
\[
\mathbb{E}^{\mathcal{A}}\left[\sum_{t=1}^{T}f_{t}(X_{t})\right] - \inf_{x\in\mathcal{X}}\left\{\sum_{t=1}^{T}f_{t}(x)\right\}
\;\leq\; 2G + \frac{\left(G^{2}+\sigma^{2} + H\right)d^{3/2}}{\sqrt{2H}}\cdot\sqrt{T}.\vspace{-0.2cm}
\]
Since the above holds for any $f\in \mathcal{F}_{s}$, we conclude that\vspace{-0.2cm}
\[
\mathcal{G}^{\mathcal{A}}_{\phi^{(0)}}\left(\mathcal{F}_{s},T\right)
\;\leq\; 2G +\frac{\left(G^{2}+\sigma^{2} + H\right)d^{3/2}}{\sqrt{2H}}\cdot\sqrt{T}.\vspace{-0.2cm}
\]
Finally, we consider the case in which there exists at least one time epoch $t$ such that $\delta_{t} > \nu$. Then, for any such time epoch we select $h'_{t} = \delta'_{t} =\min\left\{\nu, \delta_{t}\right\}$. We note that the sequence $\left\{\delta_{t}\right\}$ is converging to $0$, and therefore for any number $\nu$ there is some epoch $t_{\nu}$, independent of $T$, such that $\delta_{t}\leq \nu$ for any $t\geq t_{\nu}$. Therefore there can be no more than $t_{\nu}$ such epochs. In particular, it follows that such a case could add to the regret above no more than a constant (independent of $T$), that depends solely on $\nu$, the dimension $d$, and the second derivative bound $H$. This concludes the proof.\qed

\begin{lemma}\label{lem:sbaGrad}\textbf{\textup{(Performance of OGD in the adversarial setting)}}
Consider the feedback structure $\phi =~\phi^{(1)}$. Let $\mathcal{A}$ be the OGD algorithm given in \S4, with the selection $\eta_{t+1} = 1/Ht$ for $t=1,\ldots T-1$. Then, there exists a constant $\bar{C}$, independent of $T$ such that for any $T\geq 1$,\vspace{-0.3cm}
\[
\mathcal{G}^{\mathcal{A}}_{\phi}\left(\mathcal{F}_{s},T\right)
\;\leq\; \bar{C}\log T.\vspace{-0.1cm}
\]
\end{lemma}
\textbf{Proof.} We adapt the proof of Theorem 1 in \cite{Haz-Aga-Kal2007} to the feedback $\phi^{(1)}$. Fix $\phi = \phi^{(1)}$, $T\geq 1$, and $f\in\mathcal{F}_{s}$. Selecting $\eta_{t} = 1/Ht$ for any $t=2,\ldots T$, one has that for any $x\in\mathcal{X}$ and $f_{t}$,\vspace{-0.3cm}
\begin{equation}\label{eq:haz0}
\mathbb{E}\left[\phi^{(1)}_{t}\left(X_{t},f_{t}\right)\;|X_{t} = x\right] \;=\; \nabla f_{t}(x), \quad\text{ and }\quad \mathbb{E}\left[\left\|\phi^{(1)}_{t}\left(x,f_{t}\right)\right\|^{2}\right] \;\leq\; G^{2} + \sigma^{2}, \vspace{-0.3cm}
\end{equation}
for some $\sigma \geq 0$. Define $x^{\ast}$ as the single best action in hindsight: $x^{\ast} = \argmin_{x\in \mathcal{X}}\left\{\sum_{t=1}^{T}f_{t}(x)\right\}$. Then, by a Taylor expansion, for any $x\in \mathcal{X}$ there is a point $\tilde{x}$ on the segment between $x$ and $x^{\ast}$ such that:\vspace{-0.3cm}
\begin{eqnarray}
f_{t}(x^{\ast}) &=& f_{t}(x) + \nabla f_{t}(x)\cdot(x^{\ast} - x) + \frac{1}{2}(x^{\ast} - x)\cdot\nabla^{2}f_{t}(\tilde{x})(x^{\ast} - x)\notag\\
&\stackrel{(a)}\geq& f_{t}(x) + \nabla f_{t}(x)\cdot(x^{\ast} - x) + \frac{H}{2}\left\|x^{\ast} - x\right\|^{2},\notag
\end{eqnarray}
for any $t\in \mathcal{T}$, where $(a)$ holds by (\ref{eq:H}). Substituting $X_{t}$ in the above and taking expectation with respect to $X_{t}$, one has:\vspace{-0.3cm}
\begin{equation}\label{eq:haz1}
\mathbb{E}\left[f_{t}(X_{t})\right] - f_{t}(x^{\ast})
\;\leq\; \mathbb{E}\left[\nabla f_{t}(X_{t})\cdot(X_{t} - x^{\ast})\right] - \frac{H}{2}\mathbb{E}\left\|x^{\ast} - X_{t}\right\|^{2},\vspace{-0.3cm}
\end{equation}
for any $t\in \mathcal{T}$. By the OGD step,\vspace{-0.3cm}
\[
\left\|X_{t+1} - x^{\ast}\right\|^{2}
 \;=\; \left\|P_{\mathcal{X}}\left(X_{t} - \eta_{t+1}\phi^{(1)}_{t}\left(X_{t},f_{t}\right)\right) - x^{\ast}\right\|^{2}
 \;\stackrel{(a)} \leq\; \left\|X_{t} - \eta_{t+1}\phi^{(1)}_{t}\left(X_{t},f_{t}\right) - x^{\ast}\right\|^{2},\vspace{-0.3cm}
\]
where $(a)$ follows from a standard contraction property of the Euclidean projection operator. Taking expectation with respect to $X_{t}$, one has:\vspace{-0.3cm}
\begin{eqnarray}
\mathbb{E}\left\|X_{t+1} - x^{\ast}\right\|^{2}
&\leq& \mathbb{E}\left\|X_{t} - x^{\ast}\right\|^{2}
 + \eta_{t+1}^{2}\mathbb{E}\left\|\phi^{(1)}_{t}\left(X_{t},f_{t}\right)\right\|^{2}
 - 2\eta_{t+1} \mathbb{E}\left[\left(\phi^{(1)}_{t}\left(X_{t},f_{t}\right)\right)\cdot(X_{t} - x^{\ast})\right]\notag\\
&\stackrel{(a)}\leq& \mathbb{E}\left\|X_{t} - x^{\ast}\right\|^{2}
 + \eta_{t+1}^{2}\left(G^{2} + \sigma^{2}\right)
 - 2\eta_{t+1} \mathbb{E}\left[\left(\nabla f_{t}(X_{t})\right)\cdot(X_{t} - x^{\ast})\right],\notag
\end{eqnarray}
where $(a)$ follows from (\ref{eq:haz0}). Therefore, for any $t\in \mathcal{T}$, we get:\vspace{-0.3cm}
\begin{equation}\label{eq:haz2}
\mathbb{E}\left[\nabla f_{t}(X_{t})\cdot(X_{t} - x^{\ast})\right]
\;\leq\; \frac{\mathbb{E}\left\|X_{t} - x^{\ast}\right\|^{2} - \mathbb{E}\left\|X_{t+1} - x^{\ast}\right\|^{2}}{2\eta_{t+1}}
+ \frac{\eta_{t+1}}{2}\left(G^{2} + \sigma^{2}\right).\vspace{-0.3cm}
\end{equation}
Summing (\ref{eq:haz1}) over the horizon and using (\ref{eq:haz2}), one has:\vspace{-0.3cm}
\begin{eqnarray}\label{eq:haz3}
\sum_{t=1}^{T}\left( \mathbb{E}\left[f_{t}(X_{t})\right] - f_{t}(x^{\ast}) \right)
 &\leq& \frac{1}{2}\sum_{t=2}^{T}\mathbb{E}\left\|X_{t} - x^{\ast}\right\|^{2}\left(\frac{1}{\eta_{t+1}} - \frac{1}{\eta_{t}} -H \right)\notag\\
 & & +\; \frac{1}{2}\mathbb{E}\left\|X_{1} - x^{\ast}\right\|^{2}\left(\frac{1}{\eta_{2}} - \frac{H}{2}\right)
  - \frac{1}{2}\mathbb{E}\left\|X_{T+1} - x^{\ast}\right\|^{2}\left(\frac{1}{\eta_{T+1}} + \frac{H}{2}\right)\notag\\
  & & +\; \frac{\left(G^{2} + \sigma^{2}\right)}{2}\sum_{t=1}^{T}\eta_{t+1}\\
 &\stackrel{(a)} \leq & \frac{\left(G^{2} + \sigma^{2}\right)}{2}\sum_{t=1}^{T}\frac{1}{Ht}
 \;\leq\; \frac{\left(G^{2} + \sigma^{2}\right)}{2H}\left(1+\log T\right),\notag
\end{eqnarray}
where $(a)$ holds using $\eta_{t} = 1/Ht$. Since the above holds for any sequence of functions in $\mathcal{F}_{s}$ we have that\vspace{-0.2cm}
\[
\mathcal{G}^{\mathcal{A}}_{\phi^{(1)}}\left(\mathcal{F}_{s},T\right)
\;\leq\; \frac{\left(G^{2} + \sigma^{2}\right)}{2H}\left(1+\log T\right),\vspace{-0.2cm}
\]
which concludes the proof. \qed\vspace{-0.2cm}

\subsection{Lower bounds}\vspace{-0.2cm}
The last two results of this section, Lemma \ref{lem:lowerEGS} and Lemma \ref{lem:lowerFlax}, establish lower bounds on the best achievable performance in the adversarial setting, under the structures $(\mathcal{F}_{s},\phi^{(0)})$, and $(\mathcal{F},\phi^{(1)})$, respectively. Lemma \ref{lem:lowerEGS} provides a lower bound that (together with the upper bound in Lemma \ref{lem:sbaRegret}) establishes that the EGS algorithm is rate optimal in a setting with strongly convex cost functions and noisy cost observations. Lemma \ref{lem:lowerFlax} provides a lower bound that matches the upper bound in Lemma 3.1 in \cite{FLA2005}, establishing that the OGD algorithm (with a careful selection of step-sizes), is rate optimal in a setting with general convex cost functions and noisy gradient observations.

\begin{lemma}\label{lem:lowerEGS}
Let Assumption~2 hold. Then, there exists a constant $C$, independent of $T$ such that for any online algorithm $\mathcal{A}\in\mathcal{P}_{\phi^{(0)}}$ and for all $T\geq 1$:\vspace{-0.3cm}
\[
\mathcal{G}^{\mathcal{A}}_{\phi^{(0)}}\left(\mathcal{F}_{s},T\right)
\;\geq\; C\sqrt{T}.\vspace{-0.1cm}
\]
\end{lemma}
\textbf{Proof.} Let $\mathcal{X} = [0,1]$. Consider the quadratic functions $f^{1}$ and $f^{2}$ in (\ref{eq:sconvFuns}), used in the proof of Theorems~\ref{thm:strictGrad} and~\ref{thm:strictCost}. (note that $\delta$ will be selected differently). Fix some algorithm $\mathcal{A}\in\mathcal{P}_{\phi^{(0)}}$. Let $\tilde{f}$ be a random sequence where in the beginning of the horizon nature draws (according to a uniform discrete distribution) a cost function from $\left\{f^{1},f^{2}\right\}$, and applies it throughout the horizon. Taking expectation over the random sequence $\tilde{f}$ one has\vspace{-0.3cm}
\[
\mathcal{G}^{\mathcal{A}}_{\phi^{(0)}}(\mathcal{F}_{s},T)
\;\geq\; \frac{1}{2}\mathbb{E}_{f^{1}}^{\mathcal{A}}\left[\sum_{t=1}^{T}\left(f^{1}(X_{t}) - f^{1}(x_{1}^{\ast})\right)\right]
+ \frac{1}{2}\mathbb{E}_{f^{2}}^{\mathcal{A}}\left[\sum_{t=1}^{T}\left(f^{2}(X_{t}) - f^{2}(x_{2}^{\ast})\right)\right],\vspace{-0.3cm}
\]
where the inequality follows as in step 3 of the proof of theorem \ref{thm:weakGrad-low}. In the following we use notation described at the proof of Theorem \ref{thm:strictCost}, for the online algorithm $\mathcal{A}$. We start by bounding the Kullback-Leibler divergence between $\mathbb{P}^{\mathcal{A},\tau}_{f^{1}}$ and $\mathbb{P}^{\mathcal{A},\tau}_{f^{2}}$ for all $\tau\in \mathcal{T}$:\vspace{-0.3cm}
\begin{eqnarray}
\mathcal{K}\left(\mathbb{P}^{\mathcal{A},T}_{f^{1}}\|\mathbb{P}^{\mathcal{A},T}_{f^{2}}\right)
&\stackrel{(a)} \leq & \tilde{C}\mathbb{E}^{\mathcal{A}}_{f^{1}}\left[\sum_{t=1}^{T}\left(f^{1}(X_{t}) - f^{2}(X_{t})\right)^{2}\right]
\;=\; \tilde{C}\mathbb{E}^{\mathcal{A}}_{f^{1}}\left[\sum_{t=1}^{T}\left(\delta X_{t} - \frac{\delta}{2}\right)^{2}\right]\notag\\
&=& \tilde{C}\mathbb{E}_{f^{1}}^{\mathcal{A}}\left[\delta^{2}\sum_{t=1}^{T}\left(X_{t} - x_{1}^{\ast}\right)^{2}\right]\notag\\
&\stackrel{(b)} =& \tilde{C}\mathbb{E}_{f^{1}}^{\mathcal{A}}\left[2\delta^{2}\sum_{t=1}^{T}\left(f^{1}(X_{t}) - f^{1}(x_{1}^{\ast})\right)\right]
\label{eq:boundOnKLLem1}
\;\stackrel{(c)} \leq\; 4\tilde{C}\delta^{2}\mathcal{G}^{\mathcal{A}}_{\phi^{(0)}}(\mathcal{F}_{s},T),
\end{eqnarray}
where: $(a)$ follows from Lemma \ref{lem:KL2}; $(b)$ holds since\vspace{-0.3cm}
\[
f^{1}(x) - f^{1}(x_{1}^{\ast})
\;=\; \nabla f^{1}(x_{1}^{\ast})\cdot(x-x_{1}^{\ast}) + \frac{1}{2}\cdot\nabla f^{1}(x_{1}^{\ast})\cdot(x-x_{1}^{\ast})^{2}
\;=\; \frac{1}{2}(x-x_{1}^{\ast})^{2}\vspace{-0.3cm}
\]
for any $x\in \mathcal{X}$; and $(c)$ holds by\vspace{-0.4cm}
\begin{eqnarray}
\mathcal{G}^{\mathcal{A}}_{\phi^{(0)}}(\mathcal{F}_{s},T)
&\geq& \frac{1}{2}\mathbb{E}_{f^{1}}^{\mathcal{A}}\left[\sum_{t=1}^{T}\left(f^{1}(X_{t}) - f^{1}(x_{1}^{\ast})\right)\right]
+ \frac{1}{2}\mathbb{E}_{f^{2}}^{\mathcal{A}}\left[\sum_{t=1}^{T}\left(f^{2}(X_{t}) - f^{2}(x_{2}^{\ast})\right)\right]\notag\\
&\geq& \frac{1}{2}\mathbb{E}_{f^{1}}^{\mathcal{A}}\left[\sum_{t=1}^{T}\left(f^{1}(X_{t}) - f^{1}(x_{1}^{\ast})\right)\right].
\end{eqnarray}
Therefore, for any $x_{0}\in \mathcal{X}$, by Lemma \ref{lem:tsy} with $\varphi_{t} = \mathbbm{1}\{X_{t} > x_{0}\}$, we have:\vspace{-0.3cm}
\begin{equation}\label{eq:phierrorLem1}
\max\left\{\mathbb{P}^{\mathcal{A}}_{f^{1}}\left\{X_{\tau} > x_{0}\right\}, \mathbb{P}^{\mathcal{A}}_{f^{2}}\left\{X_{\tau} \leq x_{0}\right\}\right\}
 \;\geq\; \frac{1}{4}\exp\left\{-4\tilde{C}\delta^{2}\mathcal{G}^{\mathcal{A}}_{\phi^{(0)}}(\mathcal{F}_{s},T)\right\} \quad\quad \text{ for all } \tau\in \mathcal{T}.\vspace{-0.2cm}
\end{equation}
Set $x_{0} = \frac{1}{2}\left(x_{1}^{\ast}+x_{2}^{\ast}\right) = 1/2 + \delta/4$. Then, following step 3 in the proof of Theorem \ref{thm:strictCost}, one has:\vspace{-0.5cm}
\begin{eqnarray}
\mathcal{G}^{\mathcal{A}}_{\phi^{(0)}}(\mathcal{F}_{s},T)
&\geq& \frac{1}{2}\sum_{t=1}^{T}\left(f^{1}(x_{0}) - f^{1}(x_{1}^{\ast})\right)\mathbb{P}_{f^{1}}^{\mathcal{A}}\left\{X_{t}> x_{0}\right\}
+ \frac{1}{2}\sum_{t=1}^{T}\left(f^{2}(x_{0}) - f^{2}(x_{2}^{\ast})\right)\mathbb{P}_{f^{2}}^{\mathcal{A}}\left\{X_{t}\leq x_{0}\right\}\notag\\
&\geq& \frac{\delta^{2}}{16}\sum_{t=1}^{T}\left(\mathbb{P}_{f^{1}}^{\mathcal{A}}\left\{X_{t}> x_{0}\right\}
+ \mathbb{P}_{f^{2}}^{\mathcal{A}}\left\{X_{t}\leq x_{0}\right\}  \right)\notag\\
&\geq&\frac{\delta^{2}}{16}\sum_{t=1}^{T}\max\left\{\mathbb{P}^{\mathcal{A}}_{f^{1}}\left\{X_{t} > x_{0}\right\}, \mathbb{P}^{\mathcal{A}}_{f^{2}}\left\{X_{t} \leq x_{0}\right\}\right\}\notag\\
&\stackrel{(a)}\geq& \frac{\delta^{2}}{16}\sum_{t=1}^{T}\frac{1}{4}\exp\left\{-4\tilde{C}\delta^{2}\mathcal{G}^{\mathcal{A}}_{\phi^{(0)}}(\mathcal{F}_{s},T)\right\}
 \;=\; \frac{\delta^{2}T}{16}\exp\left\{-4\tilde{C}\delta^{2}\mathcal{G}^{\mathcal{A}}_{\phi^{(0)}}(\mathcal{F}_{s},T)\right\}\notag
\end{eqnarray}
where $(a)$ holds by (\ref{eq:phierrorLem1}). Set $\delta = \left(\frac{4}{\tilde{C}T}\right)^{1/4}$. Then, one has for $\beta = 8\sqrt{\tilde{C}/T}$:\vspace{-0.3cm}
\begin{equation}\label{eq:inequalityLem4}
\beta \mathcal{G}^{\mathcal{A}}_{\phi^{(0)}}(\mathcal{F}_{s},T)
\;\geq\; \exp\left\{-\beta \mathcal{G}^{\mathcal{A}}_{\phi^{(0)}}(\mathcal{F}_{s},T)\right\}.\vspace{-0.3cm}
\end{equation}
Let $y_{0}$ be the unique solution to the equation $y = \exp\left\{-y\right\}$. Then, (\ref{eq:inequalityLem4}) implies $\beta \mathcal{G}^{\mathcal{A}}_{\phi^{(0)}}(\mathcal{S},T) \geq y_{0}$. In particular, since $y_{0}>1/2$ this implies\vspace{-0.3cm}
\[
\mathcal{G}^{\mathcal{A}}_{\phi^{(0)}}(\mathcal{F}_{s},T)
\;\geq\; 1/\left(2\beta\right)
 \;=\; \frac{1}{16\sqrt{\tilde{C}}}\cdot \sqrt{T}.\vspace{-0.5cm}
\]
This concludes the proof. \qed

\begin{lemma}\label{lem:lowerFlax}
Let Assumption~1 hold. Then, there exists a constant $C$, independent of $T$, such that for any online algorithm $\mathcal{A}\in\mathcal{P}_{\phi^{(1)}}$ and for all $T\geq 1$:\vspace{-0.3cm}
\[
\mathcal{G}^{\mathcal{A}}_{\phi^{(1)}}\left(\mathcal{F},T\right)
\;\geq\; C\sqrt{T}.\vspace{-0.1cm}
\]
\end{lemma}
\textbf{Proof.} Fix $T\geq 1$. Let $\mathcal{X} = [0,1]$, and consider functions $f^{1}$ and $f^{2}$ that are given in (\ref{eq:convFuns}), and used in the proof of Theorem~\ref{thm:weakGrad-low} (note that $\delta$ will be selected differently). Let $\tilde{f}$ be a random sequence of cost functions, where in the beginning of the time horizon nature draws (from a uniform discrete distribution) a function from $\left\{f^{1},f^{2}\right\}$, and applies it throughout the horizon.

Fix $\mathcal{A}\in\mathcal{P}_{\phi^{(1)}}$. In the following we use notation described in the proof of Theorem~\ref{thm:weakGrad-low}, as well as in Lemma \ref{lem:lowerEGS}. Set $\delta =~1/\sqrt{16\tilde{C}T}$, where $\tilde{C}$ is the constant that appears in Assumption~1. Then:\vspace{-0.2cm}
\begin{eqnarray}
\mathcal{K}\left(\mathbb{P}^{\mathcal{A},T}_{f^{1}}\|\mathbb{P}^{\mathcal{A},T}_{f^{2}}\right)
&\stackrel{(a)} \leq& \tilde{C}\mathbb{E}_{f^{1}}^{\mathcal{A}}\left[\sum_{t=1}^{T}\left(\nabla f^{1}(X_{t}) - \nabla f^{2}(X_{t})\right)^{2}\right]\notag\\
\label{eq:boundOnKLGradLem2}
&=& \tilde{C}\mathbb{E}_{f^{1}}^{\mathcal{A}}\left[\sum_{t=1}^{T}16\delta^{2}X_{t}^{2}\right]
\;\leq\; 16\tilde{C}T\delta^{2}
\;\stackrel{(b)}\leq\; 1,
\end{eqnarray}
where $(a)$ follows from Lemma \ref{lem:KL}, and $(b)$ holds by $\delta =~1/\sqrt{16\tilde{C}T}$. Since $\mathcal{K}(\mathbb{P}^{\mathcal{A},\tau}_{1}\|\mathbb{P}^{\mathcal{A},\tau}_{2})$ is non-decreasing in $\tau$ throughout the horizon, we deduce that the Kullback-Leibler divergence is bounded by~$1$ throughout the horizon. Therefore, for any $x_0\in \mathcal{X}$, by Lemma \ref{lem:tsy} with $\varphi_{\tau} = \mathbbm{1}\{X_{\tau} \leq x_0\}$ and $\beta = 1$, one has:\vspace{-0.1cm}
\begin{equation}\label{eq:phierrorGradLem2}
\max\left\{\mathbb{P}_{f^{1}}^{\mathcal{A}}\left\{X_{\tau} \leq x_0\right\}, \mathbb{P}_{f^{2}}^{\mathcal{A}}\left\{X_{t} > x_0\right\}\right\} \geq \frac{1}{4e} \quad\quad \text{ for all } \tau\in \mathcal{T}.\vspace{-0.1cm}
\end{equation}
Set $x_{0} =\frac{1}{2}\left(x_{1}^{\ast}+x_{2}^{\ast}\right) =\frac{1}{2}$. Taking expectation over $\tilde{f}$ and following step 3 in the proof of Theorem~\ref{thm:weakGrad-low}, one has:\vspace{-0.3cm}
\begin{eqnarray}
\mathcal{G}^{\mathcal{A}}_{\phi^{(1)}}\left(\mathcal{F},T\right)
&\geq& \frac{1}{2}\mathbb{E}_{f^{1}}^{\mathcal{A}}\left[\sum_{t=1}^{T}\left(f^{1}(X_{t}) - f^{1}(x_{1}^{\ast})\right)\right]
+ \frac{1}{2}\mathbb{E}_{f^{2}}^{\mathcal{A}}\left[\sum_{t=1}^{T}\left(f^{2}(X_{t}) - f^{2}(x_{2}^{\ast})\right)\right]\notag\\
&\geq& \frac{1}{2}\sum_{t=1}^{T}\left(f^{1}(x_{0}) - f^{1}(x_{1}^{\ast})\right)\mathbb{P}_{f^{1}}^{\mathcal{A}}\left\{X_{t}> x_{0}\right\}
+ \frac{1}{2}\sum_{t=1}^{T}\left(f^{2}(x_{0}) - f^{2}(x_{2}^{\ast})\right)\mathbb{P}_{f^{2}}^{\mathcal{A}}\left\{X_{t}\leq x_{0}\right\}\notag\\
&\geq& \left(\frac{\delta}{4} + \frac{\delta^{2}}{2}\right)\sum_{t=1}^{T}\left(\mathbb{P}_{f^{1}}^{\mathcal{A}}\left\{X_{t}> x_{0}\right\}
+ \mathbb{P}_{f^{2}}^{\mathcal{A}}\left\{X_{t}\leq x_{0}\right\}  \right)\notag\\
&\geq& \left(\frac{\delta}{4} + \frac{\delta^{2}}{2}\right)\sum_{t=1}^{T}\max\left\{\mathbb{P}^{\pi}_{f^{1}}\left\{X_{t} > x_{0}\right\}, \mathbb{P}^{\mathcal{A}}_{f^{2}}\left\{X_{t} \leq x_{0}\right\}\right\}\notag\\
&\stackrel{(a)}\geq&\left(\frac{\delta}{4} + \frac{\delta^{2}}{2}\right)\sum_{t=1}^{T}\frac{1}{4}\exp\left\{-1\right\}
\; \geq \; \frac{\delta T}{16e}
\;\stackrel{(b)} =\; \frac{1}{64e\sqrt{\tilde{C}}}\cdot  \sqrt{T},\notag
\end{eqnarray}
where $(a)$ holds by (\ref{eq:phierrorGradLem2}), and $(b)$ holds by $\delta = 1/\sqrt{16\tilde{C}T}$. This concludes the proof.\qed

\vspace{-0.1cm}
\section{Numerical Results}\label{subsec:numerics}
\vspace{-0.1cm}
We illustrate the upper bounds on the regret by numerical experiments measuring the average regret that is incurred in the presence of various patterns of changing costs, and under different feedback structures and noise. We compare the performance of the restarted OGD and restarted EGS against the performance achieved by applying the respective subroutine without restarting, and with fixed step sizes. We note that policies of fixed step sizes, while having no performance guarantees relative the dynamic oracle, are considered in the SA literature as practical approach to sequential stochastic optimization of general cost functions when these may change; see, e.g., chapter 4 of \cite{Ben-Pri-Met1990}.

\textbf{Variation and feedback.}
We fix $\mathcal{X} = \left[-2,3\right]$ and consider the sequence of quadratic cost functions $f_{t}(x_{t}) =~\frac{x_{t}^{2}}{2} -~b_{t}x_{t} +~1$, where the coefficient $b_{t}$ is time-varying. In particular, for a given horizon length $T\geq 1000$, we let $\tau$ be the random time in which the cost begins to change, drawn from a discrete uniform distribution over $\left\{1,2,\ldots, \left\lfloor T/4 \right\rfloor \right\}$. Then, we consider the following variation patterns: \vspace{-0.2cm}
\begin{displaymath}
b^{shock}_{t} = \left\{
     \begin{array}{lr}
       1 & \text{ if } t\leq \tau \\
       0 & \text{ otherwise}
     \end{array}
   \right.\quad\;\;
b^{decay}_{t} = \left\{
     \begin{array}{lr}
       1 & \text{ if } t\leq \tau \\
       e^{-10(t- \tau)/T} & \text{ otherwise}
     \end{array}
   \right.\quad\;\;
b^{linear}_{t} = \left\{
     \begin{array}{lr}
       1 & \text{ if } t\leq \tau \\
       \frac{T-t}{T-\tau} & \text{ otherwise}
     \end{array}
   \right.
   \vspace{-0.2cm}
\end{displaymath}
for all $t = 1,\ldots,T$, where $\alpha$ is some decay parameter in $\left(0,1\right)$. One may observe that the variation may be bounded by the budget $V_{T} = 1$ in all the considered patterns. Let $\left\{\varepsilon_{t}\right\}_{t}^{T}$ be a sequence of independent normal random variables with zero mean and standard deviation $\sigma$. We consider the case of noisy access to the cost, where $\phi_{t}^{(0)}(x_{t},f_{t}) = f_{t}(x_{t}) + \varepsilon_{t}$ for each $t \in \mathcal{T}$, and the case of noisy access to the gradient, where $\phi_{t}^{(1)}(x_{t},f_{t}) = \nabla f_{t}(x_{t}) + \varepsilon_{t}$ for each $t \in \mathcal{T}$. At each epoch $t\in \mathcal{T}$ an action $X_{t}\in \mathcal{X}$ is selected, and then the expected cost $f_{t}(X_{t})$ is incurred, and a feedback $\phi_{t}(x_{t},f_{t})$ is observed.

\textbf{Policies and performance.} In the case of noisy gradient access we use a version of the restarted OGD policy that is considered in \S4 and \S5.1, where the first action in batch $j\geq 2$ is obtained by taking a gradient step from the last action of batch $j-1$; in other words, only the sequence of gradient steps $\left\{\eta_{t}\right\}$ is restarted. Similarly, when the feedback consists of noisy cost observations we use a variation of the restarted EGS policy that is considered in \S5.2, where only the sequences $\left\{a_{t}\right\}$, $\left\{h_{t}\right\}$, and $\left\{\delta_{t}\right\}$ are restarted.
Given the actions $\left\{X_{t}\right\}$ generated throughout $T$ epoch by a policy $\pi$ under feedback $\phi$ and a sequence $f$ of cost functions, we measure the regret relative to the dynamic oracle $R^{\pi}_{\phi}(f,T) = \sum_{t=1}^{T}\left(f_{t}(X_{t}) - f_{t}(x^{\ast}_{t})\right)$.
We denote the relative loss (in percentage) relative to the dynamic oracle by $L^{\pi}_{\phi}(f,T) = 100\cdot R^{\pi}_{\phi}(f,T) / \left(\sum_{t=1}^{T}f_{t}(x^{\ast}_{t})\right)$. We refer to policy (OGD or EGS) as ``non-restarted" when applied without restarting, and as ``fixed step size of $a$" when these apply a fixed (non-updating) step size $a$ ($\eta_{t} = a$ in the OGD; $a_{t}=a$ and $\delta_{t} = h_{t} = a^{1/4}$ in the EGS).

For each of the considered variation patterns, and for each value of $\sigma\in\left\{0.1,0.3,1\right\}$, we simulated the action paths of the policy (restated EGS for feedback $\phi^{(0)}$, restarted OGD for feedback $\phi^{(1)}$) for various values of $T\in \left\{1000, 5000,\ldots, 37000\right\}$, replicating each instance $10^{3}$ times and calculating the average regret and relative loss. Assuming the structure $R^{\pi}_{\phi}(f,T) =~cT^{\alpha}$, we estimate the coefficients $c$ and $\alpha$ from fitting the log-regret as a function of log-time.

\textbf{Results and discussion.} Table \ref{tab:OGD} details the estimated coefficients $c$ and $\alpha$ under the considered variation patterns and $\sigma$ values, for the restarted OGD and feedback $\phi^{(1)}$. Table \ref{tab:OGD} also includes the average loss (in $\%$) relative to the dynamic oracle for two representative values of $T$. Table~\ref{tab:EGS} includes the respective results for the restarted EGS and feedback $\phi^{(0)}$. In all the linear fits we observed $R^{2}>0.98$, and the standard error of the percentage loss was always below $5\%$ of the policy's average performance.

\begin{table}[!ht]\small
\begin{center}
\begin{tabular}{l*{12}{c}}
\hline
Variation pattern
& & \multicolumn{3}{c}{$b_{t}^{shock}$}
& & \multicolumn{3}{c}{$b_{t}^{decay}$}
& & \multicolumn{3}{c}{$b_{t}^{linear}$}  \\
\cmidrule(r){3-5} \cmidrule(r){7-9} \cmidrule(r){11-13}
$\sigma$                           & & 0.1  & 0.3  & 1    & &  0.1  & 0.3  & 1    & &  0.1  & 0.3  & 1    \\
\hline
$\alpha$                           & & \;0.54\; & 0.54 & \;0.54\; & &  \;0.47\; & 0.47 & \;0.52\;  & &  \;0.47\;  & 0.51  & \;0.54\; \\
$c$                                & & 0.26 & 0.32 & 1.02 & &  0.14 & 0.16 & 0.89  & &  0.05  & 0.09  & 0.94 \\
\hline
$L^{\pi}_{\phi}(f,T)$, $T=5000$:\\
Restarted                 & & 0.56  & 0.68 & 2.02 & &  0.05 & 0.17 & 1.56  & & 0.03   & 0.17  & \;1.78\; \\
Non-restarted          & & 3.02  & 3.02 & 3.08 & &  4.94 & 4.95 & 5.01  & & 5.81   & 5.81  & 5.86  \\
Fixed step size of $0.1$   & & 0.09  & 0.31 & 2.80 & & 0.05  & 0.27 & 2.82  & & 0.05   & 0.31  & 3.26  \\
Fixed step size of $0.01$ & & 0.64  & 0.66 & 0.89 & & 0.22  & 0.25 & 0.49  & & 0.18   & 0.21  & 0.49  \\
Fixed step size of $0.001$ & & 2.87  & 2.87 & 2.89 & & 2.41  & 2.41 & 2.43  & & 2.44   & 2.45  & 2.47  \\
\hline
$L^{\pi}_{\phi}(f,T)$, $T=25000$:\\
Restarted                  & & 0.26  & 0.32 & 0.94 & &  0.02 & 0.07 & \;0.71\; & &  \;0.01\; & \;0.08\; & 0.82 \\
Non-restarted         & & 1.49  & 1.49 & 1.50 & &  3.49 & 3.50 &  3.51    & & 5.41      & 5.41     &  5.41 \\
Fixed step size of $0.1$ & & 0.04  & 0.25 & 2.70 & &  0.03 & 0.25 &  2.75    & & 0.04   & 0.29  & 3.21  \\
Fixed step size of $0.01$ & & 0.13  & 0.15 & 0.38 & &  0.03 & 0.06 &  0.29    & & 0.03   & 0.06  & 0.34  \\
Fixed step size of $0.001$ & & 0.85  & 0.85 & 0.88 & &  0.43 & 0.44 &  0.46    & & 0.37   & 0.38  & 0.41  \\
\hline\vspace{-0.6cm}
\end{tabular}
\caption{\small \textbf{Performance of restarted OGD under noisy gradient observations.}\vspace{-0.4cm}}
\label{tab:OGD}
\end{center}
\end{table}

\begin{table}[!ht]\small
\begin{center}
\begin{tabular}{l*{12}{c}}
\hline
Variation pattern
& & \multicolumn{3}{c}{$b_{t}^{shock}$}
& & \multicolumn{3}{c}{$b_{t}^{decay}$}
& & \multicolumn{3}{c}{$b_{t}^{linear}$}  \\
\cmidrule(r){3-5} \cmidrule(r){7-9} \cmidrule(r){11-13}
$\sigma$                            & & 0.1   & 0.3   & 1      & &  0.1   & 0.3   & 1     & &  0.1   & 0.3   & 1     \\
\hline
$\alpha$                            & & 0.68  & 0.68  & 0.68   & &  0.67  & 0.67  & 0.68  & &  0.67  & 0.67  & 0.68  \\
$c$                                 & & 2.22  & 2.28  & 2.84   & &  2.13  & 2.18  & 2.88  & &  2.09  & 2.16  & 2.83  \\
\hline
$L^{\pi}_{\phi}(f,T)$, $T=5000$:\\
Restarted                    & & 14.45 & 14.82 & 19.02  & &  14.42 & 14.89 & 19.01 & &  15.52 & 16.06  & 21.20 \\
Non-restarted                & & 37.71 & 38.03 & 41.25  & &  29.96 & 31.30 & 33.35 & & 28.49  & 29.46  & 39.35 \\
Fixed step size of $0.1$     & & 26.58 & 27.33 & 34.99  & & 26.78  & 27.53 & 35.62 & &  28.49 & 29.46  & 39.35  \\
Fixed step size of $0.01$    & & 8.48  & 8.70  & 11.22  & & 8.10   & 8.33  & 10.87 & &  8.48  & 8.86   & 11.82  \\
Fixed step size of $0.001$   & & 5.22  & 5.30  & 5.95   & & 4.94   & 4.81  & 5.53  & &   4.96 & 5.03   & 5.80  \\
\hline
$L^{\pi}_{\phi}(f,T)$, $T=25000$:\\
Restarted                  & & 8.44  & 8.67  & 11.19  & &  8.35  & 8.58  & 11.18 & &  8.92  & 9.19   & 12.27  \\
Non-restarted              & & 31.27 & 31.31 & 32.54  & & 32.42  & 31.22 & 33.14 & & 19.10  & 19.42  & 23.56  \\
Fixed step size of $0.1$   & & 26.29 & 27.00 & 34.51  & & 26.57  & 27.30 & 35.12 & & 28.32  & 29.25  & 38.94  \\
Fixed step size of $0.01$  & & 7.89  & 8.10  & 10.59  & & 7.86   & 8.08  & 10.61 & & 8.41   & 8.66   & 11.59  \\
Fixed step size of $0.001$ & & 3.26  & 3.32  & 4.07   & &  2.84  & 2.92  & 3.69  & & 2.96   & 3.04   & 3.91  \\
\hline\vspace{-0.7cm}
\end{tabular}
\caption{\small \textbf{Performance of restarted EGS under noisy cost observations.\vspace{-0.6cm}}}
\label{tab:EGS}
\end{center}
\end{table}

\begin{figure}[!ht]
\centering
\includegraphics[height=2.0in]{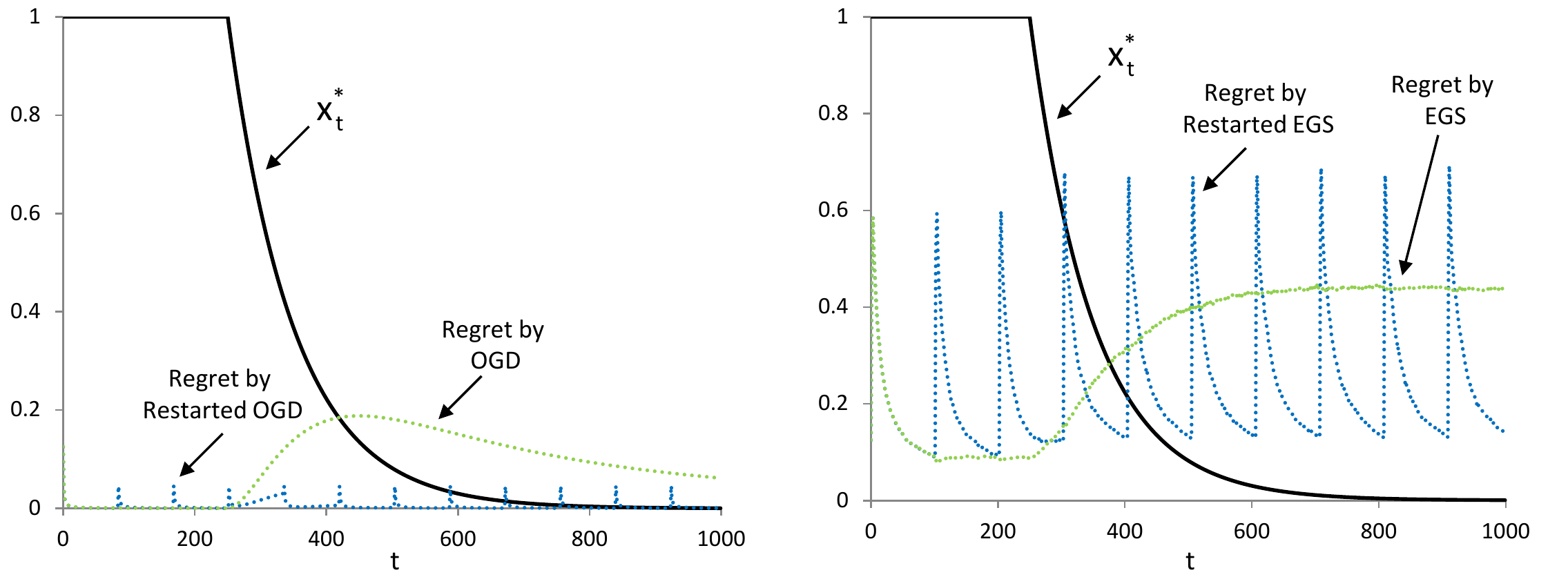}
\vspace{-0.3cm}\caption{\small \textbf{Regret in the presence of changing cost.} The cost $f_{t}(x_{t}) =~\frac{x_{t}^{2}}{2} -~b^{decay}_{t}x_{t} +~1$ is reflected by $x^{\ast}_{t} = b^{decay}_{t}$, with $T = 1,000$, $\tau = 250$, and $\sigma = 0.3$. \emph{(Left)} The average regret incurred at each epoch by the restarting procedure with OGD as a subroutine, and the one incurred at each epoch by OGD (without restarting), under feedback $\phi_{t}^{(1)}$. \emph{(Right)} The average regret incurred at each epoch by the restarting procedure with EGS as a subroutine, and the one incurred in each epoch by EGS (without restarting), under feedback $\phi_{t}^{(0)}$.
} \label{fig:simulation}\vspace{-0.2cm}
\end{figure}

Figure \ref{fig:simulation} depicts the averaged regret the restarted policies incur at each epoch, for one representative instance (decay-type variation, $T=1000$, $\sigma=0.3$). For illustration purposes, Figure \ref{fig:simulation} also includes the regret incurred at each epoch by the subroutine policies (OGD/EGS) when those applied without restarting.
The estimated values of $\alpha$  (capturing the regret rates) are consistent with the theoretical bounds ($0.5$ for restarted OGD under $\phi^{(1)}$, $0.67$ for restarted EGS under $\phi^{(0)}$). Together with the estimation of the coefficient $c$ (ranges in $[0.05,0.94]$ for restarted OGD, in $[2.09,2.88]$ for restarted EGS), these demonstrate the actual performance of the policies in a variety of cost-varying instances. One may observe that when $\sigma$ is larger (observations are more noisy) the multiplying constant typically increases. The estimated loss values indicate the extant at which each policy's performance is ``close" the one of the dynamic oracle (the restarted policies as well as the subroutine themselves, when applied without restarting, gets ``closer" to the dynamic oracle when $T$ grows).

One may observe that, not surprisingly, the restarted policies consistently outperform the OCO policies when these are not restarted. Considering policies with a fixed step size, we observe that different setting are characterized by different step sizes are considered to be the (ex-post) ``best" in different settings; the ``right" step size is effected by the variation pattern, the feedback structure, and the noisiness of the observations.
Indeed, comparing the performance of the restarting procedure to that of policies with fixed step size, one may observe that in various settings the restarting procedure is outperformed by a certain fixed step size; this occur more often for small values of $T$.
While there are various heuristics to set a-priory a fixed step size, non of those have any performance guarantee relative to the dynamic oracle for arbitrary variation (even when the variation is known to be fixed). We note that policies with a fixed step-size may perform well (and even better than known rate-optimal policies) over finite horizons even when the environment is stationary.\footnote{Repeating the numerical analysis for various \emph{stationary} settings we observed that policies with fixed step sizes may incur relative loss of less than 0.02 percent under noisy gradient access, and less than 1 percent under noisy cost access; for various (small enough) values of $T$ these policies outperformed rate optimal SA policies.} While the restarting policies were not designed and tuned in this paper to optimize practical performance, in most of the instances that are considered here they perform at least ``on par" with policies with the considered fixed step-sizes.


\end{document}